\documentclass[twoside,11pt]{article}
\usepackage{versions}  

\usepackage{jmlr2e}
\usepackage[utf8]{inputenc} 
\usepackage[T1]{fontenc}    
\usepackage{hyperref}       
\usepackage{url}            
\usepackage{booktabs}       
\usepackage{amsfonts}       
\usepackage{lipsum}
\usepackage{amsmath}
\usepackage{microtype}
\usepackage{graphicx}
\usepackage{subfigure}
\usepackage{booktabs} 
\usepackage{cite}
\usepackage{comment}
\usepackage{hyperref}
\usepackage{tikz}
\usepackage{kotex}
\usepackage{amssymb} 
\usepackage{algorithm,algpseudocode}
\usepackage{hyperref}

\newtheorem{claim}{Claim}
\newtheorem{fact}{Fact}

\DeclareMathOperator*{\argmin}{arg\,min}
\DeclareMathOperator*{\minimize}{minimize}
\DeclareMathOperator*{\maximize}{maximize}

\DeclareMathOperator*{\tr}{tr}
\DeclareMathOperator*{\diag}{diag}
\newcommand{\norm}[1]{\left\lVert#1\right\rVert}
\newcommand{\alp}{\alpha}
\newcommand{\bet}{\beta}
\newcommand{\gam}{\gamma}
\newcommand{\bE}{\mathbb{E}}
\newcommand{\bP}{\mathbb{P}}

\newcommand{\bR}{\mathbb{R}}

\newcommand{\cA}{\mathcal{A}}

\newcommand{\cC}{\mathcal{C}}

\newcommand{\cF}{\mathcal{F}}

\newcommand{\cI}{\mathcal{I}}

\newcommand{\cL}{\mathcal{L}}

\newcommand{\cO}{\mathcal{O}}
\newcommand{\cP}{\mathcal{P}}

\newcommand{\cR}{\mathcal{R}}

\newcommand{\fb}{\mathbf{b}}
\newcommand{\fc}{\mathbf{c}}
\newcommand{\fd}{\mathbf{d}}

\newcommand{\fq}{\mathbf{q}}

\newcommand{\fu}{\mathbf{u}}
\newcommand{\fv}{\mathbf{v}}
\newcommand{\fw}{\mathbf{w}}

\newcommand{\falp}{\pmb{\alpha}}
\newcommand{\fbet}{\pmb{\beta}}
\newcommand{\flam}{\pmb{\lambda}}
\newcommand{\fgam}{\pmb{\gamma}}
\newcommand{\fA}{\mathbf{A}}
\newcommand{\fB}{\mathbf{B}}
\newcommand{\fC}{\mathbf{C}}
\newcommand{\fD}{\mathbf{D}}
\newcommand{\fF}{\mathbf{F}}
\newcommand{\fG}{\mathbf{G}}
\newcommand{\fK}{\mathbf{K}}

\newcommand{\kA}{\mathfrak{A}}

\newcommand{\pg}{\pmb{g}}
\newcommand{\pf}{\pmb{f}}
\newcommand{\px}{\pmb{x}}

\usepackage[framemethod=TikZ]{mdframed}
\mdfsetup{frametitlealignment=\center}

\hypersetup{ hidelinks }
\newcommand{\revision}[1]{{#1}}

\jmlrheading{25}{2024}{1-69}{10/21; Revised
3/23}{2/24}{21-1256}{Chanwoo Park and Ernest Ryu}

\ShortHeadings{Optimal First-Order Algorithms as a Function of Inequalities}{Park and Ryu}
\firstpageno{1}

\begin{document}

\title{Optimal First-Order Algorithms as a Function of Inequalities}

\author{\name Chanwoo Park \email cpark97@mit.edu \\
       \addr Department of Electrical Engineering and Computer Science \\
       Massachusetts Institute of Technology\\
       Massachusetts, United States of America\\
       \AND
       \name Ernest K.\ Ryu \email eryu@snu.ac.kr \\
       \addr Department of Mathematical Sciences \\
       Interdisciplinary Program in Artificial Intelligence\\
       Seoul National University\\
       Seoul, Korea}
\editor{Martin Jaggi}

\maketitle

\begin{abstract}
In this work, we present a novel algorithm design methodology that finds the optimal algorithm as a function of inequalities. Specifically, we restrict convergence analyses of algorithms to use a prespecified subset of inequalities, rather than utilizing all true inequalities, and find the optimal algorithm subject to this restriction. This methodology allows us to design algorithms with certain desired characteristics. As concrete demonstrations of this methodology, we find new state-of-the-art accelerated first-order gradient methods using randomized coordinate updates and backtracking line searches.
\end{abstract}


\section{Introduction}
\label{Introduction}
Nesterov's seminal work presented the fast gradient method (FGM) with rate $\mathcal{O}(1/k^2)$ \citep{10029946121}, and Nemirovsky and Yudin established a complexity lower bound matching the rate up to a constant \citep{nemirovsky1983problem}. A rich line of research following these footsteps flourished in the following decades, and FGM became the prototypical ``optimal'' method. Recently, however, it was discovered that the Nesterov's FGM can be improved by a constant; the optimized gradient method (OGM) \citep{drori2014performance, kim2016optimized} outperforms FGM by a factor of $2$. Furthermore, the prior complexity lower bound was also improved by a constant factor to exactly match the rate of OGM \citep{drori2017exact}. Thus, the search for the exact optimal first-order gradient is now complete, and OGM, not FGM, emerges as the victor.

That FGM can be improved was, in our view, a surprising discovery, and it leads us to ask the following questions.
First, can we also improve the variants of Nesterov's FGM in related setups? FGM's acceleration has been extended to utilize randomized coordinate updates \citep{nesterov2012efficiency, allen2016even, nesterov2017efficiency} and backtracking line searches \citep{beck2009fast}.
Second, is there some sense in which Nesterov's FGM is exactly optimal?


In this work, we address these two questions by examining the inequalities used in the analyses of the algorithms.
We introduce the notion of $\mathcal{A}^\star$-optimality, which defines optimality of an algorithm \emph{conditioned} on a set of inequalities.
OGM is the $\mathcal{A}^\star$-optimal algorithm conditioned on all true inequalities and therefore is the exact optimal algorithm in the classical sense. However, other algorithms become $\mathcal{A}^\star$-optimal when conditioned on a different restrictive subset of inequalities.
By restricting convergence analyses of algorithms to use only a prespecified subset of inequalities with good properties, rather than utilizing all true inequalities, we obtain algorithms with better capacity for extensions.
Specifically, we obtain new algorithms utilizing randomized coordinate updates and backtracking linesearches that improve upon the prior state-of-the-art rates.
Moreover, we show that FGM is the optimal \revision{algorithm} roughly in the sense that it is the best algorithm that admits the use of randomized coordinate updates and backtracking line searches.

\paragraph{Contributions.}
As the main contribution of this work, we present an algorithm design methodology based on $\mathcal{A}^\star$-optimality and the performance estimation problem (PEP)  \citep{drori2014performance,taylor2017smooth} and demonstrate the strength of the methodology by finding new $\mathcal{A}^\star$-optimal algorithms that improve the state-of-the-art rates achieved by variants of FGM.
As a minor contribution, we establish the optimality of FGM in the following sense: FGM is the $\mathcal{A}^\star$-optimal algorithm that relies on a certain set of inequalities that are amenable to both randomized coordinate updates and backtracking line searches.

\subsection{Preliminaries and notations}
\label{ss:prelim}
In this section, we review standard definitions and set up the notation. 

\paragraph{Problem setting and $L$-smoothness.} 
For $L>0$, $f\colon\mathbb{R}^n\rightarrow\mathbb{R}$ is $L$-smooth if $f$ is differentiable and
\[
\|\nabla f(x)-\nabla f(y)\| \le L\|x-y\|\qquad
\forall\,x,y\in \mathbb{R}^n.
\]
Throughout this paper, we consider the problem 
\begin{align*}
\begin{array}{ll}
\underset{x\in \mathbb{R}^n}{\mbox{minimize}}  &f(x)
  \end{array}
\end{align*}
with the following assumptions
\begin{itemize}
\item [\hypertarget{A1}{(A1)}] $f$ is convex
\item [\hypertarget{A2}{(A2)}] $f$ is $L$-smooth
\item [\hypertarget{A3}{(A3)}] $f$ has a minimizer (not necessarily unique)
\item [\hypertarget{A4}{(A4)}] $\underset{x\in \mathbb{R}^n}{\inf}  f(x)>-\infty$.
\end{itemize}
We write $x_\star$ to denote a minimizer of $f$ if one exists and $f_\star= \underset{x\in \mathbb{R}^n}{\inf}  f(x)$ for the optimal value.

\paragraph{Fixed-step first-order algorithms (FSFO).}
The class of fixed-step first-order algorithms (FSFO) has the following form: given a differentiable $f$, total iteration count $N$, and starting point $x_0 \in \bR^d$, the iterates are defined by 
\begin{align}
    x_{i+1} = x_i - \frac{1}{L} \sum_{k=0}^i h_{i+1, k} \nabla f(x_k) \label{eqn::FO}
\end{align}
for $i = 0,1, \dots, \revision{N-1}$. The coefficients $\{h_{i,k}\}_{0 \leq k < i \leq N}$ may depend on $N$ and prior information about the function $f$, such as the smoothness coefficient $L$, but are otherwise predetermined.
In particular, $\{h_{i,k}\}_{0 \leq k < i \leq N}$ may not depend on function values or gradients observed throughout the algorithm.
The classical algorithms such as FGM, OGM, and the heavy-ball method are all FSFO.

\paragraph{Nesterov's fast gradient method (FGM).}
The celebrated FGM is 
\begin{align*}
  y_{k+1} &= x_{k} - \frac{1}{L}\nabla{f(x_{k})}
\\
  z_{k+1} &= z_{k} - \frac{\theta_{k}}{L}\nabla{f(x_{k})} 
\\
  x_{k+1} &= \left(1- \frac{1}{\theta_{k+1}}\right)y_{k+1} + \frac{1}{\theta_{k+1}}z_{k+1},
\end{align*}
where $z_0 = x_0$, $\theta_0 = 1$, and $\{\theta_k\}_{k=0}^\infty$ is defined as $\theta_{k+1} = \frac{1 + \sqrt{4\theta_k^2 + 1}}{2}$ for $k=0,1,\dots$
 \citep{10029946121}. FGM has the rate 
 \[
 f(y_N) - f_\star \leq \frac{2L\norm{x_0 -x_\star}^2}{N^2} + o \left(\frac{1}{N^2}\right),
 \]
 which is optimal \emph{up to a constant}.
 Many extensions and variants of FGM have been presented, including versions utilizing randomized coordinate updates \citep{lee2013efficient, allen2016even, nesterov2017efficiency} and backtracking linesearches \citep{beck2009fast}.

\paragraph{Optimized gradient method (OGM).}
Let $N$ be the total iteration count.
OGM \citep{drori2014performance,kim2016optimized} is 
\begin{align*}
  y_{k+1} &= x_{k} - \frac{1}{L}\nabla{f(x_{k})}
\\
  z_{k+1} &= z_{k} - \frac{2\theta_k}{L}\nabla{f(x_{k})} 
\\
  x_{k+1} &= \left(1-\frac{1}{\theta_{k+1}}\right)y_{k+1} + \frac{1}{\theta_{k+1}}z_{k+1},
\end{align*}
for $k = 0,1, \dots, N-2$, where $z_0=y_0 =x_0$ and $\theta_k$ is the same as with FGM.
Different from FGM, OGM has what we refer to as the \emph{last-step modification}
\begin{align*}
  x_{N} &= \left(1-\frac{1}{\tilde{\theta}_{N}}\right)y_{N} + \frac{1}{\tilde{\theta}_{N}}z_{N},
\end{align*}
where $\{\tilde{\theta}_k\}_{k=1}^\infty$ is defined as $\tilde{\theta}_k = \frac{1 + \sqrt{8\theta_{k-1}^2 +1}}{2}$. OGM exhibits the rate 
$$f(x_N) - f_\star \leq \frac{L\norm{x_0 -x_\star}^2}{N^2} + o \left(\frac{1}{N^2}\right),$$
which is faster than FGM by a factor of $2$ and is in fact exactly optimal \citep{drori2017exact}.
This remarkable discovery was made using a computer-assisted methodology, the performance estimation problem (PEP) \citep{drori2014performance,kim2016optimized}. Variants of OGM with randomized coordinate updates or backtracking linesearches had not been discovered.

\paragraph{Optimized gradient method - Gradient norm (OGM-G).}
Let $N$ be the total iteration count.
The method OGM-G has what we refer to as the \emph{first-step modification}
\begin{align*}
  y_{1} &= x_{0} - \frac{1}{L}\nabla{f(x_{0})}
\\
  z_{1} &= z_{0} - \frac{1+\tilde{\theta}_{N}}{2L}\nabla{f(x_{0})} 
\\
  x_{1} &= \frac{\theta_{N-1}^4}{\tilde{\theta}_{N}^4}y_{1} + \left(1-\frac{\theta_{N-1}^4}{\tilde{\theta}_{N}^4}\right)z_{1},
\end{align*}
where $z_0 = y_0 = x_0$ and $\{\theta_k\}_{k=0}^\infty$, $\{\tilde{\theta}_k\}_{k=1}^\infty$ are defined as OGM.
The remaining iterates of OGM-G are defined as
\begin{align*}
  y_{k+1} &= x_{k} - \frac{1}{L}\nabla{f(x_{k})}
\\
  z_{k+1} &= z_{k} - \frac{\theta_{N-k+1}}{L}\nabla{f(x_{k})} 
\\
  x_{k+1} &= \frac{\theta_{N-k-1}^4}{\theta_{N-k}^4}y_{k+1} + \left(1-\frac{\theta_{N-k-1}^4}{\theta_{N-k}^4}\right)z_{k+1},
\end{align*}
for $k =1, 2, \dots, N-1$ \citep{kim2021optimizing,lee2021geometric}.
Note that the indices of the $\theta$-coefficients are decreasing as the iteration count increases.
Different to FGM and OGM, the guarantee of OGM-G is on the gradient magnitude: 
\[
\norm{\nabla f(x_N)}^2 \leq \mathcal{O}\left(\frac{L\left(f(x_0) - f_\star\right)}{N^2}\right).
\]
An important use of OGM-G is that when combined with FGM (or OGM), one can achieve the rate \citep[Remark 2.1]{nesterov2020primal}
$$\norm{\nabla f(x_N)}^2 \leq \mathcal{O}\left(\frac{L^2\|x_0-x_\star\|^2}{N^4}\right).$$
Recently, analyses of OGM-G based on potential function approaches have been presented \citep{lee2021geometric, diakonikolas2021potential}. 

\paragraph{Coordinate-wise smoothness.}
We say $f\colon\mathbb{R}^n\rightarrow\mathbb{R}$ is coordinate-wise smooth with parameters $(L_1, \dots, L_n)$ if it is differentiable and
\begin{align}
  \| \nabla_i f(x + \delta e_i) - \nabla_i f(x)\| \leq L_i \delta \label{eqn::coordinate-wise-coco}  
\end{align}
for all $x \in \bR^n$, $\delta>0$, and $i=1,\dots,n$, 
where $\nabla_if=\frac{\partial f}{\partial x_i}e_i$ is the $i$\nobreakdash-th partial derivative \revision{vector} and $e_i$ is the $i$\nobreakdash-th unit vector for $i=1,\dots,n$.

\paragraph{Fast gradient method - randomized coordinate updates (FGM-RC).} 
There are several randomized coordinate update variants of Nesterov's FGM \citep{nesterov2012efficiency, allen2016even, nesterov2017efficiency}. We discuss the version of \citet{allen2016even}, which we call FGM-RC, as it has the smallest (best) constant.
Assume $f\colon\mathbb{R}^n\rightarrow\mathbb{R}$ is convex and coordinate-wise smooth with parameters $(L_1, L_2, \dots, L_n)$. FGM-RC is 
\begin{align*}
    &\text{Sample } i(k) \text{ from }\{1,2, \dots, n \} \text{ with } \bP (i(k) = t) = \frac{\sqrt{L_t}}{S}
    \\
    &y_{k+1} = x_k - \frac{1}{L_{i(k)}} \nabla_{i(k)} \nabla f(x_{k})
    \\
    &z_{k+1} = z_k - \frac{k+2}{2S^2}\frac{1}{p_{i(k)}}\nabla_{i(k)} f(x_{k})
    \\
    &x_{k+1} = \frac{k+1}{k+3}y_{k+1} + \frac{2}{k+3}z_{k+1}
\end{align*}
for $k = 0,1, \dots$, where $z_0 = y_0 = x_0$ , \revision{$p_{t} = \mathbb{P}(i = t)$}, and $S = \sum_{k=1}^n \sqrt{L_{k}}$.
 FGM-RC exhibits the rate
$$f(y_N) - f_\star \leq \frac{2S^2\norm{x_0 -x_\star}^2}{(N+1)^2}.$$

While Allen-Zhu et al.'s FGM-RC uses convenient rational coefficients, their algorithm can be slightly sharpened (through straightforward modifications of their presented analysis) to use the $\theta$\nobreakdash-coefficients of Nesterov.
We refer to this refinement as FGM-RC$^\sharp$:
\begin{align*}
    &\text{Sample } i(k) \text{ from }\{1,2, \dots, n \} \text{ with } \bP (i(k) = t) = \frac{\sqrt{L_t}}{S}
    \\
    &y_{k+1} = x_k - \frac{1}{L_{i(k)}} \nabla_{i(k)} f(x_{k})
    \\
    &z_{k+1} = z_k - \frac{\theta_k}{S^2}\frac{1}{p_{i(k)}}\nabla_{i(k)} f(x_{k})
    \\
    &x_{k+1} = \left(1 - \frac{1}{\theta_{k+1}}\right)y_{k+1} + \frac{1}{\theta_{k+1}}z_{k+1}
\end{align*}
for $k = 0,1,\dots$, where $z_0 = y_0 = x_0$ and $S = \sum_{k=1}^n \sqrt{L_{k}}$. FGM-RC$^\sharp$ exhibits the rate
$$f(y_N) - f_\star \leq \frac{S^2\norm{x_0 -x_\star}^2}{2\theta_{N-1}^2}.$$

Compared to first-order methods utilizing the full gradient, randomized coordinate updates methods have a lower cost per iteration and can be significantly faster.
In particular, FGM-RC$^\sharp$ can be significantly faster than FGM or OGM.

\paragraph{Fast gradient method - backtracking linesearch (FGM-BL).} \citet{beck2009fast} provides a version of Nesterov's FGM that uses backtracking linesearches (FGM-BL). Define $z_0 = y_0 = x_0$, $\eta>1$, and $L_0>0$.
Consider the backtracking linesearch that finds the smallest nonnegative integer $i_k$ such that with $\bar{L} = \eta^{i_k} L_{k-1}$ 
$$f\left(x_k - \frac{1}{\bar{L}} \revision{\nabla f(x_k)}\right) \leq f(x_k) - \frac{1}{2\bar{L}} \norm{\nabla f(x_k)}^2$$
holds, for each step $k$.
In FGM-BL, we set $L_k = \eta^{i_k} L_{k-1}$ and define
\begin{align*}
    y_{k+1} &= x_k - \frac{1}{L_k}\nabla f(x_k)
    \\
    z_{k+1} &= z_k - \frac{\theta_k}{L_k} \nabla f(x_k)
    \\
    x_{k+1} &= \left( 1- \frac{1}{\theta_{k+1}} \right) y_{k+1} + \frac{1}{\theta_{k+1}}  z_{k+1}
\end{align*}
for $k = 0,1, \dots$. FGM-BL exhibits the rate 
$$f(y_{N}) -f_\star \leq \frac{L_{N}}{2\theta_{N-1}^2}\norm{x_0 - x_\star}^2. $$

The backtracking linesearch is useful when we do not know the smoothness parameter $L$. FGM-BL obtains an estimate of $L$ while making progress with the accelerated gradient method.

\paragraph{Computer-assisted algorithm design.}
The performance estimation problem (PEP) is a computer-assisted proof methodology that analyzes the worst-case performance of optimization algorithms through semidefinite programs \citep{drori2014performance, taylor2017smooth, taylor2017exact}. 
The use of the PEP has \revision{led} to many discoveries that would have otherwise been difficult without the assistance \citep{kim2018another,taylor2018exact, taylor2019stochastic,barre2020principled,de2020worst,gu2020tight,lieder2020convergence, ryu2020operator,Dragomir2021, kim2021accelerated, yoon2021accelerated}. 
Notably, the algorithms OGM \citep{drori2014performance,kim2016optimized,kim2018generalizing}, OGM-G \citep{kim2021optimizing}, and ITEM \citep{taylor2021optimal} were obtained by using the PEP for the setup of minimizing a smooth convex (possibly strongly convex) function.
OGM and ITEM improve the rates of Nesterov's FGM by constants and have an exact matching complexity lower bound \citep{drori2017exact, drori2021oracle}.
The integral quadratic constraints (IQC) is another technique based on control-theoretic notions for computer-assisted algorithm design \citep{lessard2016analysis, hu2017dissipativity,van2017fastest,fazlyab2018analysis, seidman2019control, zhang2021unified}. This work builds upon the PEP methodology.

\subsection{Organization}
This paper is organized as follows. 
Section~\ref{sec::handy-ineq} defines the notion of handy inequalities and discusses how it will be utilized for designing accelerated algorithms with randomized coordinate updates and backtracking linesearches.
Section~\ref{sec::best-algorithm} defines the notion of $\cA^\star$-optimality.
Section~\ref{ssec::ORC-F} exhibits the main methodology of this work by using it to obtain an  $\cA^\star$-optimal algorithm and a variant utilizing randomized coordinate updates.
Section~\ref{sec::other-results} presents several other $\cA^\star$-optimal algorithms obtained with our algorithm design methodology.
The proofs of $\cA^\star$-optimality of the algorithms of Section~\ref{sec::other-results} are deferred to Sections~\ref{sec::appendix-deferred-cal} and \ref{sec::proof} of the appendix.

As the proofs of $\cA^\star$-optimality require lengthy calculations, we provide Matlab scripts verifying them. Specifically, the following scripts show that the derived analytical results agree with the numerical solutions of the SDPs:\\
 \url{https://github.com/chanwoo-park-official/A-star-map/}. 


\section{Handy inequalities for deriving variants of FSOM}
\label{sec::handy-ineq}
In this section, we define the notion of \emph{handy inequalities}.
The definition captures the empirical observation that the use of some inequalities makes the algorithm more amenable to modifications (and are therefore ``handy'') while other inequalities make the analysis brittle and not amenable to modifications.
In particular, FGM has extensions using randomized coordinate updates and backtracking linesearches, as discussed in the preliminaries, while such modifications seem difficult with OGM.
By identifying the notion of handy inequalities, we point out that the fault is in the inequalities being utilized.




\subsection{Inequalities for smooth convex functions}
In this section, we quickly review and name a few commonly used inequalities for smooth convex functions.

If $f\colon\mathbb{R}^n\rightarrow\mathbb{R}$ is convex, $L$-smoothness is equivalent to 
\begin{align}
    f(x)\ge f(y)+\left\langle \nabla f(y),x-y\right\rangle+\frac{1}{2L}\|\nabla f(x)-\nabla f(y)\|^2  \label{eqn::coco}
\end{align}
for all $x,y\in\mathbb{R}^n$.
We call \eqref{eqn::coco} the \textit{cocoercivity inequality on $(x,y)$}.
As a particular case, when $y = x- \frac{1}{L} \nabla f(x)$, the cocoercivity inequality on $(x,y)$ becomes
\begin{align*}
    f(x)\ge f(y)+ \frac{1}{2L}\norm{\nabla f(x)}^2 +\frac{1}{2L}\|\nabla f(y)\|^2,
\end{align*}
and dropping the last term leads to
\[
f(x) \geq f(y) + \frac{1}{2L} \norm{\nabla f(x)}^2.
\]
We call this the \textit{gradient-step inequality at $x$}.
Dropping $\frac{1}{2L}\norm{\nabla f(x) - \nabla f(y)}^2$ in \eqref{eqn::coco}, we get
\[
f(x) \geq f(y) + \langle \nabla f(y), x-y \rangle.
\]
We call this the \textit{convexity inequality on $(x,y)$}.
Note that the gradient-step or convexity inequalities are weaker than the cocoercivity inequality in the sense that they were obtained by dropping a nonnegative term from the cocercivity inequality.

If $f\colon\mathbb{R}^n\rightarrow\mathbb{R}$ is convex and coordinate-wise smooth with parameters $(L_1, \dots, L_n)$, then 
\begin{equation}
    f(x) \geq f(y) + \langle \nabla f(y), x-y \rangle + \frac{1}{2L_i} \norm{\nabla_i f(x) - \nabla_i f(y)}^2
      \label{eqn::coord-coco}
    \end{equation}
for all $x,y\in\mathbb{R}^n$ and $i=1,\dots,n$.
We call this the \textit{coordinate-wise cocoercivity inequality on $(x,y,i)$}.
As a particular case, when $y = x- \frac{1}{L_i} \nabla_i f(x)$, the coordinate-wise cocoercivity inequality on $(x,y,i)$ becomes
\begin{align*}
    f(x)\ge f(y)+ \frac{1}{2L_i}\norm{\nabla_i f(x)}^2 +\frac{1}{2L_i}\|\nabla_i f(y)\|^2.
\end{align*}
Dropping the last term, we get
\[
f(x) \geq f(y) + \frac{1}{2L_i}\norm{\nabla_i f(x)}^2.
\]
We call this the \textit{coordinate-wise gradient-step inequality at $(x,i)$}.

\paragraph{Proof of \eqref{eqn::coord-coco}.}
While we suspect the coordinate-wise cocoercivity inequality to be commonly known, we are unaware of a written proof to reference. We therefore quickly provide the following proof.

Let $g$ be convex, coordinate-wise smooth with parameters $(L_1, \dots, L_n)$, and $y = x- \frac{1}{L_i} \nabla_i g(x)$. Then, we have
\begin{align*}
    g_\star- g(x) &\leq g(y) - g(x) 
    \\
    &= \int_{t=0}^1 \langle \nabla g(x + t(y-x)), y-x \rangle dt 
    \\
    &= \langle \nabla g(x), y-x \rangle + \int_{t=0}^1 \langle \nabla g(x + t(y-x)) - \nabla g(x), y-x \rangle dt 
    \\
    &= \langle \nabla g(x), y-x \rangle + \int_{t=0}^1 (\nabla_i g(x + t(y_i-x_i)e_i) - \nabla_i g(x))(y_i-x_i) dt
    \\
    &\leq\langle \nabla g(x), y-x \rangle + \int_{t=0}^1 tL_i (y_i-x_i)^2 dt = -\frac{1}{2L_i}\norm{\nabla_i g(x)}^2.
\end{align*}
\revision{For all $y$, the function $f(x) - f(y) - \langle \nabla f(y), x -y \rangle$ is convex and coordinate-wise smooth with parameters $(L_1, \dots, L_n)$, as a function of $x$. Therefore, for all $y$, setting  $g(x) = f(x) - f(y) - \langle \nabla f(y), x -y \rangle$, we conclude}
\begin{align*}
    f(x) - f(y) - \langle \nabla f(y), x-y \rangle \geq \frac{1}{2L_i} \norm{\nabla_i f(x) - \nabla_i f(y)}^2.
\end{align*}

\subsection{The inequalities of FGM and OGM}
The analyses of FGM and OGM crucially differ in the inequalities they use.
The common convergence analysis of FGM defines the Lyapunov function 
\[
U_k=\theta_{k-1}^2(f(y_k) - f_\star) + \frac{L}{2}\norm{z_k - x_\star}^2
\]
and establishes the non-increasing property
\begin{align*}
    U_{k} - U_{k+1} &= \theta_k^2 \left(f(x_k) - f(y_{k+1}) - \frac{1}{2L} \norm{\nabla f(x_k)}^2\right) 
    \\
    &+\theta_{k-1}^2\left(f(y_k) - f(x_{k}) - \langle \nabla f(x_k), y_k - x_k \rangle \right)
    \\
    & + \theta_k\left(f_\star - f(x_k) - \langle \nabla f(x_k), x_\star - x_k \rangle \right)
    \\ &\geq 0.
\end{align*}
In contrast, the analysis of OGM defines the Lyapunov function
\begin{align*}
    U_k = &2\theta_{k}^2\left(f(x_k) - f_\star - \frac{1}{2L}\norm{\nabla f(x_k)}^2 \right) + \frac{L}{2}\norm{z_{k+1}- x_\star}^2
\end{align*}
and establishes the non-increasing property
\begin{align*}
U&_{k} - U_{k+1} 
\\
& = 2\theta_{k}^2\left(f(x_{k}) - f(x_{k+1}) +\langle\nabla f(x_{k+1}), x_{k+1} - x_k \rangle - \frac{1}{2L}\norm{\nabla f(x_{k}) - \nabla f(x_{k+1})}^2\right) 
\\
& + 2\theta_{k+1} \left(f_\star - f(x_{k+1})  + \langle \nabla f(x_{k+1}), x_{k+1} - x_\star \rangle- \frac{1}{2L}\norm{\nabla f(x_{k+1})}^2
\right)\\
&\ge 0.
\end{align*}

Note the difference in the inequalities being used;
FGM uses gradient-step and convexity inequalities, while OGM uses cocoercivity inequalities. 
The fact that OGM uses the stronger cocoercivity inequalities partially explains why its guarantee is stronger than FGM's guarantee.
These inequalities are all, of course, true inequalities, but the inequalities used by FGM are much more amenable to obtaining variants using randomized coordinate updates and backtracking line searches; the inequalities used by OGM do not have this property.
In Sections~\ref{sssec::CI} and \ref{sssec::LI}, we define inequalities that admit such variants to be algorithmically \emph{handy}.

\subsection{Handy inequalities for randomized coordinate updates}

\label{sssec::CI}


We now discuss the notion of handy inequalities for randomized coordinate updates.
Let us examine the analysis of FGM-RC$^\sharp$.
For $k = 0,1, \dots $, define $$U_k = \frac{\theta_{k-1}^2}{S^2}(f(y_k) - f_\star) + \frac{1}{2}\norm{z_k - x_\star}^2$$ and write $\bE_{i(k)}$ for the expectation conditioned on information up to the $k$\nobreakdash-th iteration.
Then,
\begin{align*}
    U_k - U_{k+1} & = 
    \frac{\theta_k^2}{S^2}\left(f(x_k) - f(y_{k+1}) - \frac{1}{2L_{i(k)}}\norm {\nabla_{i(k)} f(x_k)}^2 \right)
    \\
    &+  \frac{\theta_{k-1}^2}{S^2} \left(f(y_k) - f(x_{k}) - \frac{S}{\sqrt{L_{i(k)}}}\langle \nabla_{i(k)} f (x_k) , y_k -x_k \rangle\right)
    \\
    &+\frac{\theta_{k} }{S^2}\left(f_\star - f(x_{k}) - \frac{S}{\sqrt{L_{i(k)}}}\langle \nabla_{i(k)} f(x_{k}), x_\star- x_k \rangle \right),
\end{align*}
and taking the conditional expectation $\bE_{i(k)}$ gives us
\begin{align*}
    U_k - \bE_{i(k)} U_{k+1} &=\bE_{i(k)} \left[
    \frac{\theta_k^2}{S^2}\left(f(x_k) - f(y_{k+1}) - \frac{1}{2L_{i(k)}}\norm {\nabla_{i(k)} f(x_k)}^2 \right)
    \right]
    \\
    &+  \frac{\theta_{k-1}^2}{S^2} \left(f(y_k) - f(x_{k}) - \langle \nabla f(x_k) , y_k -x_k \rangle\right)
    \\
    &+\frac{\theta_{k} }{S^2}\left(f_\star - f(x_{k}) - \langle \nabla f(x_{k}), x_\star- x_k \rangle \right) 
    \\
    &\geq 0.
\end{align*}
Finally, taking the full expectation gives us
\[
\frac{\theta_k^2}{S^2}(\bE [f(y_{k+1})] -f_\star) \leq \bE U_{k+1} \leq\dots \leq U_0 \leq \frac{1}{2}\norm{x_0 - x_\star}^2.
\]

We can interpret this convergence analysis as a direct modification of FGM's analysis by taking expectations of the inequalities.
Utilizing the linearity of expectation to obtain the convexity inequality and having the coordinate-wise gradient-step inequality holds almost surely is crucial.
The gradient-step inequality and the convexity inequality are handy for randomized coordinate updates as this analysis of FGM-RC$^\sharp$ demonstrates.
The coordinate-wise cocoercivity inequality on $(x_\star,x_k)$ is also handy as the term $\frac{1}{2L_{i(k)}}\norm {\nabla_{i(k)} f(x_k)}^2$ is one we can take the expectation of.

On the other hand, the cocoercivity inequalities on $(y_k,x_k)$ or $(x_k,y_{k+1})$ do not seem to be handy as the terms $\|\nabla_{i(k)}f(x_k)-\nabla f(y_k)\|^2$ or $\|\nabla_{i(k)} f(x_k)-\nabla f(y_{k+1})\|^2$ are not easily manipulated under expectations.
For this reason, adapting OGM and its analysis to use randomized coordinate updates seems difficult.

\subsection{Handy inequalities for backtracking linesearch}
\label{sssec::LI}


Next, we discuss the notion of handy inequalities for backtracking linesearches.
Let us examine the analysis of FGM-BL.
For $k = 0,1,\dots $, define $$U_{k,L} = \frac{\theta_{k-1}^2}{L}(f(y_{k}) - f_\star) + \frac{1}{2}\norm{z_{k} - x_\star}^2.$$
Then, 
\begin{align*}
    U_{k, L_{k+1}} - U_{k+1,L_{k+1}} &= \frac{1}{L_{k+1}} \Biggl(\theta_k^2 \left(f(x_k) - f(y_{k+1}) - \frac{1}{2L_{k+1}} \norm{\nabla f(x_k)}^2\right) 
    \\
    &+\theta_{k-1}^2\left(f(y_k) - f(x_{k}) - \langle \nabla f(x_k), y_k - x_k \rangle \right)
    \\
    & + \theta_k\left(f_\star - f(x_k) - \langle \nabla f(x_k), x_\star - x_k \rangle \right)\Biggr)
    \\ &\geq 0.
\end{align*}
Note that the inequality $f(x_k) - f(y_{k+1}) - \frac{1}{2L_{k+1}} \norm{\nabla f(x_k)}^2\ge 0$ is enforced by the backtracking linesearch.
Finally, we conclude
\[
\frac{\theta_{k}^2}{L_{k+1}}(f(y_{k+1}) - f_\star) \leq U_{k+1,L_{k+1}}\leq U_{k, L_{k+1}} \leq U_{k, L_{k}} \leq \dots \leq U_{0, L_0} \leq \frac{1}{2}\norm{x_0 - x_\star}^2.
\]

We can interpret this convergence analysis as a direct modification of FGM's analysis with $L$ replaced with $L_{k+1}$, an estimate of the unknown Lipschitz parameter $L$.
The role of the backtracking linesearch is to verify the inequality involving $L_k$.

For a linesearch to be implementable, it is critical that it relies on quantities that are \emph{algorithmically observable}.
In the analysis of OGM, the inequalities
\begin{gather*}
f(x_{k}) - f(x_{k+1}) +\langle\nabla f(x_{k+1}), x_{k+1} - x_k \rangle - \frac{1}{2L}\norm{\nabla f(x_{k}) - \nabla f(x_{k+1})}^2\ge 0\\
f_\star - f(x_{k+1}) + \langle \nabla f(x_{k+1}), x_{k+1} - x_\star \rangle- \frac{1}{2L}\norm{\nabla f(x_{k+1})}^2 \ge 0
\end{gather*}
are used.
The first inequality involves algorithmically observable quantities and is therefore handy for backtracking linesearches.
However, the second is not handy as its verification requires the knowledge of $x_\star$.
The convexity inequality on $(x_\star,x_{k+1})$ is handy as it does not involve $L$ and hence does not require verification through a linesearch.

\revision{
To clarify, we define the notion of ``handy inequalities'' informally through examples. The motivation is to avoid having certain problematic terms in the analysis. In the following sections, we demonstrate cases where we succeed in generating algorithms with the desired characteristic using the notion of handy inequalities. }

\section{Optimal algorithm map}
\label{sec::best-algorithm}
In this section, we define the notion of $\cA^\star$-optimality, the notion of optimality conditioned on a set of inequalities.

\subsection{$\cA^\star$-optimality}
Define oracles $\cO=(\cO_0,\cO_1)$ to take as input a function and a point and return zero and \revision{first-order} information of the function, i.e., $\cO_0(f,x)=f(x)$ and $\cO_1(f,x)=\nabla f(x)$. 
Define the optimal oracle $\cO^\star=(\cO_x^\star,\cO_f^\star) $, which takes as input a function and returns an optimal point, if one exists, and the optimal value,
i.e., $\cO_x^\star(f)=x^\star\in\argmin f$ and $\cO_f^\star(f)=f^\star = \inf f$.
The optimal oracle $\cO^\star$ is used in the minimax formulation, but is, of course, not used in the algorithms.

Let $\kA_N$ be the class of fixed-step first-order algorithms (FSFO) with $N$ iterations, i.e., an algorithm in $\kA_N$ may access $\cO_1$ up to $N$ times.
To further specify our notation, a \revision{first-order} algorithm $\cA_N(x_0, f)\colon  \bR^k \times \cF_L \to \bR^{k \times N}$ in $\kA_N$ generates the $N$ iterates as follows: 
\begin{align*}
    x_1 &= \cA_{N,1}(x_0, \cO_1(f, x_0)) 
    \\
    x_2 &= \cA_{N,2}(x_0, \cO_1(f, x_0), \cO_1(f, x_1))
    \\
    &\vdots
    \\
    x_N &= \cA_{N,N}(x_0, \cO_1(f, x_0), \dots, \cO_1(f, x_{N-1})),
\end{align*}
\revision{where $\cA_{N, i}$ is defined for $i \in \{1,2, \dots, N\}$ as 
\begin{align*}
    \cA_{N,i}(x_0, g_0, \dots, g_{i-1}) = x_0 - h_{i, 0} g_0 - \dots - h_{i,i-1} g_{i-1}.
\end{align*}
}
Let $\cP$ be a \emph{performance criterion} that measures the performance of an algorithm $\cA$ on a function $f$. \revision{To clarify, a performance criterion only depends on $f(x^\star), \{f(x_i), \nabla f(x_i)\}_{i=0}^{N}, x^\star, \{x_i\}_{i=0}^N$. For example,  
the function-value suboptimality} 
\[
\cP(\cA_N(x_0, f), \cO, \cO^\star) =  f(x_N) - f_\star= \cO_0(f, x_N) - \cO^\star_f(f)  
\]
or the squared gradient magnitude 
\[
\cP(\cA_N(x_0, f),  \cO, \cO^\star)
= \norm{\nabla{f(x_N)}}^2
= \norm{ \cO_1(f, x_N)}^2 
\]
of the last iterate, $x_N$ are commonly used performance criteria.




Let $\cC$ be an \emph{initial condition}, a condition we impose or assume on the initial point $x_0$. \revision{To clarify, an initial condition only depends on $f(x^\star), (f(x_0), \nabla f(x_0)), x^\star, x_0$. For example,  
the initial distance to a solution }
\[
\cC(x_0, \cO, \cO^\star) = \{\|x_0-x_\star\|\le R\}
 = \{\|x_0-\cO^\star_x(f)\|\le R\}
\]
or function value suboptimality 
\[
\cC (x_0, \cO, \cO^\star)= \{ f(x_0)-f_\star\le R \}
= \{ \cO_0(f, x_0)-\cO^\star_f(f)\le R \}
\]
are commonly used \revision{initial conditions}.

Let $\cI$ be an \emph{inequality collection}, a set of inequalities the output of the oracles\\ $\mathcal{O}(f, x_0),\dots,\mathcal{O}(f, x_{N}),\mathcal{O}^\star(f)$ we assume satisfies.
In prior work, convergence analyses were permitted to use all true inequalities. 
Unique to our work, we consider analyses based on a restricted inequality collection $\cI$;
convergence proofs may use inequalities in $\cI$, a strict subset of the true inequalities.

Define the \emph{rate} (or risk) of an algorithm conditioned on $\mathcal{I}$ as
\begin{alignat*}{3}
    \cR(\cA_N, \cP, \cC, \cI) = &\sup_{x_0, f,\cO, \cO^\star} &&\cP(\cA_N(x_0, f ),  \cO,\cO^\star)
    \\
    &\text{subject to } && x_{i} = \cA_{N,i} (x_0,  g_0, \dots, g_{i-1}), \quad  i \in \{1,2,\dots, N \}
    \\
    & &&(x_0, g_0, f_0, x_\star, f_\star) \revision{\text{ satisfies } \cC(x_0, \cO, \cO^\star)} \\
    & && \{ (x_i, g_i, f_i)\}_{i=0}^N \text{  and  }  (x_\star, f_\star) \text{  satisfy } \cI
    \\
    & &&(f_{i}, g_i) = \cO(f, x_{i}), \quad  i \in \{0, \dots, N\}
    \\
    & &&(x_\star, f_\star) = \cO^\star(f).
\end{alignat*}
Note that we impose constraints on $f$ only through the output of the oracles $ \cO$ and $\cO^\star$.
Define the minimax optimal rate conditioned on $\mathcal{I}$ as
\begin{equation}
\begin{aligned}
    \cR^\star(\kA_N, \cP, \cC, \cI) &= \inf_{\cA_N\in \kA_N} \cR(\cA_N, \cP, \cC, \cI).
\end{aligned}\label{eqn::a*map}
\end{equation}
If this infimum is attained, write $\cA^\star_N$ to denote the optimal algorithm, and (with some abuse of notation) say the algorithm is $\cA^\star$-optimal conditioned on $\mathcal{I}$.
Conversely, write $\cA^\star_N(\cP, \cC, \cI)$ to denote the $\cA^\star$-optimal algorithm, and refer to this as the $\cA^\star$-map.
(An $\cA^\star$-optimal algorithm may or may not be \revision{unique}.)


\subsection{Optimality of OGM}
A series of work on OGM \citep{drori2014performance, kim2016optimized, drori2017exact} established that OGM is the exact optimal first-order gradient method.
Using our notation, we can express these prior results as
$$\operatorname{OGM} = \cA^\star_N\left(f(x_{N})- f_\star, \norm{x_0 - x_\star} \leq R, \cI_{L-\operatorname{smooth}}\right)$$
and
\begin{equation*}
\begin{aligned}
    \cR^*\left(\kA_N,f(x_{N})- f_\star, \norm{x_0 - x_\star} \leq R, \cI_{L-\operatorname{smooth}}\right) 
    =\frac{LR^2}{2\tilde{\theta}_N^2},
\end{aligned}\label{eqn::rate-ogm}
\end{equation*}
where
\begin{align*}
  \cI_{L-\operatorname{smooth}} = &\left\{f(x_i) \geq f(x_j) \revision{-} \langle \nabla f(x_j), x_j - x_i \rangle + \frac{1}{2L} \norm{\nabla f(x_i) - \nabla f(x_j)}^2\right\}_{i,j=0}^{\revision{N}} 
  \\
  &\qquad \qquad \bigcup \left\{ f_\star \geq f(x_k) \revision{-} \langle \nabla f(x_k), x_k -x_\star \rangle + \frac{1}{2L}\norm{\nabla f(x_k)}^2 \right\}_{k=0}^N. 
\end{align*}

\subsection{Inequality collection selection}
When considering $\cA^\star$-optimal algorithms, the choice of the inequality collection represents a tradeoff.
{ On one extreme, if the inequality collection is empty, the performance criterion is the function-value suboptimality or the square gradient magnitude, and the initial condition is the initial distance condition or the function value suboptimality condition (i.e.,
$\cI=\emptyset$, $\cP(\cA_N, \cO, \cO^\star) = \cO_0(f,x_n) - \cO_f^\star(f) \text{ or }\norm{\cO_1(f,x_N)}^2$, and $\cC(x_0, \cO, \cO^\star) = \{\norm{x_0 - \cO_x^\star(f)} \leq R\} \text{ or } \{\cO_0(f, x_0) - \cO_f^\star(f) \leq R \}$), no convergence analysis can be done, and the ``algorithm'' that does not move from the starting point is $\mathcal{A}^\star$-optimal.}
On the other extreme, using all true inequalities in the smooth convex minimization setup makes OGM $\mathcal{A}^\star$-optimal.
The inequality collections that we consider in later sections include handy inequalities that have the capacity to admit randomized coordinate updates or backtracking linesearches while being sufficiently powerful to establish good rates.\footnote{\revision{
As a relevant negative result, we tried but did not succeed in finding a randomized coordinate update version of OGM-G. We tried to modify co-coercivity inequality to the randomized coordinate version. Proof of the OGM-G uses co-coercivity inequality on $(x_k, x_{k+1})$. Since we choose a random direction for each iterate, it is hard to utilize in $\norm{\nabla f(x_k) - \nabla f(x_{k+1})}^2$ term for different direction partial differentiation.  
We considered several inequality collections that are handy for randomized coordinate updates, but the $\cA^\star$-optimal algorithms conditioned on those inequality collections exhibited $\cO(1/k)$ rates, i.e., the handy inequalities we considered were not sufficiently powerful to establish the accelerated $\cO(1/k^2)$ rate of OGM-G.
}}



\section{ORC-F (Optimized randomized coordinate updates - function value)}
\label{ssec::ORC-F}
In this section, we present ORC-F$_\flat$, an $\cA^\star$-optimal algorithm.
We first state the theorem precisely describing the $\cA^\star$-optimality result, while deferring the proof to the end of this section.
We then provide a direct Lyapunov analysis of ORC-F$_\flat$ and modify this Lyapunov analysis to obtain ORC-F, a randomized coordinate update version \revision{of} ORC-F$_\flat$.

\subsection{Main results}
Optimized randomized coordinate updates - function value$_\flat$ (\textbf{ORC-F$_\flat$}) is defined as
\begin{align*}
    y_{k+1} &= x_k -\frac{1}{L}\nabla f(x_k)
    \\
    z_{k+1} &= z_{k} - \frac{\varphi_{k+1} - \varphi_k}{L}\nabla{f(x_{k})}
    \\
    x_{k+1} &=\frac{\varphi_{k+1}}{\varphi_{k+2}} y_{k+1} +  \left(1 - \frac{\varphi_{k+1}}{\varphi_{k+2}}\right) z_{k+1} 
\end{align*}
for $k = 0, 1, \dots$ where $y_0 = z_0 = x_0$, \revision{$ \varphi_0 = 0$, and the strictly increasing sequence} $\{\varphi_k\}_{k=0}^\infty$ is defined by $\left(2\varphi_{k+1} - \varphi_k\right) = (\varphi_{k+1} - \varphi_k)^2 $  for $k =0, 1 \dots$. 
\begin{theorem}[$\mathcal{A}^\star$-optimality of ORC-F$_\flat$]
\label{thm::ORC-F-optimal}
ORC-F$_\flat$ is $\cA^\star$-optimal in the sense that
$$\operatorname{ORC-F_\flat} = \cA^\star_N\left(f(y_{N+1})- f_\star, \norm{x_0 - x_\star} \leq R, \cI_{\operatorname{ORC-F_\flat}}\right)$$
and has the minimax optimal rate
$$\cR^\star\left(\kA_N,f(y_{N+1}) - f_\star, \norm{x_0 - x_\star} \leq R, \cI_{\operatorname{ORC-F_\flat}}\right) = \frac{LR^2}{2\varphi_{N+1}}
$$
with respect to the inequalities
\begin{align*}
  \cI_{\operatorname{ORC-F_\flat}} = &\biggl\{  f(x_k) \geq f(y_{k+1}) + \frac{1}{2L} \norm{\nabla f(x_k)}^2 \biggr\}_{k=0}^N  
  \\
  &\bigcup \biggl\{f(y_k) \geq f(x_k) + \langle \nabla f(x_k), y_{k} - x_{k} \rangle\biggr\}_{k=1}^{N}
  \\
  &\bigcup \biggl\{f_\star \geq f(x_k) + \langle \nabla f(x_k), x_\star - x_k \rangle  + \frac{1}{2L} \norm{\nabla f(x_{k})}^2\biggr\}_{k=0}^N.
\end{align*}
\end{theorem}
Note that the inequalities in $\cI_{\operatorname{ORC-F_\flat}}$ are handy for randomized coordinate updates. We defer the proof of Theorem~\ref{thm::ORC-F-optimal} to Section~\ref{ssec::proof-thm1}.

The following corollary is a consequence of Theorem~\ref{thm::ORC-F-optimal}, but we state it separately and present a standalone proof so that we can modify it for the proof of Theorem~\ref{thm::ORC-F-optimal2}.

\begin{corollary}
\label{cor::ORCFb}
Assume (\hyperlink{A1}{A1}), (\hyperlink{A2}{A2}), and (\hyperlink{A3}{A3}). ORC-F$_\flat$'s $y_k$-sequence exhibits the rate
\begin{align*}
f(y_{k+1}) -f_\star &\leq \frac{L\norm{x_0 - x_\star}^2}{2\varphi_{k+1}}
\end{align*}
for $k=1,2,\dots$.
\end{corollary}

\begin{proof}
For $k = 0, 1, 2, \dots$, define 
$$U_k = {\varphi_{k}}(f(y_k) -f_\star) + \frac{L}{2}\norm{z_k - x_\star}^2. $$
Then we have 
\begin{align*}
    U_{k} - U_{k+1} &= \varphi_k (f(y_k) - f_\star) + \frac{L}{2}\norm{z_k - x_\star}^2 -  \varphi_{k+1} (f(y_{k+1}) - f_\star) - \frac{L}{2}\norm{z_{k+1} - x_\star}^2
    \\
    &= \varphi_k (f(y_k) - f_\star)-\varphi_{k+1} (f(y_{k+1}) - f_\star)  + \frac{L}{2}\langle z_k - z_{k+1}, z_k + z_{k+1} - 2x_\star \rangle 
    \\
    &=\varphi_k (f(y_k) - f_\star)-\varphi_{k+1} (f(y_{k+1}) - f_\star) 
    \\
    &+\frac{L}{2} \left\langle \frac{\varphi_{k+1} - \varphi_k}{L} \nabla f(x_k), 2z_k -  \frac{\varphi_{k+1} - \varphi_k}{L} \nabla f(x_k) - 2x_\star \right\rangle
    \\
    &= \varphi_{k+1}\left(f(x_k) - f(y_{k+1}) - \frac{1}{2L} \norm{\nabla f(x_k)}^2\right) 
    \\
    &+\varphi_{k}\left(f(y_k) - f(x_{k}) - \langle \nabla f(x_k), y_k - x_k \rangle \right)
    \\
    & + (\varphi_{k+1} - \varphi_{k})\left(f_\star - f(x_k) - \langle \nabla f(x_k), x_\star - x_k \rangle - \frac{1}{2L} \norm{\nabla f(x_k)}^2 \right)
    \\
    &+\frac{2\varphi_{k+1} - \varphi_k }{2L} \norm{\nabla f(x_k)}^2 + \varphi_k \langle \nabla f(x_k), y_k - x_k \rangle + (\varphi_{k+1} - \varphi_k) \langle \nabla f(x_k), x_\star - x_k \rangle
    \\
    & - \frac{(\varphi_{k+1} - \varphi_k)^2}{2L} \norm{\nabla f(x_k)}^2 + (\varphi_{k+1} - \varphi_k) \langle \nabla f(x_k), z_k - x_\star \rangle.
\end{align*}
Since 
\begin{align*}
    (\varphi_{k+1} - \varphi_k)\langle \nabla f(x_k), z_k - x_\star \rangle &= (\varphi_{k+1} - \varphi_k)\langle \nabla f(x_k), z_k - x_k \rangle + (\varphi_{k+1} - \varphi_k)\langle \nabla f(x_k), x_k - x_\star \rangle  
    \\
    &=  \varphi_k\langle \nabla f(x_k), x_k - y_k \rangle +(\varphi_{k+1} - \varphi_k)\langle \nabla f(x_k), x_k - x_\star \rangle  
\end{align*}
and $(2\varphi_{k+1} - \varphi_k) = (\varphi_{k+1} - \varphi_k)^2$, we get 
\begin{align*}
    U_{k} - U_{k+1} &= \varphi_{k+1}\left(f(x_k) - f(y_{k+1}) - \frac{1}{2L} \norm{\nabla f(x_k)}^2\right) 
    \\
    &+\varphi_{k}\left(f(y_k) - f(x_{k}) - \langle \nabla f(x_k), y_k - x_k \rangle \right)
    \\
    & + (\varphi_{k+1} - \varphi_{k})\left(f_\star - f(x_k) - \langle \nabla f(x_k), x_\star - x_k \rangle - \frac{1}{2L} \norm{\nabla f(x_k)}^2 \right)
    \\ &\geq 0.
\end{align*}
We conclude $\revision{\varphi_{k+1}(f(y_{k+1}) - f^\star)} \leq U_{k+1} \leq \dots \leq U_{-1} = \frac{L}{2}\norm{x_0 - x_\star}^2$.
\end{proof}
Note that the proof only utilized inequalities in $\cI_{\operatorname{ORC-F_\flat}}$.

\paragraph{Randomized coordinate updates version.}
Assume $f$ is \revision{a} coordinate-wise smooth function with parameters $(L_1, \dots, L_n)$. Define $S  = \sum_{i=1}^n {\sqrt{L_i}}$.
At iteration $k$, select the coordinate $i(k)$ with probability $ \bP(i(k) = t ) = \frac{\sqrt{L_{t}}}{S}$.
Define optimized randomized coordinate updates - function value (\textbf{ORC-F}), a randomized coordinate updates version of ORC-F$_\flat$, as
\begin{align*}
    y_{k+1} &= x_k -\frac{1}{L_{i(k)}}\nabla_{i(k)} f(x_k)
    \\
    z_{k+1} &= z_{k} - \frac{\varphi_{k+1} - \varphi_k}{S\sqrt{L_{i(k)}}}\nabla_{i(k)}{f(x_{k})}
    \\
    \revision{x_{k+1}} &= \revision{\frac{\varphi_{k+1}}{\varphi_{k+2}} y_{k+1} + \left(1- \frac{\varphi_{k+1}}{\varphi_{k+2}}\right) z_{k+1}}
\end{align*}
for $k=0,1,\dots$. 
\begin{theorem}
\label{thm::ORC-F-optimal2}
Assume (\hyperlink{A1}{A1}) and (\hyperlink{A3}{A3}).
    Assume $f$ is a coordinate-wise smooth function with parameters $(L_1, \dots, L_n)$. Then ORC-F exhibits the rate as 
    $$\bE \left[ f(y_{k+1}) \right] - f_\star \leq \frac{S^2 \norm{x_0 - x_\star}^2}{2\varphi_{k+1}}$$
    for $k = 0,1,\dots$. 
\end{theorem}
\begin{proof}
For $k = 0, 1, \dots$, define 
$$U_k = \frac{\varphi_k}{S^2}(f(y_k) - f_\star) + \frac{1}{2}\norm{z_k - x_\star}^2$$ and define $\bE_{i(k)}$ as the expectation conditioned on $i(0),\dots,i(k-1)$. Then, we have 
\begin{align*}
    U_k - U_{k+1} &=  \frac{\varphi_k}{S^2}(f(y_k) - f_\star) + \frac{1}{2}\norm{z_k - x_\star}^2 -  \frac{\varphi_{k+1}}{S^2}(f(y_{k+1}) - f_\star) - \frac{1}{2}\norm{z_{k+1} - x_\star}^2
    \\
    &= \frac{\varphi_k}{S^2}(f(y_k) - f_\star)  - \frac{\varphi_{k+1}}{S^2}(f(y_{k+1}) - f_\star) + \frac{1}{2}\langle z_k - z_{k+1}, z_k + z_{k+1} - 2x_\star \rangle
    \\
    &=\frac{\varphi_k}{S^2}(f(y_k) - f_\star)  - \frac{\varphi_{k+1}}{S^2}(f(y_{k+1}) - f_\star)
    \\&+ \frac{1}{2}\left\langle  \frac{\varphi_{k+1} - \varphi_k}{S \sqrt{L_{i(k)}}}\nabla_{i(k)} f(x_k) , 2z_k -\frac{\varphi_{k+1} - \varphi_k}{S \sqrt{L_{i(k)}}}\nabla_{i(k)} f(x_k)- 2x_\star \right\rangle
    \\
    &=\frac{(\varphi_{k+1} - \varphi_k)}{S^2} \left( f_\star - f(x_k) - \frac{S}{\sqrt{L_{i(k)}}} \langle \nabla_{i(k)} f(x_k), x_\star - x_k \rangle - \frac{1}{2L_{i(k)}} \norm{\nabla_{i(k)} f(x_k)}^2\right)
    \\
    &+ \frac{\varphi_k}{S^2}\left(f(y_k) - f(x_k) - \frac{S}{\sqrt{L_{i(k)}}} \langle \nabla_{i(k)} f(x_k), y_k - x_k \rangle\right)
    \\
    &+\frac{\varphi_{k+1}}{S^2} \left(f(x_k) -f(y_{k+1}) - \frac{1}{2L_{i(k)}} \norm{\nabla_{i(k)} f(x_k)}^2\right) 
    \\
    &+ \frac{2\varphi_{k+1} - \varphi_k}{2S^2L_{i(k)}}\norm{\nabla_{i(k)} f(x_k)}^2 + \frac{\varphi_k}{S\sqrt{L_{i(k)}}}\langle \nabla_{i(k)} f(x_k), y_k - x_k \rangle+ \frac{\varphi_{k+1} - \varphi_k}{S\sqrt{L_{i(k)}}}\langle \nabla_{i(k)} f(x_k), x_\star - x_k \rangle  
    \\
    &- \frac{(\varphi_{k+1} - \varphi_k)^2}{2S^2 L_{i(k)}} \norm{\nabla_{i(k)}f(x_k)}^2 + 
     \frac{\varphi_{k+1} - \varphi_k}{S \sqrt{L_{i(k)}}}\langle \nabla_{i(k)} f(x_k), z_k - x_\star \rangle.
\end{align*}
Since 
\begin{align*}
    \frac{\varphi_{k+1} - \varphi_k}{S \sqrt{L_{i(k)}}}\langle \nabla_{i(k)} f(x_k), z_k - x_\star \rangle &= \frac{\varphi_{k+1} - \varphi_k}{S \sqrt{L_{i(k)}}}\langle \nabla_{i(k)}  f(x_k), z_k - x_k \rangle +\frac{\varphi_{k+1} - \varphi_k}{S \sqrt{L_{i(k)}}}\langle \nabla_{i(k)}  f(x_k), x_k - x_\star \rangle  
    \\
    &= \frac{\varphi_{k+1} - \varphi_k}{S \sqrt{L_{i(k)}}}\langle \nabla_{i(k)} f(x_k), x_k - x_\star \rangle   + \frac{\varphi_k}{S \sqrt{L_{i(k)}}}\langle\nabla_{i(k)}  f(x_k), x_k - y_k \rangle 
\end{align*}
and $(2\varphi_{k+1} - \varphi_k) = (\varphi_{k+1} - \varphi_k)^2$, we get 
\begin{align*}
    U_{k} - U_{k+1} &=\frac{(\varphi_{k+1} - \varphi_k)}{S^2} \left( f_\star - f(x_k) - \frac{S}{\sqrt{L_{i(k)}}} \langle \nabla_{i(k)} f(x_k), x_\star - x_k \rangle - \frac{1}{2L_{i(k)}} \norm{\nabla_{i(k)} f(x_k)}^2\right)
    \\
    &+ \frac{\varphi_k}{S^2}\left(f(y_k) - f(x_k) - \frac{S}{\sqrt{L_{i(k)}}} \langle \nabla_{i(k)} f(x_k), y_k - x_k \rangle\right)
    \\
    &+\frac{\varphi_{k+1}}{S^2} \left(f(x_k) -f(y_{k+1}) - \frac{1}{2L_{i(k)}} \norm{\nabla_{i(k)} f(x_k)}^2\right).
\end{align*}
and taking the conditional expectation, we have
\begin{align*}
    U_k - \bE_{i(k)} U_{k+1} &=  \bE_{i(k)} \left[ \frac{(\varphi_{k+1} - \varphi_k)}{S^2} \left( f_\star - f(x_k) - \langle \nabla f(x_k), x_\star - x_k \rangle - \frac{1}{2L_{i(k)}} \norm{\nabla_{i(k)} f(x_k)}^2\right) \right]
    \\
    &+ \frac{\varphi_k}{S^2}\left(f(y_k) - f(x_k) - \langle \nabla f(x_k), y_k - x_k \rangle\right)
    \\
    &+\bE_{i(k)} \left[\frac{\varphi_{k+1}}{S^2} \left(f(x_k) -f(y_{k+1}) - \frac{1}{2L_{i(k)}} \norm{\nabla_{i(k)} f(x_k)}^2\right) \right]
    \\
    &\geq 0.
\end{align*}
Taking the full expectation, we have $\bE U_k \leq \dots \leq U_{0}$, and we conclude the statement of the theorem.
\end{proof}

\paragraph{Discussion.}
ORC-F has the smallest (best) constant among the randomized coordinate updates methods, to the best of our knowledge.
In particular, ORC-F's rate is slightly faster \revision{than} that of FGM-RC of \citet{allen2016even} or FGM-RC$^\sharp$ since $\theta_k^2 \leq \varphi_{k+1}$ for $k=0,1,\dots$, which follows from induction.
However, the improvement is small as the leading-term constant is the same, i.e., $\theta_k^2/ \varphi_{k+1}\rightarrow 1$ as $k\rightarrow\infty$.

\subsection{\revision{Brief review of \citet{taylor2017smooth}}}
\label{ss::notation}
\revision{This section is closely follows \citet{taylor2017smooth} with minor changes including different in notation and the addition of some variables.}
Consider the underlying space $\mathbb{R}^d$ with $d\ge N+2$.
The assumption that $d$ is sufficiently large is made to obtain dimension-independent results. See \citep[Section 3.3]{taylor2017smooth} for further discussion of this matter. Denote $y_{k+1} = x_k - \frac{1}{L}\nabla f(x_k)$ and $y_0= x_0$. Define $f_{i,0}= f(x_i)$, $f_{i,1} = f(y_i)$, $g_i = \nabla f(x_i)$, and
\begin{equation}
\begin{aligned}
\fG &= \left(
\begin{array}{ccccc}
    \norm{x_0 - x_\star}^2 & \langle g_0, x_{0} - x_\star \rangle & \langle g_1, x_{0} - x_\star \rangle & \dots &\langle g_N, x_0 - x_\star \rangle \\
    \langle g_0, x_0 - x_\star \rangle & \norm{g_0}^2 & \langle g_1, g_0 \rangle & \dots & \langle g_N, g_0 \rangle \\
    \vdots & \vdots & \vdots & \ddots & \vdots  \\
    \langle g_N, x_0 - x_\star \rangle & \langle g_N, g_0 \rangle & \langle g_N, g_1 \rangle  & \dots & \norm{g_N}^2\\
\end{array}
\right) , 
\\
\revision{\fF_0}&= \left(
\begin{array}{c}
    f_{0,0} - f_\star \\
    f_{1,0} - f_\star \\ 
    \vdots 
    \\
    \revision{f_{N+1,0} - f_\star} \\
\end{array}
\right) , \qquad \fF_1= \left(
\begin{array}{c}
    f_{0,1} - f_\star \\
    f_{1,1} - f_\star \\ 
    \vdots 
    \\
    f_{N+1,1} - f_\star \\
\end{array}
\right).
\end{aligned}
\end{equation}
Note that $\fG\succeq 0$, $\fF_0\succeq 0$, and $\fF_1\succeq 0$, i.e.,
$\fG$ is positive semidefinite and $\fF_0$ and $\fF_1$ are elementwise nonnegative.
Since $d\ge N+2$, 
given $\fG$ and $\fF_0\succeq 0$, we can take \revision{the Cholesky factorization%
\footnote{
\revision{Since $\fG$ is not strictly positive definite, the ``Cholesky factorization'' is not unique.
In fact, any factorization of the form $\fG=MM^\intercal$ suffices.
See \citep[Section 3]{taylor2017smooth} for further discussion on this matter.}
}of $\fG$ to recover the triplet $\{ (x_i, g_i, f_i) \}_{i=0}^N$. }
Define
\begin{align}
  \px_0 = e_1 \in \bR^{N+2}, \quad \pg_i = e_{i+2} \in \bR^{N+2}, \quad \revision{\pf_i = e_{i+1} \in \bR^{N+2}} \label{eqn::egf}  
\end{align}
for $i = 0, 1, \dots, N$, where $e_i$ are standard basis of $\bR^{N+2}$ or $\bR^{N+1}$.
\revision{We define FSFO \eqref{eqn::FO} with $(h_{i,j})$ and $(s_{i,j})$ as }
\begin{equation}
\begin{aligned}
  \px_{i+1} &= \px_i - \sum_{k=0}^i\frac{h_{i+1, k}}{L}\pg_k
  \\
  &= \px_0 - \sum_{k=0}^i \frac{s_{i+1, k}}{L}\pg_k
\end{aligned}
\label{eqn::fsfo2}
\end{equation}
for $i = 0,1, \dots, N-1$.
(Note that $\px_i$ appears in the first expression while $\px_0$ appears in the second.)
With this new notation, we can write
\begin{align*}
    &f_{i,j} - f_\star = \pf_i^\intercal \fF_j , \qquad \qquad \revision{i = 0,1, \dots, N+1}, \quad j = 0, 1
    \\
    &\langle g_i, g_j\rangle = \pg_i^\intercal \fG \pg_j, \qquad \qquad i,j = 0,1, \dots, N
    \\
    &\norm{x_i - x_\star}^2 = \px_i^\intercal \fG \px_i, \qquad \qquad i = 0, 1, \dots, N
    \\
    &\langle g_i, x_j - x_\star \rangle = \pg_i^\intercal \fG \px_j, \qquad \qquad i,j = 0,1, \dots, N.
\end{align*}
This allows us to express the optimization of the algorithm as an optimization problem with variables $\fG, \fF_0, \fF_1$.

Let $\cF_{L}$ be the class of $L$-smooth convex functions. Let $I$ be an index set and consider the set of triplets $S = \{(x_i,g_i,f_i)\}_{i \in I}$, where $x_i,g_i \in \bR^d$ and $f_i \in \bR$ for all $i \in I$. We say $S$ is  $\cF_{L}$-interpolable if and only if there exists a function $f \in \cF_{L}$ that $g_i \in \partial f(x_i)$ and $f(x_i) = f_i$ for all $i \in I$.
\begin{fact}\emph{\citep[Theorem 4]{taylor2017smooth}}
    $S$ is $\cF_{L}$-interpolable if and only if 
    $$f_i - f_j - \langle g_j, x_i - x_j \rangle \geq \frac{1}{2L} \norm{g_i - g_j}^2,\qquad \forall\,i,j\in I.$$
\end{fact}

\subsubsection{Strong duality for the PEP}
We propose \revision{a} general PEP form. The original sdp-PEP of \citet{taylor2017smooth} is
\begin{alignat*}{4}
&\maximize_{\fG, \fF_0} \quad  &&b^\intercal \fF_0  + \operatorname{Tr}(C\fG)
\\
&\mbox{subject to}\quad 
    &&0 \geq (\pf_j - \pf_i)^\intercal \fF_0 + \operatorname{Tr}(\fG ((\px_i - \px_j) \pg_j^\intercal + \revision{\frac{1}{2L} (\pg_i -\pg_j)(\pg_i - \pg_j)^\intercal} ))  \qquad &&i,j \in \{0,1, \dots, N\}
    \\
    & && 1 \geq \operatorname{Tr}(\fG \px_0 \px_0^\intercal) 
    \\
    & && 0 \preceq \fG,
\end{alignat*}
for $b \in \bR^{N+1}$ and $C$ is a nonnegative definite matrix. For further details, refer to \citep[Theorem 5]{taylor2017smooth}. This original sdp-PEP is induced by the $\cF_L$-interpolable condition. \revision{We extend this sdp-PEP to replace the constraints with relaxed inequalities. Our general sdp-PEP is 
\begin{alignat*}{4}
&\maximize_{\fG, \fF_0, \fF_1} \quad  &&b_0^\intercal \fF_0  + b_1^\intercal \fF_1 +  \operatorname{Tr}(C\fG)
\\
&\mbox{subject to}\quad 
    &&\text{conditions corresponding to inequality collection $\mathcal{I}$}
    \\
    & && 1 \geq \operatorname{Tr}(\fG \px_0 \px_0^\intercal) 
    \\
    & && 0 \preceq \fG,
\end{alignat*}
for $b_0, b_1 \in \bR^{N+1}$ and $C$ is a nonnegative definite matrix.
Specific instances of this general sdp-PEP are considered in subsequent section.
We call the convex-dual problem of general sdp-PEP as dual-sdp-PEP \citep{taylor2017smooth}. 
Strong duality holds between the primal and dual SDPs.
\begin{fact} 
\label{fact::2}
    Assume the stepsizes of \eqref{eqn::fsfo2} satisfy
    $ s_{k, k-1} \neq 0$  for $k=1,\dots,N$. In addition, inequality collection corresponds to the algorithms in Sections 4.3, 5.1, and 5.2. Then, the strong duality holds between general sdp-PEP and dual-sdp-PEP. 
\end{fact}
The formal proof Fact~\ref{fact::2}, which we omit for the sake of brevity, follows from the same reasoning as that of \citep[Theorem 5]{taylor2017smooth}.
}
\subsection{\revision{Proof of Theorem~\ref{thm::ORC-F-optimal}}}
\label{ssec::proof-thm1}
\revision{In this section, we prove Theorem~\ref{thm::ORC-F-optimal}, i.e., $\mathcal{A}^\star$-optimality of ORC-F$_\flat$, using the PEP machinery.
To verify the lengthy calculations, we provide Matlab scripts verifying the analytical solution of the SDP: \\\url{https://github.com/chanwoo-park-official/A-star-map/}. }
\\To obtain ORC-F as an $\mathcal{A}^\star$-optimal algorithm, set $f(y_{N+1}) - f_\star$ to be the performance measure and $\norm{x_0 - x_\star} \leq R$ to be the initial condition.
Since the constraints and the objective of the problem \revision{are} homogeneous, we assume $R = 1$ without loss of generality. For the argument of homogeneous, we refer \revision{to} \citep{drori2014performance, kim2016optimized,taylor2017smooth}. We use the set of inequalities that are handy for randomized coordinate updates:
\begin{align*}
  \cI_{\operatorname{ORC-F}_\flat} = &\biggl\{  f_{k,0} \geq f_{k+1,1} + \frac{1}{2L} \norm{g_k}^2 \biggr\}_{k=0}^N  \bigcup \biggl\{f_{k,1} \geq f_{k,0} + \langle g_k, y_{k} - x_{k} \rangle\biggr\}_{k=1}^{N} 
  \\
  &\qquad   \bigcup \biggl\{f_\star \geq f_{k,0} + \langle g_k, x_\star - x_k \rangle + \frac{1}{2L} \norm{g_k}^2 \biggr\}_{k=0}^N.
\end{align*}
For calculating $\cR(\cA_N, \cP, \cC, \cI_{\operatorname{ORC-F}_\flat}) $ with fixed $\cA_N$, define the PEP with $\cI_{\operatorname{ORC-F}_\flat}$ as
\revision{\begin{equation}
\begin{aligned}
&\cR(\cA_N, \cP, \cC, \cI_{\operatorname{ORC-F}_\flat}) 
\\
&= \left(
\begin{array}{lllll}
    &\maximize \quad & f_{N+1,1}&  - f_\star
    \\
    &\text{subject to} \quad &1 \quad \geq& \norm{x_0 - x_\star}^2 
    \\
    & &f_{k,0} \geq& f_{k+1,1} + \frac{1}{2L} \norm{g_k}^2, \quad &k\in \{0,1, \dots, N\}
    \\
    & &f_{k,1} \geq& f_{k,0} + \langle g_k, y_{k} - x_{k} \rangle, &k \in \{1, \dots, N\}
    \\
    & &f_\star \mkern11mu \geq& f_{k,0} + \langle g_k, x_\star - x_k \rangle +\frac{1}{2L}\norm{g_k}^2,  &k \in \{0,1, \dots, N\}
    \\
    & & x_k, y_k &\text{ are following the algorithm } \cA_N.
\end{array}\right)
\end{aligned}
\label{eqn::ORC-primal}
\end{equation}}
Using the notation of Section~\ref{ss::notation}, we reformulate \revision{the} problem of computing the risk $\cR(\cA_N, \cP, \cC, \cI_{\operatorname{ORC-F}_\flat})$ as the following SDP:

\begin{alignat*}{4}
&\maximize_{\fG, \fF_0, \fF_1} \quad  &&\pf_{N+1}^\intercal \fF_1
\\
&\mbox{subject to}\quad 
    &&1 \geq \px_0^\intercal \fG \px_0  
    \\
    & && 0 \geq \pf_{k+1}^\intercal\fF_1 - \pf_{k}^\intercal \fF_0 + \frac{1}{2L} \pg_k^\intercal \fG \pg_k, \qquad &&k \in \{0,1, \dots, N\}
    \\
    & && 0 \geq \pf_k^\intercal(\fF_0 - \fF_1) + \pg_k^\intercal \fG (\px_{k-1} - \px_k ) - \frac{1}{L} \pg_{k-1}^\intercal \fG \pg_k, \qquad &&k \in \{1,2, \dots, N\}
    \\
    & && 0 \geq \pf_k^\intercal \fF_0 - \pg_k^\intercal \fG \px_k + \frac{1}{2L}\pg_k^\intercal \fG \pg_k, \qquad && k \in \{0,1, \dots, N\}
    \\
    & && \revision{\fG \succcurlyeq 0, \fF_0 \geq 0, \fF_1 \geq 0.}
\end{alignat*}
For above transformation, $d\geq N+2 $ is used \citep{taylor2017smooth}. \revision{The Lagrangian} of the optimization problem becomes 
\begin{align*}
    \Lambda(\fF_0, \fF_1&, \fG, \flam, \fbet, \falp, \tau)
    \\
    &=  - \pmb{f}_{N+1}^\intercal \fF_1 + \tau (\px_0^\intercal \fG \px_0 -1 ) + \sum_{k=0}^{N} \alpha_{k}\left( \pf_{k+1}^\intercal\fF_1 - \pf_{k}^\intercal \fF_0 + \frac{1}{2L} \pg_k^\intercal \fG \pg_k
    \right)
    \\& + \sum_{k=1}^N \lambda_k \left(\pf_k^\intercal(\fF_0 - \fF_1) + \pg_k^\intercal \fG (\px_{k-1} - \px_k ) - \frac{1}{L} \pg_{k-1}^\intercal \fG \pg_k\right) 
    \\
    & +\sum_{k=0}^{N} \beta_k \left( \pf_k^\intercal \fF_0 - \pg_k^\intercal \fG \px_k + \frac{1}{2L} \pg_k^\intercal \fG \pg_k\right) 
\end{align*}
with dual variables $\flam = (\lambda_1, \dots, \lambda_N) \in \bR_{+}^{N}$, $\fbet = (\beta_0, \dots, \beta_N) \in \bR_{+}^{N+1}$, $\falp = (\alp_0, \dots, \alp_N) \in \bR_{+}^{N+1}$, and $\tau \geq 0$. 
Then the dual formulation of PEP problem is
\begin{equation}
\begin{alignedat}{2}
    &\maximize_{(\flam, \fbet, \falp, \tau) \geq \pmb{0}}\quad  &&{-\tau} 
    \\
     &\mbox{subject to} &&\pmb{0} = -\sum_{k=0}^{N} \alpha_k \pmb{f}_k + \sum_{k=1}^N \lambda_k \pmb{f}_k +\sum_{k=0}^{N} \beta_k \pmb{f}_k
     \\
     & &&\pmb{0} = -\pmb{f}_{N+1} - \sum_{k=1}^N \lambda_k \pmb{f}_k + \sum_{k=0}^{N} \alpha_k \pmb{f}_{k+1}
     \\
     & &&0 \preceq S(\flam, \fbet, \falp, \tau),
\end{alignedat} \label{eqn::ORC-F-dualPEP}    
\end{equation}
where $S$ is defined as 
\begin{align*}
    S(\flam, \fbet, \falp, \tau) &=  \tau \pmb{x}_0\pmb{x}_0^\intercal+ \sum_{k=0}^{N} (\alpha_k + \beta_k) \left( \frac{1}{2L} \pmb{g}_{k}\pmb{g}_{k}^\intercal\right) + \sum_{k=0}^{N}\frac{\beta_k}{2}\left(-\pmb{g}_k \pmb{x}_k^\intercal-\pmb{x}_k \pmb{g}_k^\intercal  \right)
    \\
    &+ \sum_{k=1}^N \frac{\lambda_k}{2} \left(\pmb{g}_k (\pmb{x}_{k-1} - \pmb{x}_{k})^\intercal +  (\pmb{x}_{k-1} - \pmb{x}_{k})\pmb{g}_k^\intercal - \frac{1}{L} \pmb{g}_{k-1}\pmb{g}_{k}^\intercal - \frac{1}{L} \pmb{g}_{k}\pmb{g}_{k-1}^\intercal\right).
\end{align*}
Using the strong duality result of Fact~\ref{fact::2} and a continuity argument that we justify at the end of this proof, we proceed with
\begin{align}
    \argmin_{h_{i,j}}\maximize_{\fG, \fF_0, \fF_1} \mkern7mu \pf_{N+1}^\intercal \fF_1  = \argmin_{h_{i,j}}\minimize_{(\flam, \fbet, \falp,\tau) \geq \pmb{0}}\mkern7mu {\tau} \label{eqn::ORC-F-strongdual}
\end{align}
\revision{i.e., it is sufficient to obtain $h_{i,j}$'s argmin value of \eqref{eqn::ORC-F-dualPEP}. We omitted \eqref{eqn::ORC-F-strongdual}'s constraints for ease of writing. } Note that \eqref{eqn::ORC-F-dualPEP} finds the optimal proof for the algorithm. 
Minimizing \eqref{eqn::ORC-F-dualPEP} with respect to $(h_{i,j})$ corresponds to optimizing the algorithm:
\begin{alignat}{3}
    &\minimize_{h_{i,j}}&&\minimize_{(\flam, \fbet, \falp,\tau) \geq \pmb{0}}\quad && {\tau}  \label{eqn::ORC-F-hij-wkdual}
    \\
    & && \text{subject to} \quad && \pmb{0}  = -\sum_{k=0}^{N} \alpha_k \pmb{f}_k + \sum_{k=1}^N \lambda_k \pmb{f}_k +\sum_{k=0}^{N} \beta_k \pmb{f}_k \label{eqn::ORC-F-const1}
     \\
     & && &&\pmb{0}  = -\pmb{f}_{N+1} - \sum_{k=1}^N \lambda_k \pmb{f}_k + \sum_{k=0}^{N} \alpha_k \pmb{f}_{k+1} \label{eqn::ORC-F-const2}
     \\
     & && &&0 \preceq S(\flam, \fbet, \falp, \tau). \label{eqn::ORC-F-Smatrix}
\end{alignat}
We note that $\pf_i$ is \revision{a} standard unit vector mentioned in \eqref{eqn::egf} (not a variable), we can write \eqref{eqn::ORC-F-const1} and \eqref{eqn::ORC-F-const2} as 
\begin{equation}
\begin{alignedat}{2}
 &&\left(\begin{array}{ll}
\beta_k = \alpha_k - \lambda_k  = \lambda_{k+1} - \lambda_k, \qquad & k \in \{1, \dots, N-1\}
    \\
    \beta_0 = \alpha_0 = \lambda_1
    \\
    \beta_N = \alpha_N - \lambda_N = 1 - \lambda_N. 
\end{array}\right) 
\\
 &&\left(\begin{array}{ll}
\alpha_N = 1 
    \\
    \alpha_k = \lambda_{k+1}, \quad  & k \in \{0,1, \dots, N-1\}
\end{array}\right) 
\end{alignedat}\label{eqn::ORC-F-constlinear}
\end{equation}

We consider \eqref{eqn::ORC-F-Smatrix} with \eqref{eqn::ORC-F-constlinear} and FSFO's $h_{i,j}$. To be specific, we substitute $\falp$ and $\fbet$ to $\flam$ in $S(\flam, \fbet, \falp, \tau)$. To show the dependency of $S$ to $(h_{i,j})$ since $\px_k$ are represented with $(h_{i,j})$, we will explicitly write $S$ as $S(\flam, \tau;(h_{i,j}))$. Then, we get
\begin{align*}
    S(\flam, \tau;(h_{i,j}))     &=\tau \pmb{x}_0\pmb{x}_0^\intercal +\frac{\lambda_1}{2L} \pg_0\pg_0^\intercal  - \sum_{k=1}^{N-1} \frac{2\lambda_k - \lambda_{k+1}}{2L}\pmb{g}_k\pmb{g}_k^\intercal  + \frac{2-2\lambda_N}{2L} \pmb{g}_{N}\pmb{g}_{N}^\intercal + \sum_{k=1}^N\frac{\lambda_k}{2L} (\pmb{g}_{k-1}-\pmb{g}_{k})(\pmb{g}_{k-1}-\pmb{g}_{k})^\intercal
    \\&+ \sum_{k=1}^{N-1} \sum_{t=0}^{k-1} \left( \frac{\lambda_k}{2} \frac{h_{k,t}}{L} + \frac{\lambda_{k+1} - \lambda_{k}}{2}\sum_{j=t+1}^{k}\frac{h_{j,t}}{L} \right) \left(\pmb{g}_k\pmb{g}_t^\intercal + \pmb{g}_t\pmb{g}_k^\intercal\right) 
    \\&+ \sum_{t=0}^{N-1} \left( \frac{\lambda_N}{2} \frac{h_{N,t}}{L} + \frac{1 - \lambda_{N}}{2}\sum_{j=t+1}^{N}\frac{h_{j,t}}{L} \right) \left(\pmb{g}_N\pmb{g}_t^\intercal + \pmb{g}_t\pmb{g}_N^\intercal\right)
    \\&- \sum_{k=1}^{N-1} \frac{\lambda_{k+1} - \lambda_{k}}{2} \left(\pmb{x}_0\pmb{g}_k^\intercal + \pmb{g}_k\pmb{x}_0^\intercal\right) - \frac{\lambda_1}{2}\left(\pmb{x}_0\pmb{g}_0^\intercal + \pmb{g}_0\pmb{x}_0^\intercal\right) - \frac{1 - \lambda_N}{2}\left(\pmb{x}_0\pmb{g}_N^\intercal + \pmb{g}_N\pmb{x}_0^\intercal\right).
\end{align*}

Using the fact that $\px_0, \pg_i, \pf_i$ are unit vectors, we can represent $S(\flam, \tau;(h_{i,j}))$ with $\fgam(\flam) = -L\fbet = -L(\lambda_1, \lambda_2 - \lambda_1, \dots, 1-\lambda_N) = (\hat{\fgam}(\flam), \gamma_{N}(\flam))
$ and  $\tau' = 2L\tau$ as 
\begin{align*}
    S(\flam, \tau';(h_{i,j}))  &= \frac{1}{L}\left(
\begin{array}{ccc}
    \frac{1}{2}\tau' & \frac{1}{2}\hat{\fgam}(\flam)^\intercal & \frac{1}{2}\gamma_{N}(\flam) \\
    \frac{1}{2}\hat{\fgam}(\flam) & Q(\flam; (h_{i,j})) & q(\flam; (h_{i,j})) 
    \\
    \frac{1}{2}\gamma_{N}(\flam) & q(\flam; (h_{i,j}))^\intercal & \frac{2-\lambda_N}{2} \\
\end{array}
\right) \succeq 0.
\end{align*}
\revision{Here, $Q$ and $\fq$ are defined as 
\begin{align*}
     Q(\flam ;(h_{i,j})) &=  \frac{\lambda_1}{2} \pg'_0 \pg_0^{'\intercal} + \sum_{k=1}^{N-1} \frac{\lambda_{k+1}-2\lambda_k}{2}\pmb{g}'_k\pmb{g}_k^{'\intercal}  + \sum_{k=1}^{N-1}\frac{\lambda_k}{2} (\pmb{g}'_{k-1}-\pmb{g}'_{k})(\pmb{g}'_{k-1}-\pmb{g}_{k})^{'\intercal} + \frac{\lambda_N}{2}\pg'_{N-1}\pg_{N-1}^{'\intercal}
    \\&+ \sum_{k=1}^{N-1} \sum_{t=0}^{k-1} \left( \frac{\lambda_k}{2} h_{k,t} + \frac{\lambda_{k+1} - \lambda_{k}}{2}\sum_{j=t+1}^{k}h_{j,t} \right) \left(\pmb{g}'_k\pmb{g}_t^{'\intercal} + \pmb{g}'_t\pmb{g}_k^{'\intercal}\right)
\end{align*}
and 
\begin{align*}
    \fq(\flam;(h_{i,j})) &= -\frac{\lambda_N}{2}\pg'_{N-1} + \sum_{t=0}^{N-1} \left( \frac{\lambda_N}{2} h_{N,t} + \frac{1 - \lambda_{N}}{2}\sum_{j=t+1}^{N}h_{j,t} \right) \pmb{g}'_t
    \\
    &= \sum_{t=0}^{N-2} \left( \frac{\lambda_N}{2} h_{N,t} + \frac{1 - \lambda_{N}}{2}\sum_{j=t+1}^{N}h_{j,t} \right) \pmb{g}'_t +  \left(\frac{1}{2}{h_{N,N-1}} -\frac{\lambda_N}{2}\right)\pmb{g}'_{N-1}
\end{align*}
where $\pg'_k = e_{k+1} \in \bR^{N+1}$.}
Note that \eqref{eqn::ORC-F-hij-wkdual} is equivalent to 
\begin{alignat}{3}
    &\minimize_{h_{i,j}}&&\minimize_{(\flam,\tau') \geq \pmb{0}}\quad && {\tau'}  \label{eqn::ORC-F-hij-wkdual-2}
    \\
    & && \text{subject to} \quad && \left(
\begin{array}{ccc}
    \frac{1}{2}\tau' & \frac{1}{2}\hat{\fgam}(\flam)^\intercal & \frac{1}{2}\gamma_{N}(\flam) \\
    \frac{1}{2}\hat{\fgam}(\flam) & Q(\flam; (h_{i,j})) & q(\flam; (h_{i,j})) 
    \\
    \frac{1}{2}\gamma_{N}(\flam) & q(\flam; (h_{i,j}))^\intercal & \frac{2-\lambda_N}{2} \\
\end{array}
\right) \succeq 0 \nonumber
\end{alignat}
and dividing this optimized value with $2L$ gives the optimized value of \eqref{eqn::ORC-F-hij-wkdual}. Using Schur complement \citep{golub1996matrix}, (we already know $0\leq\lambda_N \leq 1$ by \eqref{eqn::ORC-F-constlinear}) \eqref{eqn::ORC-F-hij-wkdual-2} can be converted to the problem as 
\begin{alignat}{3}
    &\minimize_{h_{i,j}}&&\minimize_{(\flam, \tau') \geq \pmb{0}}\quad && {\tau'}  \label{eqn::ORC-F-hij-final-problem}
    \\
    & && \text{subject to} \quad && \left(
\begin{array}{cc}
    Q - \frac{2\fq \fq^\intercal}{2-\lambda_N} & \frac{1}{2}( \hat{\fgam}(\flam) - \frac{2\fq \gam_N(\flam)}{2-\lambda_N}) \\
    \frac{1}{2}( \hat{\fgam}(\flam) - \frac{2\fq \gam_N(\flam)}{2-\lambda_N})^\intercal & \frac{1}{2} (\tau' - \frac{\gam_N(\flam)^2}{2-\lambda_N}) \\
\end{array}
\right) \succeq 0. \label{eqn::ORC-F-final-matrix}
\end{alignat}
So far we simplified SDP. We will have \revision{three steps}: finding variables that \revision{make} \eqref{eqn::ORC-F-final-matrix}'s left hand side zero, showing that the solution from \revision{the} first step satisfies Karuch-Kuhn-Tucker (KKT) condition, and finally showing that obtained algorithm is equivalent to ORC-F$_\flat$. 

\begin{claim}
\label{claim::1}
    There is a point that makes \eqref{eqn::ORC-F-final-matrix}'s left-hand side zero. 
\end{claim}
\begin{proof}
\revision{Defining the positive sequence $\{\varphi_k\}_{k=0}^\infty$ as 
\begin{align*}
    2\varphi_{k+1} - \varphi_k = (\varphi_{k+1} - \varphi_k)^2
\end{align*}
for $k = 0,1, \dots$, \revision{$ \varphi_0 = 0$, and $\{\varphi_k\}_{k=0}^\infty$ is a  strictly increasing sequence. } Defining $\{r_{k,t}\}_{k=1,2, \dots, N,  t=0, \dots, k-1}$ as $$r_{k,t} = \lambda_k h_{k,t}- \frac{\gamma_k}{L}\sum_{j=t+1}^k h_{j,t}.$$ 
Then, if $r_{i,j}$ is determined, \citep[Theorem 3]{drori2014performance} indicates this uniquely determine $h_{i,j}$. We set $(\lambda_k)_{k=0}^N$ and $(r_{N, k})_{k=0}^{N-1}$ as 
\begin{equation}
\begin{aligned}
    &\lambda_k = \frac{\varphi_k}{\varphi_{N+1}}, \qquad k \in \{0,1,\dots, N\}
    \\
    &r_{N,k} = \frac{(\varphi_{k+1} - \varphi_k)(\varphi_{N+1} - \varphi_N)}{\varphi_{N+1}}, \qquad k \in \{0,1,\dots, N-2 \}
    \\
    &r_{N, N-1} - \lambda_N = \frac{(\varphi_{N+1} - \varphi_{N-1})(\varphi_{N+1} - \varphi_N)}{\varphi_{N+1}}.
\end{aligned}     \label{eqn::ORC-F-solution1}
\end{equation}
Moreover, we set
\begin{equation}
\begin{aligned}
 &r_{k,t} = \frac{1}{\varphi_{N + 1}} (\varphi_{k+1} - \varphi_k)(\varphi_{t+1} - \varphi_t), \qquad k \in \{1,2, \dots, N-1\}, \quad t\in \{0,1, \dots ,k-2\}
 \\
 &r_{k,k-1} - \lambda_k= \frac{1}{\varphi_{N + 1}} (\varphi_{k+1} - \varphi_k)(\varphi_{k} - \varphi_{k-1}), \qquad k \in \{ 1,2 \dots ,N-1\}.
\end{aligned}
 \label{eqn::ORC-F-solution2}
\end{equation}
In addition, we set $\hat{\fgam}$ as 
\begin{align*}
    \gamma_t &= \frac{\gamma_N}{2-\lambda_N} r_{N, t}, \qquad \qquad \qquad t \in \{0,1, \dots, N-2 \}
    \\
    \gamma_{N-1} &= \frac{\gamma_N(r_{N, N-1} - \lambda_N)}{2-\lambda_N}.
\end{align*}
Lastly, we set $\tau'$ as
\begin{align}
    \tau' = \frac{L^2}{\varphi_{N+1}},
     \label{eqn::ORC-F-solution3}
\end{align}
and $\tau = \frac{L}{2\varphi_{N+1}}$. These variables make \eqref{eqn::ORC-F-final-matrix}'s left-hand side zero. }
\end{proof}

\begin{claim}
\eqref{eqn::ORC-F-solution1}, \eqref{eqn::ORC-F-solution2} and \eqref{eqn::ORC-F-solution3} are an optimal solution of \eqref{eqn::ORC-F-hij-final-problem}.
\end{claim}
\begin{proof}
Let we represent $S$ with the variable $(r_{i,j})$. We will denote this as $\fA$. To be specific, 
\begin{align*}
   \fA(\flam, \fbet, \falp, \tau',(r_{i,j})) = S(\flam, \fgam, \tau'; (h_{i,j})) = \left(
\begin{array}{ccc}
    \frac{1}{2}\tau' & -\frac{L}{2}\hat{\fbet}^\intercal & -\frac{L}{2}\bet_N \\
    -\frac{L}{2}\hat{\fbet}^\intercal & Q(\flam; (r_{i,j})) & \fq((r_{i,j})) 
    \\
    -\frac{L}{2}\bet_N & \fq((r_{i,j}))^\intercal & \frac{2-\lambda_N}{2} \\
\end{array}
\right)  \succeq 0 .
\end{align*}
Here, $\fbet = (\hat{\fbet}^\intercal, \beta_N)^\intercal$, \begin{align*}
     Q(\flam;\revision{(r_{i,j})}) &=  \frac{\lambda_1}{2} \pg_0 \pg_0^\intercal + \sum_{k=1}^{N-1} \frac{\lambda_{k+1}-2\lambda_k}{2}\pmb{g}_k\pmb{g}_k^\intercal  + \sum_{k=1}^{N-1}\frac{\lambda_k}{2} (\pmb{g}_{k-1}-\pmb{g}_{k})(\pmb{g}_{k-1}-\pmb{g}_{k})^\intercal
    \\&+ \frac{\lambda_N}{2}\pg_{N-1}\pg_{N-1}^\intercal+ \sum_{k=1}^{N-1} \sum_{t=0}^{k-1} \left( \frac{r_{k,t}}{2} \right) \left(\pmb{g}_k\pmb{g}_t^\intercal + \pmb{g}_t\pmb{g}_k^\intercal\right),
\end{align*}
and 
\begin{align*}
    \fq((r_{i,j})) &= \sum_{t=0}^{N-1} \frac{r_{N,t}}{2}\pmb{g}_t - \frac{\lambda_N}{2} \pg_{N-1}.
\end{align*}

 
\revision{ Define a linear SDP relaxation of \eqref{eqn::ORC-F-hij-wkdual-2} as 
\begin{equation}
\begin{alignedat}{4}
&\minimize_{r_{i,j}}&&\minimize_{(\flam, \fbet, \falp, \tau') \geq \pmb{0}}\quad && {\tau'} 
    \\
    & &&\text{subject to} \quad && \fA(\flam, \fbet, \falp, \tau', (r_{i,j})) \succeq 0.
    \\
    & && &&\fB(\flam, \fbet, \falp, \tau') = \left(\flam, \fbet, \falp, \tau'\right) \geq 0
    \\
    & && &&\fC(\flam, \fbet, \falp) = (-\alpha_0 + \beta_0, -\alpha_1 + \lambda_1 + \beta_1, \dots, -\alpha_N + \lambda_N + \beta_N) = 0
    \\
    & && &&\fD(\flam, \fbet, \falp) = (-\lambda_1 + \alpha_0, -\lambda_2 + \alpha_1, \dots, -\lambda_N + \alpha_{N-1}, \alpha_N -1) = 0.
\end{alignedat}
\label{eqn::ORC-F-relaxation}
\end{equation}
\citep[Theorem 3]{drori2014performance} indicates that
 if we prove the choice in the previous claim satisfies KKT condition of \eqref{eqn::ORC-F-relaxation}, then this is also an optimal solution for the original problem \revision{since $(r_{i,j})$ uniquely determines $h_{i,j}$.}}
The Lagrangian of the minimization problem is 
\begin{align*}
    \cL(\flam, &\fbet, \falp, \tau',(r_{i,j}), \fK, \fb, \fc, \fd) 
    \\
    &= \frac{1}{2}\tau' - \tr\left\{\fA(\flam, \fbet, \falp, \tau',(r_{i,j})) \fK\right\}- \fb^\intercal \fB(\flam, \fbet, \falp, \tau') - \fc^\intercal  \fC(\flam, \fbet, \falp) - \fd^\intercal\fD(\flam, \fbet, \falp)
\end{align*}
and the KKT conditions of the minimization problems are 
    \begin{align*}
        &\fA(\flam, \fbet, \falp, \tau';(r_{i,j}))\succeq 0, \fB(\flam, \fbet, \falp, \tau') \geq 0, \fC(\flam, \fbet, \falp) =0, \fD(\flam, \fbet, \falp) = 0,
        \\
        &\nabla_{(\flam, \fbet, \falp, \tau',(r_{i,j}))}\cL(\flam, \fbet, \falp, \tau',(r_{i,j}), \fK, \fb, \fc, \fd)  = 0,
        \\
        &\fK \succeq 0, \fb \geq 0,
        \\
        &\tr\left\{\fA(\flam, \fbet, \falp, \tau',(r_{i,j})) \fK\right\}=0,  \fb^\intercal \fB(\flam, \fbet, \falp, \tau')=0,
    \end{align*}
where $\fK$ is a symmetric matrix. Here, $\fb = (\fu, \fv, \fw, s)$. We re-index $K$'s column and row starting from -1 (so $K$'s rows and columns index are $\{-1,0,1. \dots, N\}$). Now, we will show that there exist a dual optimal solution $(\fK, \fb, \fc, \fd)$ that $(\flam, \fbet, \falp, \tau',(r_{i,j}), \fK, \fb, \fc, \fd)$ satisfies KKT condition, which proves a pair $(\flam, \fbet, \falp, \tau', (r_{i,j}))$ is an optimal solution for primal problem. The stationary condition $\nabla_{(\flam, \fbet, \falp, \tau',(r_{i,j}))}\cL(\flam, \fbet, \falp, \tau',(r_{i,j}), \fK, \fb, \fc, \fd)  = 0$ can be rewritten as 
\begin{equation}
    \begin{aligned}
    &\frac{\partial \cL}{\partial \lambda_k} = -\frac{1}{2}\left(2 K_{k-1, k-1} - K_{k-1, k} - K_{k, k-1} - K_{k,k}\right) - u_k - c_k +d_{k-1} = 0, \qquad k \in \{1,2, \dots, N\}
        \\
    &\frac{\partial \cL}{\partial \beta_k} = \frac{L}{2}\left(K_{-1, k}  + K_{k, -1}\right) - v_k - c_k = 0, \qquad k  \in \{0,1, \dots, N\}
    \\
    &\frac{\partial \cL}{\partial \alpha_k} = -w_k + c_k - d_k = 0,\qquad k \in  \{0, 1, \dots, N\}
    \\
    &\frac{\partial \cL}{\partial \tau'} = \frac{1}{2} -\frac{1}{2}K_{-1, -1} - s =0
    \\
    &\frac{\partial \cL}{\partial r_{k,t}} = -\frac{1}{2}(K_{k,t} + K_{t,k}) = 0,\qquad  k \in \{1,2, \dots, N\}, \quad t \in \{0,1,\dots, k-1\}.
    \end{aligned}\label{eqn::ORC-F-KKT-stationary}
\end{equation}
We already know that $\fB(\flam, \fbet, \falp, \tau') \neq 0$, we can set $\fb = 0$. Then, \eqref{eqn::ORC-F-KKT-stationary} reduces to 
    \begin{align*}
    &K_{k,t} = 0, \qquad  k\in \{1,2, \dots, N\}, \quad t \in \{0,1,\dots, k-1\}
    \\
    &-\frac{1}{2}\left(2{K_{k-1, k-1}} - K_{k,k}\right)  - c_k + d_{k-1} = 0, \qquad  k \in \{1,2, \dots, N\}
    \\
    &LK_{-1, k} - c_k = 0, \qquad k  \in \{0,1, \dots, N\}
    \\
    & c_k - d_k = 0,\qquad k  \in \{0,1, \dots, N\}
     \\
    &K_{-1, -1}  =1.
    \end{align*}
Then, we have
\begin{equation*}
\begin{aligned}
\fK &= \left(
\begin{array}{cccccc}
    1& \frac{c_0}{L} & \frac{c_1}{L} & \dots & \frac{c_{N-1}}{L} &\frac{c_N}{L} \\
    \frac{c_0}{L} & K_{0,0} & 0 & \dots&0 & 0\\
    \vdots & \vdots & \vdots & \ddots & \vdots &\vdots  \\
    \frac{c_{N-1}}{L} & 0 & 0 & \dots & K_{N-1, N-1} & 0\\
    \frac{c_{N}}{L} & 0& 0 & \dots &0 & K_{N,N}\\
\end{array}
\right) \succeq 0
\end{aligned}
\end{equation*}
 and since $\tr\left\{\fA(\flam, \fbet, \falp, \tau',(r_{i,j})) \fK\right\}=0$ with  $\fA \succeq 0$, we can replace this condition by $\fA(\flam, \fbet, \falp, \tau',(r_{i,j})) \fK=0$.
Then the KKT condition for the given $(\flam, \fbet, \tau', (r_{i,j}))$ reduces to 
\begin{equation*}
    \begin{aligned}
        &\frac{1}{2}\tau' - \frac{1}{2}\fbet^\intercal\fc = 0
        \\
        &\frac{1}{2L}\tau'\fc - \frac{L}{2}\diag(K_{0,0}, \dots, K_{N-1, N-1}, K_{N,N}) \fbet = 0 
        \\
        &-\frac{1}{2}\fbet \fc^\intercal + \left(\begin{array}{cc}Q & \fq\\ \fq^\intercal & \frac{2 - \lambda_N}{2} \end{array}\right)\diag(K_{0,0}, \dots, K_{N-1, N-1}, K_{N,N}) = 0.
    \end{aligned}
\end{equation*}
By solving the equation, we have $c_i = (\varphi_{i+1} - \varphi_{i}) K_{i,i} $ for $i = 0,1, \dots, N$ and we have $ 1 =  \sum_{i=0}^N \frac{c_i^2}{L^2K_{i,i}}$ by the first above equation. Therefore, $K \succeq 0$.
\end{proof}

\begin{claim}
The obtained algorithm is ORC-F$_\flat$. 
\end{claim}

\begin{proof}
By calculating $(h_{i,j})$ of ORC-F$_\flat$, we can prove the equivalence of the obtained solution and ORC-F$_\flat$. Indeed, ORC-F$_\flat$ is obtained by using \citep{lee2021geometric}'s auxiliary sequences.
We will show that obtained $(\hat{h}_{i,j})$ satisfies \begin{align*}
    x_0 - \sum_{i=1}^{k+1} \sum_{j=0}^{i-1} \frac{\hat{h}_{i.j}}{L} \nabla f(x_j) &= \frac{\varphi_{k+1}}{\varphi_{k+2}} \left( x_0 - \sum_{i=1}^k \sum_{j=0}^{i-1} \frac{\hat{h}_{i,j}}{L} \nabla f(x_j) - \frac{1}{L} \nabla f(x_k)\right)
    \\&+ \left(1- \frac{\varphi_{k+1}}{\varphi_{k+2}} \right)\left(x_0 - \sum_{j=0}^k \frac{\varphi_{j+1} - \varphi_j}{L} \nabla f(x_j)\right),
\end{align*}
which is re-written form of ORC-F$_\flat$. Comparing $\nabla f(x_j)$'s each coefficient, we should prove 
\begin{align*}
    \sum_{i = j+1}^{k+1} \hat{h}_{i,j} &= \frac{\varphi_{k+1}}{\varphi_{k+2}}\sum_{i = j+1}^{k} \hat{h}_{i,j} + \left(1- \frac{\varphi_{k+1}}{\varphi_{k+2}}\right) (\varphi_{j+1} - \varphi_j) \qquad j \in \{0,1,\dots, k-1\}
    \\
    \hat{h}_{k+1,k} &= \frac{\varphi_{k+1}}{\varphi_{k+2}} + \left(1- \frac{\varphi_{k+1}}{\varphi_{k+2}}\right) (\varphi_{k+1} - \varphi_k),
\end{align*}
which is exactly equal to the recursive rule of \citep[Theorem 3]{drori2014performance}. 
\end{proof}

It now remains to justify \eqref{eqn::ORC-F-strongdual}.
Write  $\|(s_{i,j})-(s_{i,j}')\|_\infty\le \varepsilon$ if $\max_{i,j} |s_{i,j} - s_{i,j}'| \leq \varepsilon$. 
Denote the optimal value of \eqref{eqn::ORC-primal} as $p((s_{i,j}))$, i.e.,
\[
p((s_{i,j}))=\maximize_{\fG, \fF_0, \fF_1} \mkern7mu \pf_{N+1}^\intercal \fF_1.
\]
We show that $p$ is a continuous function. 
If $p$ is continuous, by Fact~\ref{fact::2},
\begin{align*}
    \argmin_{s_{i,i-1} }\maximize_{\fG, \fF_0, \fF_1} \mkern7mu \pf_{N+1}^\intercal \fF_1  =
    \argmin_{s_{i,i-1} \neq 0}\maximize_{\fG, \fF_0, \fF_1} \mkern7mu \pf_{N+1}^\intercal \fF_1  = \argmin_{s_{i,i-1} \neq 0}\minimize_{(\flam, \fbet, \falp,\tau) \geq \pmb{0}}\mkern7mu {\tau}.
\end{align*}
Since our analytic solution for $\argmin_{(s_{i,j})}\minimize_{(\flam, \fbet, \falp,\tau) \geq \pmb{0}}\mkern7mu {\tau}$ satisfies $s_{i,i-1} \neq 0$, the strong duality claim \eqref{eqn::ORC-F-strongdual} is justified.

Finally, we establish continuity of $p((s_{i,j}))$ with the following claim.

\begin{claim}
\label{claim::4}
Assume $\|(s_{i,j})-(s_{i,j}')\|_\infty\le \varepsilon$. Let $\{x_i, y_i\}_{i=0}^N$ be points with the FSFO with coefficients $(s_{i,j})$ and $\{x_i', y_i'\}_{i=0}^N$ be points with the FSFO with coefficients $(s_{i,j}')$. Assume $\{(x_i, y_i, g_i, f_{i,0}, f_{i+1, 1})\}_{i=0}^N$ satisfies $\cI_{\operatorname{ORC}}$. \revision{Assume} $\norm{x_0 - x_\star} \leq R$. 
Moreover, we can find $\{(x_i', y_i', g_i', f_{i,0}', f_{i+1, 1}')\}_{i=0}^N = \{(x_i', y_i', g_i, f_{i,0} - \frac{i^2}{L}C\varepsilon, f_{i+1, 1} - \frac{i^2}{L}C\varepsilon)\}_{i=0}^N$ that satisfies  $\cI_{\operatorname{ORC}}$, where \revision{
$C=C(\{(s_{i,j}), R,L\})$
is a constant continuously depending only on $(\{(s_{i,j}), R,L\})$.
}
\end{claim}
\begin{proof}
A continuous function $C$ only depending on  $\{(s_{i,j}), R,L\}$ bounds $ \max_{i=0}^N \norm{g_i}^2\le C$.
We first show that such a constant exists.
Under the initial condition $\norm{x_0 - x_\star}^2 \leq R^2$, 
\begin{align*}
     \frac{1}{4L}\norm{g_k}^2 + L\norm{x_k - x_\star}^2 - \frac{1}{2L}\norm{g_k}^2 &\geq \langle g_k,  x_k - x_\star \rangle - \frac{1}{2L}\norm{g_k}^2\\
     &\geq f_\star - f_{k,0} - \langle g_k, x_\star - x_k \rangle - \frac{1}{2L}\norm{g_k}^2 \geq 0,
\end{align*}
where the first inequality follows from Young's inequality,
the second inequality follows from the fact that $f_\star$ is the optimal value,
and the third inequality is the cocoercivity inequality on $(x_\star, y_k)$.
So $ L\norm{x_k - x_\star}^2 \ge \frac{1}{4L}\norm{g_k}^2$,
and if $\norm{x_k - x_\star}^2$ is bounded by a continuous function then $\norm{g_k}^2$ is also bounded by a continuous function. Inductively, if $\norm{g_0},\dots,\norm{g_{k-1}}$ are bounded by a continuous function ($\norm{g_k}$ not included), then
\begin{align*}
    \norm{x_k - x_\star} &= \norm{x_0 - x_\star - \sum_{i=0}^{k-1}{\frac{s_{k, i}}{L} g_i}}
    \\
    &\leq \norm{x_0 - x_\star} + \sum_{i=0}^{k-1} \frac{|s_{k,i}|}{L} \norm{g_i}
\end{align*}
indicates that $\norm{x_k - x_\star}$ is bounded by a continuous function.
Chaining these arguments inductively while making sure to check that the ``continuous function'' only depend on $\{(s_{i,j}), R,L\}$, we conclude $ \max_{i=0}^N \norm{g_i}^2\le C$.

Without loss of generality, assume $f_\star = f_\star' = 0$. Then,
\begin{align*}
\left(
\begin{array}{lllll}
    & \Gamma_{k, 1} =  f_{k,0} - f_{k+1, 1} - \frac{1}{2L} \norm{g_k}^2 \geq 0 , \qquad &k\in \{0,1, \dots, N\}
    \\
    & \Gamma_{k,2} = f_{k,1} - f_{k,0} - \langle g_k, y_{k} - x_{k} \rangle \geq 0   \qquad &k \in \{1, \dots, N\}
    \\
    & \Gamma_{k,3} = - f_{k,0} - \langle g_k, x_\star - x_k \rangle - \frac{1}{2L}\norm{g_k}^2 \geq 0, \qquad &k \in \{0,1, \dots, N\}
\end{array}\right).
\end{align*}
We will show
\begin{align*}
\left(
\begin{array}{lllll}
    &  f'_{k,0} - f'_{k+1, 1} - \frac{1}{2L} \norm{g_k'}^2 \geq 0 , \qquad &k\in \{0,1, \dots, N\}
    \\
    & f'_{k,1} - f'_{k,0} - \langle g_k', y'_{k} - x'_{k} \rangle \geq 0   \qquad &k \in \{1, \dots, N\}
    \\
    &  - f'_{k,0} - \langle g_k', x_\star - x_k' \rangle - \frac{1}{2L}\norm{g_k'}^2 \geq 0, \qquad &k \in \{0,1, \dots, N\}
\end{array}\right).
\end{align*}
This is equivalent to 
\begin{align*}
\left(
\begin{array}{lllll}
    &  f_{k,0} - f_{k+1, 1} - \frac{1}{2L} \norm{g_k}^2 \geq 0 , \qquad &k\in \{0,1, \dots, N\}
    \\
    & f_{k,1} - f_{k,0} + \frac{(2k - 1)C}{L}\varepsilon - \langle g_k, (y'_{k} - x'_{k}) - (y_k - x_k) + (y_k - x_k) \rangle \geq 0   \qquad &k \in \{1, \dots, N\}
    \\
    &  - f_{k,0} + \frac{k^2}{L}C\varepsilon - \langle g_k, (x_\star - x_k) + (x_k - x_k')  \rangle - \frac{1}{2L}\norm{g_k}^2 \geq 0, \qquad &k \in \{0,1, \dots, N\}
\end{array}\right).
\end{align*}
This can be reduced as
\begin{align*}
\left(
\begin{array}{lllll}
    &  \Gamma_{k,1} \geq 0 , \qquad &k\in \{0,1, \dots, N\}
    \\
    & \Gamma_{k,2} + \frac{(2k - 1)C}{L}\varepsilon - \langle g_k, (y'_{k} - x'_{k}) - (y_k - x_k)\rangle \geq 0   \qquad &k \in \{1, \dots, N\}
    \\
    &  \Gamma_{k,3} + \frac{k^2}{L}C\varepsilon - \langle g_k, (x_k - x_k')  \rangle  \geq 0, \qquad &k \in \{0,1, \dots, N\}
\end{array}\right).
\end{align*}
For the second one, 
\begin{align*}
    &\langle g_k, (y_k' - y_k) - (x_k' - x_k) \rangle 
    \\
    &= \langle g_k, (x_{k-1}' - x_{k-1})  - (x_k' - x_k) \rangle 
    \\
    &= \langle g_k, -\frac{1}{L}\sum_{i=0}^{k-2} (s_{k-1,i}' - s_{k-1,i}) g_{i}  + \frac{1}{L}\sum_{i=0}^{k-1} (s_{k,i}' - s_{k,i}) g_{i} \rangle 
    \\
    &\leq \frac{1}{L} \norm{g_k}\left( \sum_{i=0}^{k-2} \varepsilon\norm{g_i} + \sum_{i=0}^{k-1} \varepsilon \norm{g_{i}}\right)
    \\
    &\leq \frac{2k - 1}{L}C \varepsilon.
\end{align*}
For the third one, 
\begin{align*}
    &\langle g_k, (x_k' - x_k) \rangle = \langle g_k, -\frac{1}{L}\sum_{i=0}^{k-1} (s_{k,i}' - s_{k,i}) g_{i} \rangle \leq \frac{1}{L} \norm{g_k}\left(\sum_{i=0}^{k-1} \varepsilon \norm{g_{i}}\right)
    \leq \frac{k}{L}C\varepsilon,
\end{align*}
which shows the claim.
\end{proof}

Finally, we prove continuity of $p((s_{i,j}))$. 
For $0 < \varepsilon<1$, fix any $(s_{i,j})$ and $(s_{i,j}')$ that $\|(s_{i,j})-(s_{i,j}')\|_\infty\le \varepsilon$. Define $\cA((s_{i,j}))$ and $\cA((s'_{i,j}))$ as the algorithms corresponding to $(s_{i,j})$ and $(s_{i,j}')$, respectively. For any $\{(x_i, y_i, g_i, f_{i,0}, f_{i+1, 1})\}_{i=0}^N$ generated by $\cA((s_{i,j}))$ and satisfying $\cI_{\operatorname{ORC-F}}$ and $\norm{x_0 - x_\star} \leq R$, there exists $\{(x_i', y_i', g_i', f_{i,0}', f_{i+1, 1}')\}_{i=0}^N$ generated by $\cA((s_{i,j}'))$ satisfying $\cI_{\operatorname{ORC-F}}$ and $\norm{x_0 - x_\star} \leq R$, such that the performance measures difference satisfies $(f_{N+1,1} - f_\star) - (f'_{N+1,1} - f'_\star) = N^2 C((s_{i,j}), R,L) \varepsilon$, where $C$ depends only on $\{(s_{i,j}), R,L\}$. Therefore, $p((s_{i,j})) - p((s'_{i,j})) \leq N^2 C((s_{i,j}), R,L) \varepsilon$. Conversely, for $\{(x_i', y_i', g_i', f_{i,0}', f_{i+1, 1}')\}_{i=0}^N$ generated by $\cA((s'_{i,j}))$ and satisfying $\cI_{\operatorname{ORC-F}}$ and $\norm{x_0 - x_\star} \leq R$, there exists $\{(x_i, y_i, g_i, f_{i,0}, f_{i+1, 1})\}_{i=0}^N$ generated by $\cA((s_{i,j}))$ and satisfying $\cI_{\operatorname{ORC-F}}$ and $\norm{x_0 - x_\star} \leq R$, such that 
the performance measures difference satisfies $(f_{N+1,1} - f_\star) - (f'_{N+1,1} - f'_\star) = -N^2 C((s'_{i,j}), R,L) \varepsilon$. Therefore, $p((s'_{i,j})) - p((s_{i,j})) \leq N^2 C((s'_{i,j}), R,L) \varepsilon$. 
Since $C$ is a continuous function of $(s'_{i,j})$, we conclude $|p((s'_{i,j})) - p((s_{i,j}))| \rightarrow 0$ as $\varepsilon\rightarrow 0$.


To summarize, the algorithm's performance criterion $f(y_{N+1}) - f_\star$ is bounded as 
$$f(y_{N+1}) - f_\star \leq  \frac{L}{2\varphi_{N+1}}\norm{x_0 - x_\star}^2.$$
Overall, we showed that ORC-F$_\flat$ is $\cA^\star$-optimal in the sense that
$\cI_{\operatorname{ORC-F}_\flat}$
$$\operatorname{ORC-F}_\flat = \cA^\star_N(f(y_{N+1})- f_\star, \norm{x_0 - x_\star} \leq R, \cI_{\operatorname{ORC-F}_\flat}).$$
Furthermore,
\begin{align*}
    \cR(\text{ORC-F}_\flat, &f(y_{N+1})- f_\star, \norm{x_0 - x_\star} \leq R, \cI_{\operatorname{ORC-F}_\flat}) \nonumber
    \\
    &= \cR^\star(\kA_N, f(y_{N+1})- f_\star, \norm{x_0 - x_\star} \leq R, \cI_{\operatorname{ORC-F}_\flat})  = \frac{LR^2}{2\varphi_{N+1}}.\nonumber
\end{align*}

\section{Other results}
\label{sec::other-results}
We also have two more $\cA^\star$-optimal algorithms; OBL-F$_\flat$ and FGM, which will be explained in this section. Moreover, we give a conjecture about one $\cA^\star$-optimal algorithm; OBL-G$_\flat$.  
\subsection{OBL-F}
\label{ssec::OBL-F}
Optimized backtracking linesearch - function value$_\flat$ (\textbf{OBL-F$_\flat$}) is defined as
\begin{align*}
  y_{k+1} &= x_{k} - \frac{1}{L}\nabla{f(x_{k})}
\\
  z_{k+1} &= z_{k} - {\frac{k+1}{L}}\nabla{f(x_{k})} 
\\
  x_{k+1} &= \left( 1-\frac{2}{{k+3}} \right)y_{k+1} + \frac{2}{k+3}z_{k+1}
\end{align*}
for $k=0,1,\dots$ where $y_0 = z_0 = x_0$. The \emph{last-step modification} for OBL-F$_\flat$ on secondary sequence is written as
\begin{align*}
  \tilde{x}_{k}&=\frac{1}{\sqrt{\frac{k(k+1)}{2}}+1} \left(\sqrt{\frac{k(k+1)}{2}}y_{k} + z_k\right)
\end{align*}
where $k=0,1,\dots$.
\begin{theorem}[$\mathcal{A}^\star$-optimality of OBL-F$_\flat$]
\label{thm::OBL-F-optimal}
OBL-F$_\flat$ is $\cA^\star$-optimal in the sense that 
$$\operatorname{OBL-F_\flat} = \cA^\star_N(f(x_{N})- f_\star, \norm{x_0 - x_\star} \leq R, \cI_{\operatorname{OBL-F_\flat}})$$
and has the minimax optimal rate 
$$\cR^\star(\kA_N, f(x_{N}) - f_\star, \norm{x_0 - x_\star} \leq R, \cI_{\operatorname{OBL-F_\flat}}) = \frac{LR^2}{k(k+1) + \sqrt{2k(k+1)}}$$
with respect to the inequalities 
\begin{align*}
  \cI_{\operatorname{OBL-F_\flat}} = &\biggl\{f(x_{k-1}) \geq f(x_k) + \langle \nabla f(x_k), x_{k-1} - x_{k} \rangle + \frac{1}{2L}\norm{\nabla f(x_{k-1}) - \nabla f(x_k)}^2\biggr\}_{k=1}^{N} 
  \\
  &\qquad   \bigcup \biggl\{f_\star \geq f(x_k) + \langle \nabla f(x_k), x_\star - x_k \rangle \biggr\}_{k=0}^N.
\end{align*}
\end{theorem}
Note that the inequalities in $\cI_{\operatorname{OBL-F}_\flat}$ are handy for backtracking linesearches.
We defer the proof of Theorem~\ref{thm::OBL-F-optimal} to Appendix~\ref{sec::proof}.

The following corollary is a consequence of Theorem~\ref{thm::OBL-F-optimal}, but we state it separately and present a standalone proof so that we can modify it for the proof of Theorem~\ref{thm::OBL-F}.
\begin{corollary}
\label{cor::OBLGb}
Assume (\hyperlink{A1}{A1}), (\hyperlink{A2}{A2}), and (\hyperlink{A3}{A3}). OBL-F$_\flat$'s $\tilde{x}_k$-sequence and $y_k$-sequence exhibit the rate
\begin{align*}
f(\tilde{x}_{k}) -f_\star &\leq \frac{L\norm{x_0 - x_\star}^2}{k(k+1) + \sqrt{2k(k+1)}}
\end{align*}
and
\begin{align*}
f(y_{k+1}) -f_\star &\leq \frac{L\norm{x_0 - x_\star}^2}{(k+1)(k+2)}
\end{align*}
for $k=1,2,\dots$.
\end{corollary}
\begin{proof}
Let $x_{-1} = x_0$. For $k=-1,0,1,\dots $, define
\begin{align*}
    U_k = \frac{(k+1)(k+2)}{2}\left(f(x_k) - f_\star - \frac{1}{2L}\norm{\nabla f(x_k)}^2 \right) + \frac{L}{2}\norm{z_{k+1} - x_\star}^2
\end{align*}
and 
\begin{align*}
    \tilde{U}_{k} =& \left(\sqrt{\frac{k(k+1)}{2}} +\frac{k(k+1)}{2} \right) \left(f(\tilde{x}_{k}) - f_\star\right) +\frac{L}{2}\norm{ {z}_{k} - \frac{1}{L}\frac{k(k+1)}{2} \nabla f(\tilde{x}_{k}) - x_\star}^2.
\end{align*}
Then we have $U_{k+1} \stackrel{(*)}{\leq} U_{k}$ and $\tilde{U}_k \leq U_{k-1}$, which implies 
\begin{align*}
 &\left(\sqrt{\frac{k(k+1)}{2}} +\frac{k(k+1)}{2} \right) \left(f(\tilde{x}_{k}) - f_\star\right) \leq \tilde{U}_k \leq U_{k-1} \leq \dots \leq U_{-1} = \frac{L}{2}\norm{x_0 - x_\star}^2
\end{align*}
and 
\begin{align*}
 &\frac{(k+1)(k+2)}{2}\left(f(y_{k+1}) - f_\star\right) \leq U_k \leq U_{k-1} \leq \dots \leq U_{-1} = \frac{L}{2}\norm{x_0 - x_\star}^2.
\end{align*}
To complete the proof, it remains to justify the (*) part. We defer the calculations to Appendix~\ref{ss::cor5pf}.
\end{proof}
Note that the proof only utilized inequalities in $\cI_{\operatorname{OBL-F_\flat}}$.
The Lyapunov function in this proof was inspired by the Lyapunov function used in the analysis of OGM in \citep{park2021factor}.

\paragraph{Backtracking linesearch version.}
Define optimized backtracking linesearch - function value (\textbf{OBL-F}), a line backtracking version of OBL-F$_\flat$, as follows.
Initialize $L_0$ and $\eta>1$. For $k=0,1,\dots$, we define $x_{k+1}, y_{k+1}, z_{k+1}$ as
\begin{align*}
  y_{k+1} &= x_{k} - \frac{1}{L_{k+1}}\nabla{f(x_{k})}
\\
  z_{k+1} &= z_{k} - {\frac{k+1}{L_{k+1}}}\nabla{f(x_{k})} 
\\
  x_{k+1} &= \left( 1-\frac{2}{{k+3}} \right)y_{k+1} + \frac{2}{k+3}z_{k+1}
\end{align*}
with $L_{k+1} = \eta^{i_{k+1}} L_{k}$ where $y_0 = z_0 = x_0$. The backtracking linesearch finds the smallest $i_{k+1}$ such that 
\begin{align*}
    &\left(f(x_{k}) - f(x_{k+1}) - \frac{1}{2L_{k+1}}\norm{\nabla f(x_{k})- \nabla f(x_{k+1})}^2 +\langle \nabla f(x_{k+1}), x_{k+1} - x_{k} \rangle \right) \geq 0.
\end{align*}
 
\begin{theorem}
\label{thm::OBL-F}
    Assume (\hyperlink{A1}{A1}), (\hyperlink{A2}{A2}), and (\hyperlink{A3}{A3}). OBL-F exhibits the rate as 
    \begin{align*}
    f(y_{N+1}) - f_\star \leq \frac{L_N}{(N+1)(N+2)} \left(\norm{x_0 - x_\star}^2 +\sum_{k \in K}\frac{(k+1)(k+2)}{2}\left(\frac{1}{L_{k}^2} - \frac{1}{L_{k+1}^2} \right)\norm{\nabla f(x_{k+1})}^2 \right).
\end{align*}
\end{theorem}
\begin{proof}
Let $x_{-1} = x_0$. For $k = -1, 0, 1, \dots,$ define 
\begin{align*}
    U_k = \frac{(k+1)(k+2)}{2L_k}\left(f(x_k) - f_\star - \frac{1}{2L_k}\norm{\nabla f(x_k)}^2 \right) + \frac{1}{2}\norm{z_{k+1} - x_\star}^2.
\end{align*}
Then we have
\begin{align*}
U&_{k} - U_{k+1} \stackrel{(*)}{\geq}  \frac{(k+1)(k+2)}{4}\left(\frac{1}{L_{k+1}^2} - \frac{1}{L_k^2} \right) \norm{\nabla f(x_k)}^2 .
\end{align*}

If $L_{N} \geq L$ which $L$ is smoothness constant of $L$, then if we set $y_{N+1}$ as gradient $1/L_{N}$-step of $x_N$ (i.e. $y_{N+1} = x_N - \frac{1}{L_N}x_N$). We define $K$ as the set of having smooth factor-jump, then by the above relationship, we have 
\begin{align*}
    &\frac{1}{2}\norm{x_0 - x_\star}^2 
    \\
    &\geq \frac{(N+1)(N+2)}{2L_{N}}\left(f(x_N) - f_\star - \frac{1}{2L_{N}}\norm{\nabla f(x_N)}^2 \right)+ \sum_{k \in K}\frac{(k+1)(k+2)}{4}\left(\frac{1}{L_{k+1}^2} - \frac{1}{L_k^2} \right) \norm{\nabla f(x_k)}^2
    \\
    &\geq \frac{(N+1)(N+2)}{2L_{N}}\left(f(y_{N+1}) - f_\star \right) + \sum_{k \in K}\frac{(k+1)(k+2)}{4}\left(\frac{1}{L_{k+1}^2} - \frac{1}{L_k^2} \right) \norm{\nabla f(x_k)}^2
\end{align*}
which indicates 
\begin{align*}
    f(y_{N+1}) - f_\star \leq \frac{L_N}{(N+1)(N+2)} \left(\norm{x_0 - x_\star}^2 + \sum_{k \in K}\frac{(k+1)(k+2)}{2}\left(\frac{1}{L_{k}^2} - \frac{1}{L_{k+1}^2} \right) \norm{\nabla f(x_k)}^2 \right).
\end{align*}
Note that $K$ would be a sparse set (informally) that is subset of $\{1,2, \dots, N \}$. The justification of $(*)$ is deferred to Appendix~\ref{ss::thm6pf}.
\end{proof}

\paragraph{Discussion.}
The rates of OGM, OGM-simple \citep{park2021factor}, OBL-F$_\flat$, and OBL-F all have the same leading-term constants, i.e. the limit of convergence rate's ratio when $k\rightarrow\infty$ is $1$.
We clarify that although OBL-F$_\flat$ and OGM-simple \citep{park2021factor} are similar in their forms, the two algorithms are distinct.

\subsection{$\cA^\star$-optimality of FGM}
\label{ssec::FGM}
A question that motivated this work was whether FGM is an exactly optimal algorithm in some sense.
Here, we provide the answer that FGM is $\cA^\star$-optimal conditioned on a set of inequalities that are handy for both randomized coordinate updates and backtracking linesearches.
Indeed FGM does admit the variants FGM-RC$^\sharp$ and FGM-BL as discussed in Section~\ref{ss:prelim}.

\begin{theorem}[$\cA^\star$-optimality of FGM]
\label{thm::FGM-optimal}
    Nesterov's FGM is $\cA^\star$-optimal in the sense that
    $$\operatorname{FGM} = \cA^\star_N(f(y_{N+1})- f_\star, \norm{x_0 - x_\star} \leq R, \cI_{\operatorname{FGM}})$$
    and has the minimax optimal rate 
    $$\cR^\star(\kA_N, f(y_{N+1}) - f_\star, \norm{x_0 - x_\star} \leq R, \cI_{\operatorname{FGM}}) = \frac{LR^2}{2\theta_N^2}$$ with respect to the inequalities 
    \begin{align*}
      \cI_{\operatorname{FGM}} = &\biggl\{  f(x_k) \geq f(y_{k+1}) + \frac{1}{2L} \norm{\nabla f(x_k)}^2 \biggr\}_{k=0}^N  \bigcup \biggl\{f(y_{k}) \geq f(x_{k}) + \langle \nabla f(x_{k}), y_{k} - x_{k} \rangle\biggr\}_{k=1}^{N} 
      \\
      &\qquad   \bigcup \biggl\{f_\star \geq f(x_{k}) + \langle \nabla f(x_{k}), x_\star - x_k \rangle \biggr\}_{k=0}^N.
    \end{align*}
\end{theorem}
We defer the proof of Theorem~\ref{thm::FGM-optimal} to Section~\ref{sec::proof}.

\subsection{OBL-G}
\label{ssec::OBL-G}
Optimized backtracking linesearch - gradient norm$_\flat$ (\textbf{OBL-G$_\flat$}) is defined as
\begin{align*}
    y_{k+1} &= x_k - \frac{1}{L} \nabla f(x_k)
    \\
    z_{k+1} &= z_{k}-\frac{1}{L}\frac{N-k+1}{2}\nabla f(x_k)
    \\    
    x_{k+1} &= \frac{N-k-2}{N-k+2} y_{k+1}+ \frac{4}{N-k+2}  z_{k+1}
\end{align*}
for $k = 1, 2, \dots, N-1$ where $y_0 = z_0 = x_0$, and 
\begin{align*}
    y_{1} &= x_0 - \frac{1}{L} \nabla f(x_0)
    \\
    z_{1} &= z_{0}-\frac{1}{L}\frac{1 + \sqrt{\frac{N(N+1)}{2}}}{2}\nabla f(x_k)
    \\    
    x_{1} &= \frac{N-2}{N+2} y_{k+1}+ \frac{4}{N+2}  z_{k+1}.
\end{align*}

The PEP characterizing OBL-G$_\flat$ turns out to be bi-convex (hence non-convex) and this non-convexity prevents us from establishing $\cA^\star$-optimality of OBL-G$_\flat$.
This non-convexity was also present in the prior work of OGM-G by \citet{kim2021optimizing}, as we further discuss in Section~\ref{sss:discussion-oblg}.
Nevertheless, numerical evidence indicates that OBL-G$_\flat$ is likely $A^\star$-optimal, so we state the following claim as a conjecture.

\begin{conjecture}[$\cA^\star$-optimality of OBL-G$_\flat$]
\label{thm::OBL-G-optimal}
OBL-G$_\flat$ is $\cA^\star$-optimal in the sense that 
$$\operatorname{OBL-G}_\flat = \cA^\star_N( \norm{\nabla{f(x_{N})}}^2, f(x_0) - f_\star \leq \frac{1}{2}LR^2, \cI_{\operatorname{OBL-G}_\flat})$$
and has the minimax optimal rate 
$$\cR^\star (\kA_N, \norm{\nabla{f(x_{N})}}^2, f(x_0) - f_\star \leq \frac{1}{2}LR^2, \cI_{\operatorname{OBL-G_\flat}}) = 2L^2R^2 \frac{N^2 + N - \sqrt{2N(N+1)}}{N^2(N+1)^2 - 2\sqrt{2N(N+1)}}$$
with respect to the inequalities
\begin{align*}
  \cI_{\operatorname{OBL-G_\flat}} = &\left\{f(x_k) \geq f(x_{k+1}) + \langle \nabla f(x_{k+1}), x_{k+1} - x_k \rangle + \frac{1}{2L} \norm{\nabla f(x_k) - \nabla f(x_{k+1})}^2\right\}_{k=0}^{N-1}  
  \\
  & \bigcup\left\{ f(x_N) \geq f(x_k) + \langle \nabla f(x_k), x_k -x_N \rangle \right\}_{k=0}^N 
  \\
  & \bigcup \biggl\{ f(x_N) \geq f_\star + \frac{1}{2L}\norm{ \nabla f(x_{N})}^2\biggr\},
\end{align*}
\end{conjecture}
Note that the inequalities in $\cI_{\operatorname{OBL-G_\flat}}$are handy for backtracking linesearches.

Since Conjecture~\ref{thm::OBL-G-optimal} is just a conjecture, the following convergence rate of OBL-G$_\flat$ must be established as a standalone result.\footnote{
This conjecture was recently resolved on page 7 of \citep{kim2023time}, employing the concept of H-duality.}
Again, we will modify this proof later for the proof of Theorem~\ref{thm::OBL-G-convrate2}.
\begin{theorem}
\label{thm::OBL-G-convrate}
Assume (\hyperlink{A1}{A1}) (\hyperlink{A2}{A2}), and (\hyperlink{A4}{A4}).
OBL-G$_\flat$'s $x_k$-sequence exhibits the rate
\begin{align*}
\norm{\nabla f(x_N)}^2 \leq 4L \frac{N^2 + N - \sqrt{2N(N+1)}}{N^2(N+1)^2 - 2\sqrt{2N(N+1)}}(f(x_0) - f_\star) \leq  \frac{4L}{N^2}(f(x_0) - f_\star).
\end{align*}
\end{theorem}
\begin{proof}
For $k=1,2\dots, N-1$, define
\begin{align*}
     U_k&=\frac{1}{(N-k+1)(N-k+2)}\left(\frac{1}{2L}\|\nabla f(x_k)\|^2+f(x_k)-f(x_N)-\langle \nabla f(x_k) , x_k-y_k \rangle \right)\\
    &\qquad +\frac{4L}{(N-k)(N-k+1)(N-k+2)(N-k+3)} \langle z_{k}-y_k, z_{k}-x_N\rangle.   
\end{align*}
and 
\begin{align*}
    U_N = \frac{1}{4L}\norm{\nabla f(x_N)}^2,\qquad \qquad      U_0 = \frac{N(N+1) - \sqrt{2N(N+1)}}{(N-1)N(N+1)(N+2)}\left(f(x_0) - f(x_N)\right).
\end{align*}
Then, we have $U_k \stackrel{(*)}\geq U_{k+1}$, which implies 
\begin{align*}
 &\frac{1}{4L}\norm{\nabla f(x_N)}^2 = U_N \leq \dots \leq U_0 = \frac{N(N+1) - \sqrt{2N(N+1)}}{(N-1)N(N+1)(N+2)}\left(f(x_0) - f(x_N)\right)
 \\
 &\qquad \leq\frac{N(N+1) - \sqrt{2N(N+1)}}{(N-1)N(N+1)(N+2)}\left(f(x_0) - f(x_\star) - \frac{1}{2L}\norm{\nabla f(x_N)}^2 \right).
\end{align*}
To complete the proof, it remains to justify the (*) part. We defer the calculations to Appendix~\ref{ss::thm9pf}.
\end{proof}

Note that the proof only utilized inequalities in $\cI_{\operatorname{OBL-G_\flat}}$. The Lyapunov function in this proof was inspired by the Lyapunov function used in the analysis of OGM-G in \citep{lee2021geometric}. 

\paragraph{Backtracking linesearch version.}
Define optimized backtracking linesearch - gradient norm (\textbf{OBL-G}), a line backtracking version of OBL-G$_\flat$, as follows. Initialize $L_0$ and $\eta>1$. For $k = 1,2,\dots,0$ we define $x_{k+1}, y_{k+1}, z_{k+1}$ as
\begin{align*}
    y_{k+1} &= x_k - \frac{1}{L_{k+1}} \nabla f(x_k)
    \\
    z_{k+1} &= z_{k}-\frac{1}{L_{k+1}}\frac{N-k+1}{2}\nabla f(x_k)
    \\    
    x_{k+1} &= \frac{N-k-2}{N-k+2} y_{k+1}+ \frac{4}{N-k+2}  z_{k+1}
\end{align*}
and  
\begin{align*}
    y_{1} &= x_0 - \frac{1}{L_1} \nabla f(x_0)
    \\
    z_{1} &= z_{0}-\frac{1}{L_1}\frac{1 + \sqrt{\frac{N(N+1)}{2}}}{2}\nabla f(x_k)
    \\    
    x_{1} &= \frac{N-2}{N+2} y_{1}+ \frac{4}{N+2}  z_{1}.
\end{align*}
with $L_{k+1} = \eta^{i_{k+1}} L_{k}$ where $y_0 = z_0 = x_0$. The backtracking linesearch finds the smallest $i_{k+1}$ such that
\begin{align*}
    &\left(f(x_{k}) - f(x_{k+1}) - \frac{1}{2L_{k+1}}\norm{\nabla f(x_{k})- \nabla f(x_{k+1})}^2 +\langle \nabla f(x_{k+1}), x_{k+1} - x_{k} \rangle \right) \geq 0.
\end{align*}

\begin{theorem}
\label{thm::OBL-G-convrate2}
    Assume (\hyperlink{A1}{A1}), (\hyperlink{A2}{A2}), and (\hyperlink{A4}{A4}). OBL-G exhibits the rate as 
    \begin{align*}
     &\frac{1}{4L_{N}^2}\norm{\nabla f(x_N)}^2 
     \\
     &\qquad \leq  - \sum_{k \in I} \frac{1}{(N-k)(N-k+1)}\left(\frac{1}{L_k} - \frac{1}{L_{k+1}}\right) \left(f(x_k) - \frac{1}{2}\left(\frac{1}{L_k} + \frac{1}{L_{k+1}}\right) \norm{\nabla{f(x_{k})}}^2 - f(x_N)\right)
     \\
     &\qquad \qquad + \frac{1}{(N+1)^2}\left(f(x_0) - f(x_N)\right).
\end{align*}
\end{theorem}
\begin{proof}
Let $x_{-1} = x_0$. For $k=1,2\dots, N-1$, define
\begin{align*}
     U_k&=\frac{1}{(N-k+1)(N-k+2)L_k}\left(\frac{1}{2L_k}\|\nabla f(x_k)\|^2+f(x_k)-f(x_N)-\langle \nabla f(x_k) , x_k-y_k \rangle \right)\\
    &\qquad +\frac{4}{(N-k)(N-k+1)(N-k+2)(N-k+3)} \langle z_{k}-y_k, z_{k}-x_N\rangle
\end{align*}
and
\begin{align*}
    U_N = \frac{1}{4L_N^2}\norm{\nabla f(x_N)}^2,\qquad \qquad      U_0 = \frac{1}{L_0}\frac{N(N+1) - \sqrt{2N(N+1)}}{(N-1)N(N+1)(N+2)}\left(f(x_0) - f(x_N)\right).
\end{align*}
Then, we have 
\begin{align*}
U&_{k} - U_{k+1} &\stackrel{(*)}{\geq} \frac{1}{(N-k)(N-k+1)}\left(\frac{1}{L_k} - \frac{1}{L_{k+1}}\right) \left(f(x_k) - \frac{1}{2}\left(\frac{1}{L_k} + \frac{1}{L_{k+1}}\right) \norm{\nabla{f(x_{k})}}^2 - f(x_N)\right) ,
\end{align*}
which indicates
\begin{align*}
 &\frac{1}{4L_{N}^2}\norm{\nabla f(x_N)}^2 +  \sum_{k \in K} \frac{1}{(N-k)(N-k+1)}\left(\frac{1}{L_k} - \frac{1}{L_{k+1}}\right) \left(f(x_k) - \frac{1}{2}\left(\frac{1}{L_k} + \frac{1}{L_{k+1}}\right) \norm{\nabla{f(x_{k})}}^2 - f(x_N)\right) 
 \\
 &\leq \dots \leq \frac{1}{L_0}\frac{N(N+1) - \sqrt{2N(N+1)}}{(N-1)N(N+1)(N+2)}\left(f(x_0) - f(x_N)\right).
\end{align*}
where $K$ is defined as the set of having smooth factor-jump. Note that $K$ would be a sparse set (informally) that is a subset of $\{1,2, \dots, N \}$. The justification of $(*)$ is deferred to Appendix~\ref{ss::thm10pf}.
\end{proof}

\subsubsection{Discussion}
\label{sss:discussion-oblg}
The prior PEP formulations of OGM-G by \citet{kim2021optimizing} and of APPM by \citet{kim2021accelerated} share the bi-convex structure we encounter with OBL-G$_\flat$.
APPM is an accelerated algorithm for reducing the magnitude of the output of a maximal monotone operator, and the bi-convexity seems to arise from using the squared gradient magnitude, rather than the function-value suboptimality, as the performance measure.
Both OGM-G and APPM were obtained by numerically solving the bi-convex PEP.

More specifically, Kim and Fessler obtained OGM-G by solving a PEP formulation using the inequalities
\begin{align*}
  \cI_{\operatorname{OGM-G}} = &\left\{f(x_k) \geq f(x_{k+1}) + \langle \nabla f(x_{k+1}), x_{k+1} - x_k \rangle + \frac{1}{2L} \norm{\nabla f(x_k) - \nabla f(x_{k+1})}^2\right\}_{k=0}^{N-1} 
  \\
  & \bigcup \left\{ f(x_N) \geq f(x_k) + \langle \nabla f(x_k), x_k -x_N \rangle + \frac{1}{2L}\norm{\nabla f(x_k) - \nabla f(x_N)}^2 \right\}_{k=0}^N 
  \\
  &\bigcup \left\{f(x_N) \geq f_\star + \frac{1}{2L}\norm{ \nabla f(x_N)  }^2\right\}.
\end{align*}
When the bi-convex optimization problem was solved through alternating minimization, the iterates would converge to OGM-, from many different starting points. Based on this numerical evidence, we presume OGM-G is $\cA^\star$-optimal, i.e.,
$$\operatorname{OGM-G} \stackrel{?}{=} \cA^\star_N(\norm{\nabla f(x_N)}^2, f(x_0) - f_\star \leq \frac{1}{2}LR^2, \cI_{\operatorname{OGM-G}}).$$
\citet{kim2021optimizing} did prove
\begin{align*}
&\cR^\star(\norm{\nabla f(x_N)}^2, f(x_0) - f_\star \leq \frac{1}{2}LR^2, \cI_{\operatorname{OGM-G}})\\
&\qquad\qquad\qquad\qquad\le\cR(\operatorname{OGM-G},\norm{\nabla f(x_N)}^2, f(x_0) - f_\star \leq \frac{1}{2}LR^2, \cI_{\operatorname{OGM-G}})
= \frac{L^2R^2}{\tilde{\theta}_N^2},
\end{align*}
so the conjecture is that the inequality holds with equality.
Our numerical experiments for finding OBL-G$_\flat$ exhibit this same behavior, so we conjecture that OBL-G$_\flat$ is also $\cA^\star$-optimal.

The rates of OGM-G and OBL-G$_\flat$ have the same leading-term constants, i.e. the limit of convergence rate's ratio when $k\rightarrow\infty$ is $1$.
Moreover, OBL-G$_\flat$ turns out to be a ``memory-saving algorithm'' in the sense of \citep{zhou2021practical}, i.e., the coefficients of the algorithm have a non-inductive form and therefore do not need to be pre-computed.
We clarify that although OBL-G$_\flat$ and M-OGM-G \citep{zhou2021practical} are similar in their forms, the two algorithms are distinct.
In fact, the rate OBL-G$_\flat$ has a leading-term constant that is smaller (better) by a factor of $2$ compared to that of M-OGM-G \citep{zhou2021practical}.

\section{Conclusion}
In this work, we presented an algorithm design methodology based on the notion of $\cA^\star$-optimality and handy inequalities. We demonstrated the effectiveness of this methodology by finding new algorithms utilizing randomized coordinate updates and backtracking linesearches that improve upon the prior state-of-the-art rates.

By making the dependence on inequalities explicit, the notion of $\cA^\star$-optimality provides a more fine-grained understanding of the optimality algorithms, and we expect this idea to be broadly applicable to the analysis and design of optimization algorithms.
Investigating $\cA^\star$-optimal algorithms for setups with stochastic gradients \citep{taylor2019stochastic} and monotone operators and splitting methods \citep{bauschke2011convex, ryu2016primer, ryu2020operator, ryu2022LSCOMO} are interesting directions of future work.

\section*{Acknowledgements}
CP was supported by an undergraduate research internship in the second half of the 2020 Seoul National University College of Natural Sciences and the 2021 Student-Directed Education Regular Program. EKR was supported by the Institute of Information \& communications Technology Planning \& Evaluation (IITP) grant funded by the Korea government (MSIT) [NO.2021-0-01343, Artificial Intelligence Graduate School Program (Seoul National University)] and the Samsung Science and Technology Foundation (Project Number SSTF-BA2101-02). We thank Gyumin Roh and Shuvomoy Das Gupta for
providing valuable feedback.



\clearpage
%
\section*{Conflict of interest}
 The authors declare that they have no conflict of interest.

\vskip 0.2in
\bibliography{ineqtomethod}

\begin{thebibliography}{40}
\providecommand{\natexlab}[1]{#1}
\providecommand{\url}[1]{\texttt{#1}}
\expandafter\ifx\csname urlstyle\endcsname\relax
  \providecommand{\doi}[1]{doi: #1}\else
  \providecommand{\doi}{doi: \begingroup \urlstyle{rm}\Url}\fi

\bibitem[Allen-Zhu et~al.(2016)Allen-Zhu, Qu, Richt{\'a}rik, and
  Yuan]{allen2016even}
Zeyuan Allen-Zhu, Zheng Qu, Peter Richt{\'a}rik, and Yang Yuan.
\newblock Even faster accelerated coordinate descent using non-uniform
  sampling.
\newblock \emph{ICML}, 2016.

\bibitem[Bauschke and Combettes(2011)]{bauschke2011convex}
Heinz~H Bauschke and Patrick~L Combettes.
\newblock Convex analysis and monotone operator theory in {H}ilbert spaces.
\newblock \emph{New York: Springer-Verlag}, 2011.

\bibitem[Beck and Teboulle(2009)]{beck2009fast}
Amir Beck and Marc Teboulle.
\newblock A fast iterative shrinkage-thresholding algorithm for linear inverse
  problems.
\newblock \emph{SIAM Journal on Imaging Sciences}, 2\penalty0 (1):\penalty0
  183--202, 2009.

\bibitem[De~Klerk et~al.(2020)De~Klerk, Glineur, and Taylor]{de2020worst}
Etienne De~Klerk, Francois Glineur, and Adrien~B Taylor.
\newblock Worst-case convergence analysis of inexact gradient and newton
  methods through semidefinite programming performance estimation.
\newblock \emph{SIAM Journal on Optimization}, 30\penalty0 (3):\penalty0
  2053--2082, 2020.

\bibitem[Diakonikolas and Wang(2022)]{diakonikolas2021potential}
Jelena Diakonikolas and Puqian Wang.
\newblock Potential function-based framework for making the gradients small in
  convex and min-max optimization.
\newblock \emph{SIAM Journal on Optimization}, 32\penalty0 (3):\penalty0
  1668--1697, 2022.

\bibitem[Drori(2017)]{drori2017exact}
Yoel Drori.
\newblock The exact information-based complexity of smooth convex minimization.
\newblock \emph{Journal of Complexity}, 39:\penalty0 1--16, 2017.

\bibitem[Drori and Taylor(2022)]{drori2021oracle}
Yoel Drori and Adrien Taylor.
\newblock On the oracle complexity of smooth strongly convex minimization.
\newblock \emph{Journal of Complexity}, 68:\penalty0 101590, 2022.

\bibitem[Drori and Teboulle(2014)]{drori2014performance}
Yoel Drori and Marc Teboulle.
\newblock Performance of first-order methods for smooth convex minimization: a
  novel approach.
\newblock \emph{Mathematical Programming}, 145\penalty0 (1):\penalty0 451--482,
  2014.

\bibitem[Fazlyab et~al.(2018)Fazlyab, Ribeiro, Morari, and
  Preciado]{fazlyab2018analysis}
Mahyar Fazlyab, Alejandro Ribeiro, Manfred Morari, and Victor~M Preciado.
\newblock Analysis of optimization algorithms via integral quadratic
  constraints: {N}onstrongly convex problems.
\newblock \emph{SIAM Journal on Optimization}, 28\penalty0 (3):\penalty0
  2654--2689, 2018.

\bibitem[Gu and Yang(2020)]{gu2020tight}
Guoyong Gu and Junfeng Yang.
\newblock Tight sublinear convergence rate of the proximal point algorithm for
  maximal monotone inclusion problems.
\newblock \emph{SIAM Journal on Optimization}, 30\penalty0 (3):\penalty0
  1905--1921, 2020.

\bibitem[Hu and Lessard(2017)]{hu2017dissipativity}
Bin Hu and Laurent Lessard.
\newblock Dissipativity theory for {N}esterov’s accelerated method.
\newblock \emph{ICML}, 2017.

\bibitem[Kim(2021)]{kim2021accelerated}
Donghwan Kim.
\newblock Accelerated proximal point method for maximally monotone operators.
\newblock \emph{Mathematical Programming}, 190\penalty0 (1-2):\penalty0 57--87,
  2021.

\bibitem[Kim and Fessler(2016)]{kim2016optimized}
Donghwan Kim and Jeffrey~A Fessler.
\newblock Optimized first-order methods for smooth convex minimization.
\newblock \emph{Mathematical Programming}, 159\penalty0 (1):\penalty0 81--107,
  2016.

\bibitem[Kim and Fessler(2018{\natexlab{a}})]{kim2018another}
Donghwan Kim and Jeffrey~A Fessler.
\newblock Another look at the fast iterative shrinkage/thresholding algorithm
  ({F}{I}{S}{T}{A}).
\newblock \emph{SIAM Journal on Optimization}, 28\penalty0 (1):\penalty0
  223--250, 2018{\natexlab{a}}.

\bibitem[Kim and Fessler(2018{\natexlab{b}})]{kim2018generalizing}
Donghwan Kim and Jeffrey~A Fessler.
\newblock Generalizing the optimized gradient method for smooth convex
  minimization.
\newblock \emph{SIAM Journal on Optimization}, 28\penalty0 (2):\penalty0
  1920--1950, 2018{\natexlab{b}}.

\bibitem[Kim and Fessler(2021)]{kim2021optimizing}
Donghwan Kim and Jeffrey~A Fessler.
\newblock Optimizing the efficiency of first-order methods for decreasing the
  gradient of smooth convex functions.
\newblock \emph{Journal of Optimization Theory and Applications}, 188\penalty0
  (1):\penalty0 192--219, 2021.

\bibitem[Kim et~al.(2023)Kim, Ozdaglar, Park, and Ryu]{kim2023time}
Jaeyeon Kim, Asuman Ozdaglar, Chanwoo Park, and Ernest~K Ryu.
\newblock Time-reversed dissipation induces duality between minimizing gradient
  norm and function value.
\newblock \emph{NeurIPS}, 2023.

\bibitem[Lee et~al.(2021)Lee, Park, and Ryu]{lee2021geometric}
Jongmin Lee, Chanwoo Park, and Ernest~K Ryu.
\newblock A geometric structure of acceleration and its role in making
  gradients small fast.
\newblock \emph{NeurIPS}, 2021.

\bibitem[Lee and Sidford(2013)]{lee2013efficient}
Yin~Tat Lee and Aaron Sidford.
\newblock Efficient accelerated coordinate descent methods and faster
  algorithms for solving linear systems.
\newblock \emph{FOCS}, 2013.

\bibitem[Lessard et~al.(2016)Lessard, Recht, and Packard]{lessard2016analysis}
Laurent Lessard, Benjamin Recht, and Andrew Packard.
\newblock Analysis and design of optimization algorithms via integral quadratic
  constraints.
\newblock \emph{SIAM Journal on Optimization}, 26\penalty0 (1):\penalty0
  57--95, 2016.

\bibitem[Lieder(2021)]{lieder2020convergence}
Felix Lieder.
\newblock On the convergence rate of the {H}alpern-iteration.
\newblock \emph{Optimization Letters}, 15\penalty0 (2):\penalty0 405--418,
  2021.

\bibitem[Nemirovsky and Yudin(1983)]{nemirovsky1983problem}
Arkadi~Semenovich Nemirovsky and David~Borisovich Yudin.
\newblock \emph{Problem Complexity and Method Efficiency in Optimization.}
\newblock 1983.

\bibitem[Nesterov(1983)]{10029946121}
Yurii Nesterov.
\newblock A method for unconstrained convex minimization problem with the rate
  of convergence $\mathcal{O}(1/k^2)$.
\newblock \emph{Proceedings of the USSR Academy of Sciences}, 269:\penalty0
  543--547, 1983.

\bibitem[Nesterov(2012)]{nesterov2012efficiency}
Yurii Nesterov.
\newblock Efficiency of coordinate descent methods on huge-scale optimization
  problems.
\newblock \emph{SIAM Journal on Optimization}, 22\penalty0 (2):\penalty0
  341--362, 2012.

\bibitem[Nesterov and Stich(2017)]{nesterov2017efficiency}
Yurii Nesterov and Sebastian~U Stich.
\newblock Efficiency of the accelerated coordinate descent method on structured
  optimization problems.
\newblock \emph{SIAM Journal on Optimization}, 27\penalty0 (1):\penalty0
  110--123, 2017.

\bibitem[Nesterov et~al.(2020)Nesterov, Gasnikov, Guminov, and
  Dvurechensky]{nesterov2020primal}
Yurii Nesterov, Alexander Gasnikov, Sergey Guminov, and Pavel Dvurechensky.
\newblock Primal--dual accelerated gradient methods with small-dimensional
  relaxation oracle.
\newblock \emph{Optimization Methods and Software}, pages 1--38, 2020.

\bibitem[Park et~al.(2023)Park, Park, and Ryu]{park2021factor}
Chanwoo Park, Jisun Park, and Ernest~K Ryu.
\newblock Factor-2 acceleration of accelerated gradient methods.
\newblock \emph{Applied Mathematics \& Optimization}, 88\penalty0 (3):\penalty0
  77, 2023.

\bibitem[Ryu and Boyd(2016)]{ryu2016primer}
Ernest~K Ryu and Stephen Boyd.
\newblock Primer on monotone operator methods.
\newblock \emph{Applied and Computational Mathematics}, 15\penalty0
  (1):\penalty0 3--43, 2016.

\bibitem[Ryu and Yin(2022)]{ryu2022LSCOMO}
Ernest~K Ryu and Wotao Yin.
\newblock \emph{Large-Scale Convex Optimization: Algorithms \& Analyses via
  Monotone Operators}.
\newblock Cambridge University Press, 2022.

\bibitem[Ryu et~al.(2020)Ryu, Taylor, Bergeling, and
  Giselsson]{ryu2020operator}
Ernest~K Ryu, Adrien~B Taylor, Carolina Bergeling, and Pontus Giselsson.
\newblock Operator splitting performance estimation: Tight contraction factors
  and optimal parameter selection.
\newblock \emph{SIAM Journal on Optimization}, 30\penalty0 (3):\penalty0
  2251--2271, 2020.

\bibitem[Seidman et~al.(2019)Seidman, Fazlyab, Preciado, and
  Pappas]{seidman2019control}
Jacob~H Seidman, Mahyar Fazlyab, Victor~M Preciado, and George~J Pappas.
\newblock A control-theoretic approach to analysis and parameter selection of
  {D}ouglas--{R}achford splitting.
\newblock \emph{IEEE Control Systems Letters}, 4\penalty0 (1):\penalty0
  199--204, 2019.

\bibitem[Taylor and Drori(2023)]{taylor2021optimal}
Adrien Taylor and Yoel Drori.
\newblock An optimal gradient method for smooth strongly convex minimization.
\newblock \emph{Mathematical Programming}, 199\penalty0 (1-2):\penalty0
  557--594, 2023.

\bibitem[Taylor and Bach(2019)]{taylor2019stochastic}
Adrien~B Taylor and Francis Bach.
\newblock Stochastic first-order methods: non-asymptotic and computer-aided
  analyses via potential functions.
\newblock \emph{COLT}, 2019.

\bibitem[Taylor et~al.(2017{\natexlab{a}})Taylor, Hendrickx, and
  Glineur]{taylor2017exact}
Adrien~B Taylor, Julien~M Hendrickx, and Fran{\c{c}}ois Glineur.
\newblock Exact worst-case performance of first-order methods for composite
  convex optimization.
\newblock \emph{SIAM Journal on Optimization}, 27\penalty0 (3):\penalty0
  1283--1313, 2017{\natexlab{a}}.

\bibitem[Taylor et~al.(2017{\natexlab{b}})Taylor, Hendrickx, and
  Glineur]{taylor2017smooth}
Adrien~B Taylor, Julien~M Hendrickx, and Fran{\c{c}}ois Glineur.
\newblock Smooth strongly convex interpolation and exact worst-case performance
  of first-order methods.
\newblock \emph{Mathematical Programming}, 161\penalty0 (1-2):\penalty0
  307--345, 2017{\natexlab{b}}.

\bibitem[Taylor et~al.(2018)Taylor, Hendrickx, and Glineur]{taylor2018exact}
Adrien~B Taylor, Julien~M Hendrickx, and Fran{\c{c}}ois Glineur.
\newblock Exact worst-case convergence rates of the proximal gradient method
  for composite convex minimization.
\newblock \emph{Journal of Optimization Theory and Applications}, 178\penalty0
  (2):\penalty0 455--476, 2018.

\bibitem[Van~Scoy et~al.(2017)Van~Scoy, Freeman, and Lynch]{van2017fastest}
Bryan Van~Scoy, Randy~A Freeman, and Kevin~M Lynch.
\newblock The fastest known globally convergent first-order method for
  minimizing strongly convex functions.
\newblock \emph{IEEE Control Systems Letters}, 2\penalty0 (1):\penalty0 49--54,
  2017.

\bibitem[Yoon and Ryu(2021)]{yoon2021accelerated}
TaeHo Yoon and Ernest~K Ryu.
\newblock Accelerated algorithms for smooth convex-concave minimax problems
  with $\mathcal{O}(1/k^2)$ rate on squared gradient norm.
\newblock \emph{ICML}, 2021.

\bibitem[Zhang et~al.(2021)Zhang, Bao, Lessard, and Grosse]{zhang2021unified}
Guodong Zhang, Xuchan Bao, Laurent Lessard, and Roger Grosse.
\newblock A unified analysis of first-order methods for smooth games via
  integral quadratic constraints.
\newblock \emph{Journal of Machine Learning Research}, 22\penalty0
  (103):\penalty0 1--39, 2021.

\bibitem[Zhou et~al.(2022)Zhou, Tian, So, and Cheng]{zhou2021practical}
Kaiwen Zhou, Lai Tian, Anthony Man-Cho So, and James Cheng.
\newblock Practical schemes for finding near-stationary points of convex
  finite-sums.
\newblock \emph{AISTATS}, 2022.

\end{thebibliography}

\clearpage 
\appendix
\section{Deferred calculations}
\label{sec::appendix-deferred-cal}
\subsection{Missing part of Corollary~\ref{cor::OBLGb}}
\label{ss::cor5pf}
For $U_k$, we have 
\begin{align*}
U&_{k} - U_{k+1} 
\\&= \frac{(k+1)(k+2)}{2}\left(f(x_{k}) - f_\star - \frac{1}{2L}\norm{\nabla f(x_{k})}^2\right)-\frac{(k+2)(k+3)}{2} \left(f(x_{k+1}) - f_\star - \frac{1}{2L}\norm{\nabla f(x_{k+1})}^2\right) 
\\
& \qquad+\frac{L}{2}\norm{z_{k+1} - x_\star}^2 - \frac{L}{2}\norm{ z_{k+2} - x_\star}^2
\\
&=\frac{(k+1)(k+2)}{2}\left(f(x_{k}) - f_\star - \frac{1}{2L}\norm{\nabla f(x_{k})}^2\right)-\frac{(k+2)(k+3)}{2} \left(f(x_{k+1}) - f_\star - \frac{1}{2L}\norm{\nabla f(x_{k+1})}^2\right)
\\
& \qquad-\langle (k+2) \nabla f(x_{k+1}), x_\star - z_{k+1} \rangle - \frac{(k+2)^2}{2L} \norm{\nabla f(x_{k+1})}^2
\\
&=\frac{(k+1)(k+2)}{2}\left(f(x_{k}) - f_\star - \frac{1}{2L}\norm{\nabla f(x_{k})}^2\right)-\frac{(k+2)(k+3)}{2} \left(f(x_{k+1}) - f_\star + \frac{1}{2L}\norm{\nabla f(x_{k+1})}^2\right)
\\
& \qquad-\langle (k+2) \nabla f(x_{k+1}), x_\star - z_{k+1} \rangle + \frac{k+2}{2L} \norm{\nabla f(x_{k+1})}^2
\\
&=\frac{(k+1)(k+2)}{2}\left(f(x_{k}) - f(x_{k+1}) - \frac{1}{2L}\norm{\nabla f(x_{k})}^2- \frac{1}{2L}\norm{\nabla f(x_{k+1})}^2\right)
\\
&\qquad - (k+2)\left(f(x_{k+1}) - f_\star + \frac{1}{2L}\norm{\nabla f(x_{k+1})}^2\right)
\\
& \qquad \qquad -\langle (k+2) \nabla f(x_{k+1}), x_\star - z_{k+1} \rangle + \frac{k+2}{2L} \norm{\nabla f(x_{k+1})}^2
\\
&=\frac{(k+1)(k+2)}{2}\left(f(x_{k}) - f(x_{k+1}) - \frac{1}{2L}\norm{\nabla f(x_{k})}^2- \frac{1}{2L}\norm{\nabla f(x_{k+1})}^2\right)
\\ 
&\qquad - (k+2)\left(f(x_{k+1}) - f_\star\right)
\\
& \qquad\qquad-\langle (k+2) \nabla f(x_{k+1}), x_\star - x_{k+1} \rangle-\langle (k+2) \nabla f(x_{k+1}), x_{k+1} - z_{k+1} \rangle 
\\
&=\frac{(k+1)(k+2)}{2}\left(f(x_{k}) - f(x_{k+1}) - \frac{1}{2L}\norm{\nabla f(x_{k})- \nabla f(x_{k+1})}^2 +\langle \nabla f(x_{k+1}), x_{k+1} - x_{k} \rangle \right)
\\
&\qquad  +(k+2)\left(f_\star - f(x_{k+1}) - \langle \nabla f(x_{k+1}), x_\star - x_{k+1} \rangle \right)
\\
&\geq 0 
\end{align*}
which completes the proof of Corollary~\ref{cor::OBLGb}.

\clearpage
\subsection{Missing part of Theorem~\ref{thm::OBL-F}}
\label{ss::thm6pf}
For $U_k$, we have 

\begingroup
\allowdisplaybreaks
\begin{align*}
U&_{k} - U_{k+1} 
\\&= \frac{(k+1)(k+2)}{2L_{k}}\left(f(x_{k}) - f_\star - \frac{1}{2L_{k}}\norm{\nabla f(x_{k})}^2\right)-\frac{(k+2)(k+3)}{2L_{k+1}} \left(f(x_{k+1}) - f_\star - \frac{1}{2L_{k+1}}\norm{\nabla f(x_{k+1})}^2\right) 
\\
& \qquad+\frac{1}{2}\norm{z_{k+1} - x_\star}^2 - \frac{1}{2}\norm{ z_{k+2} - x_\star}^2
\\
&=\frac{(k+1)(k+2)}{2L_{k}}\left(f(x_{k}) - f_\star - \frac{1}{2L_{k}}\norm{\nabla f(x_{k})}^2\right)-\frac{(k+2)(k+3)}{2L_{k+1}} \left(f(x_{k+1}) - f_\star\right)
\\
& \qquad-\frac{1}{L_{k+1}}\langle (k+2) \nabla f(x_{k+1}), x_\star - z_{k+1} \rangle - \frac{(k+2)^2}{2L_{k+1}^2} \norm{\nabla f(x_{k+1})}^2 + \frac{(k+2)(k+3)}{4L_{k+1}^2} \norm{\nabla f(x_{k+1})}^2
\\
&=\frac{(k+1)(k+2)}{2L_{k+1}}\left(f(x_{k}) - f(x_{k+1}) - \frac{1}{2L_{k}}\norm{\nabla f(x_{k})}^2\right)
\\
&\qquad  + \left( \frac{(k+1)(k+2)}{2L_{k}} - \frac{(k+1)(k+2)}{2L_{k+1}} \right)\left(f(x_{k}) - f_\star - \frac{1}{2L_{k}}\norm{\nabla f(x_{k})}^2\right) 
\\
& \qquad \qquad -\frac{1}{L_{k+1}}\langle (k+2) \nabla f(x_{k+1}), x_\star - z_{k+1} \rangle - \frac{(k+2)^2}{2L_{k+1}^2} \norm{\nabla f(x_{k+1})}^2 + \frac{(k+2)(k+3)}{4L_{k+1}^2} \norm{\nabla f(x_{k+1})}^2
\\
& \qquad \qquad \qquad -\frac{k+2}{L_{k+1}}\left(f(x_{k+1}) - f_\star \right)
\\
&=\frac{(k+1)(k+2)}{2L_{k+1}}\left(f(x_{k}) - f(x_{k+1}) - \frac{1}{2L_{k}}\norm{\nabla f(x_{k})}^2\right)
\\ 
&\qquad + \frac{(k+2)}{L_{k+1}}\left( f_\star -f(x_{k+1}) \right)  - \frac{(k+2)^2}{2L_{k+1}^2} \norm{\nabla f(x_{k+1})}^2 + \frac{(k+2)(k+3)}{4L_{k+1}^2} \norm{\nabla f(x_{k+1})}^2
\\
& \qquad \qquad - \frac{1}{L_{k+1}}\langle (k+2) \nabla f(x_{k+1}), x_\star - x_{k+1} \rangle- \frac{1}{L_{k+1}}\langle (k+2) \nabla f(x_{k+1}), x_{k+1} - z_{k+1} \rangle 
\\
&\qquad \qquad \qquad + \left( \frac{(k+1)(k+2)}{2L_{k}} - \frac{(k+1)(k+2)}{2L_{k+1}} \right)\left(f(x_{k}) - f_\star - \frac{1}{2L_{k}}\norm{\nabla f(x_{k})}^2\right) 
\\
&=\frac{(k+1)(k+2)}{2L_{k+1}}\left(f(x_{k}) - f(x_{k+1}) - \langle \nabla f(x_{k+1}), x_k - x_{k+1}\rangle  - \frac{1}{2L_{k+1}}\norm{\nabla f(x_k) - \nabla f(x_{k+1})}^2 \right)
\\
&\qquad +\frac{(k+1)(k+2)}{2L_{k+1}} \left(  - (\frac{1}{2L_{k}}- \frac{1}{2L_{k+1}})\norm{\nabla f(x_{k})}^2 + \frac{1}{2L_{k+1}} \norm{\nabla f(x_{k+1})}^2\right)+ \frac{(k+2)(k+3)}{4L_{k+1}^2} \norm{\nabla f(x_{k+1})}^2
\\ 
&\qquad\qquad  + \frac{(k+2)}{L_{k+1}}\left( f_\star -f(x_{k+1}) - \langle \nabla f(x_{k+1}), x_\star - x_{k+1} \rangle \right)  - \frac{(k+2)^2}{2L_{k+1}^2} \norm{\nabla f(x_{k+1})}^2 
\\
&\qquad \qquad \qquad + \left( \frac{(k+1)(k+2)}{2L_{k}} - \frac{(k+1)(k+2)}{2L_{k+1}} \right)\left(f(x_{k}) - f_\star - \frac{1}{2L_{k}}\norm{\nabla f(x_{k})}^2\right) 
\\
&\geq \frac{(k+1)(k+2)}{2L_{k+1}} \left(  - (\frac{1}{2L_{k}}- \frac{1}{2L_{k+1}})\norm{\nabla f(x_{k})}^2 + \frac{1}{2L_{k+1}} \norm{\nabla f(x_{k+1})}^2\right)
\\ 
&\qquad  - \frac{(k+2)^2}{2L_{k+1}^2} \norm{\nabla f(x_{k+1})}^2 + \frac{(k+2)(k+3)}{4L_{k+1}^2} \norm{\nabla f(x_{k+1})}^2
\\
&\qquad \qquad + \left( \frac{(k+1)(k+2)}{2L_{k}} - \frac{(k+1)(k+2)}{2L_{k+1}} \right)\left(- \frac{1}{2L_{k}}\norm{\nabla f(x_{k})}^2\right) 
\\
&= \frac{(k+1)(k+2)}{4}\left(\frac{1}{L_{k+1}^2} - \frac{1}{L_k^2} \right) \norm{\nabla f(x_k)}^2  
\end{align*}
which completes the proof of Theorem~\ref{thm::OBL-F}.
\endgroup

\clearpage
\subsection{Missing part of Theorem~\ref{thm::OBL-G-convrate}}
\label{ss::thm9pf}
For $k = 1,2,\dots, N-2$, we have
\begin{align*}
    &\frac{4L}{(N-k)(N-k+1)(N-k+2)(N-k+3)} \left\langle z_{k}-y_k, z_{k}-x_N \right\rangle 
    \\ &\qquad - \frac{4L}{(N-k-1)(N-k)(N-k+1)(N-k+2)} \langle z_{k+1}-y_{k+1}, z_{k+1}-x_N\rangle
    \\
    &= \frac{4L}{(N-k)(N-k+1)(N-k+2)}\Biggl( \left\langle \frac{1}{N-k+3} \left( z_k - y_k \right),  z_{k}-x_N \right\rangle 
    \\
    &\qquad -\left\langle \frac{1}{N-k-1} \left( z_{k+1} - y_{k+1} \right),  z_{k+1}-x_N\right\rangle \Biggr)
    \\
    &= \frac{4L}{(N-k)(N-k+1)(N-k+2)}\Biggl( \left\langle \frac{1}{N-k+3} \left( z_k - y_k \right),  z_{k}-x_N \right\rangle 
    \\
    &\qquad  -\left\langle \frac{1}{N-k-1} \left( z_{k+1} - y_{k+1} \right),  z_{k} - \frac{1}{L}\frac{N-k+1}{2}\nabla f(x_k) -x_N\right\rangle \Biggr)
    \\
    &= \frac{4L}{(N-k)(N-k+1)(N-k+2)}\Biggl( \left\langle \frac{1}{N-k-1} \left( z_k - x_k \right),  z_{k}-x_N \right\rangle 
    \\
    &\qquad -\left\langle \frac{1}{N-k-1} \left( z_{k+1} - y_{k+1} \right),  z_{k} - \frac{1}{L}\frac{N-k+1}{2}\nabla f(x_k) -x_N\right\rangle \Biggr)
    \\
    &= \frac{4L}{(N-k-1)(N-k)(N-k+1)(N-k+2)}\Biggl( \left\langle  z_k - x_k ,  z_{k}-x_N \right\rangle 
    \\
    &\qquad -\left\langle  z_{k}  - x_{k} -\frac{1}{L}\frac{N-k-1}{2} \nabla f(x_k),  z_{k}  -x_N - \frac{1}{L}\frac{N-k+1}{2}\nabla f(x_k)\right\rangle \Biggr)
    \\
    &= \frac{4}{(N-k-1)(N-k)(N-k+1)(N-k+2)} \times 
    \\
    &\qquad \left\langle \nabla f(x_k), (N-k)z_k - \frac{N-k-1}{2} x_N - \frac{N-k+1}{2} x_k - \frac{1}{L}\frac{N-k-1}{2}\frac{N-k+1}{2} \nabla f(x_k) \right\rangle
    \\
    &= \frac{4}{(N-k-1)(N-k)(N-k+1)(N-k+2)} \times 
    \\
    &\qquad \left\langle \nabla f(x_k), (N-k)(z_k-x_k) + \frac{N-k-1}{2} (x_k- x_N) - \frac{1}{L}\frac{N-k-1}{2}\frac{N-k+1}{2} \nabla f(x_k) \right\rangle
    \\
    &= \frac{4}{(N-k-1)(N-k)(N-k+1)(N-k+2)} \times 
    \\
    &\qquad \left\langle \nabla f(x_k), \frac{(N-k)(N-k-1)}{4}(x_k-y_k) + \frac{N-k-1}{2} (x_k- x_N) - \frac{1}{L}\frac{N-k-1}{2}\frac{N-k+1}{2} \nabla f(x_k) \right\rangle.        
\end{align*}
Therefore, for $U_k$, we have 
\begin{align*}
U&_{k} - U_{k+1} 
\\
&=\frac{1}{(N-k+1)(N-k+2)}\left(\frac{1}{2L}\|\nabla f(x_k)\|^2+f(x_k)-f(x_N)-\langle \nabla f(x_k) , x_k-y_k \rangle \right)
\\
& - \frac{1}{(N-k)(N-k+1)}\left(\frac{1}{2L}\|\nabla f(x_{k+1})\|^2+f(x_{k+1})-f(x_N)-\langle \nabla f(x_{k+1}) , x_{k+1}-y_{k+1} \rangle \right)
\\
& + \frac{4}{(N-k-1)(N-k)(N-k+1)(N-k+2)} \times 
\\
&\left\langle \nabla f(x_k), \frac{(N-k)(N-k-1)}{4}(x_k-y_k) + \frac{N-k-1}{2} (x_k- x_N) - \frac{1}{L}\frac{N-k-1}{2}\frac{N-k+1}{2} \nabla f(x_k) \right\rangle   
\\
&= \frac{1}{(N-k)(N-k+1)} \left(f(x_k) - f(x_{k+1}) + \langle \nabla f(x_{k+1}), x_{k+1} - y_{k+1} \rangle - \frac{1}{2L} \norm{\nabla f(x_k)}^2 - \frac{1}{2L} \norm{\nabla f(x_{k+1})}^2 \right) 
\\
&+ \frac{2}{(N-k)(N-k+1)(N-k+2)} \left(f(x_N) - f(x_k) + \langle \nabla f(x_k), x_k - x_N \rangle \right),
\end{align*}
which completes the proof of Theorem~\ref{thm::OBL-G-convrate}.
\subsection{Missing part of Theorem~\ref{thm::OBL-G-convrate2}}
\label{ss::thm10pf}
For $U_k$, we have 
\begin{align*}
U&_{k} - U_{k+1} 
\\
&=\frac{1}{(N-k+1)(N-k+2)L_k}\left(\frac{1}{2L_k}\|\nabla f(x_k)\|^2+f(x_k)-f(x_N)-\langle \nabla f(x_k) , x_k-y_k \rangle \right)
\\
&  - \frac{1}{(N-k)(N-k+1)L_{k+1}}\left(\frac{1}{2L_{k+1}}\|\nabla f(x_{k+1})\|^2+f(x_{k+1})-f(x_N)-\langle \nabla f(x_{k+1}) , x_{k+1}-y_{k+1} \rangle \right)
\\
&  + \frac{4}{(N-k-1)(N-k)(N-k+1)(N-k+2)L_k} \times 
\\
&\qquad \left\langle \nabla f(x_k), \frac{(N-k)(N-k-1)}{4}(x_k-y_k) + \frac{N-k-1}{2} (x_k- x_N) - \frac{1}{L_k}\frac{N-k-1}{2}\frac{N-k+1}{2} \nabla f(x_k) \right\rangle.    
\\
&= \frac{1}{(N-k)(N-k+1)L_{k+1}} \biggl(f(x_k) - f(x_{k+1}) + \langle \nabla f(x_{k+1}), x_{k+1} - y_{k+1} \rangle 
\\
&\qquad - \frac{1}{2L_{k+1}} \norm{\nabla f(x_k)}^2 - \frac{1}{2L_{k+1}} \norm{\nabla f(x_{k+1})}^2 \biggr)
\\
&+ \frac{2}{(N-k)(N-k+1)(N-k+2)L_k} \left(f(x_N) - f(x_k) + \langle \nabla f(x_k), x_k - x_N \rangle \right)
\\
&+ \frac{1}{(N-k)(N-k+1)}\left(\frac{1}{L_k} - \frac{1}{L_{k+1}}\right) \left(f(x_k) - \frac{1}{2}\left(\frac{1}{L_k} + \frac{1}{L_{k+1}}\right) \norm{\nabla{f(x_{k})}}^2 - f(x_N)\right) 
\\ &\geq \frac{1}{(N-k)(N-k+1)}\left(\frac{1}{L_k} - \frac{1}{L_{k+1}}\right) \left(f(x_k) - \frac{1}{2}\left(\frac{1}{L_k} + \frac{1}{L_{k+1}}\right) \norm{\nabla{f(x_{k})}}^2 - f(x_N)\right) ,
\end{align*}
which completes the proof of Theorem~\ref{thm::OBL-G-convrate2}.
\section{Proofs of $\mathcal{A}^\star$-optimality}
\label{sec::proof}
In this section, we prove Theorems~\ref{thm::OBL-F-optimal} and \ref{thm::FGM-optimal}, and discuss Conjecture~\ref{thm::OBL-G-optimal} using the PEP machinery.
Again, to verify the lengthy calculations, we provide Matlab scripts verifying the analytical solutions of the SDPs: \url{https://github.com/chanwoo-park-official/A-star-map/}.

\subsection{Proof of $\mathcal{A}^\star$-optimality of OBL-F$_\flat$}
To obtain OBL-F as an $\mathcal{A}^\star$-optimal algorithm, set $f(x_{N}) - f_\star$ to be the performance measure and $\norm{x_0 - x_\star} \leq R$ to be the initial condition. Since the constraints and the objective of the problem \revision{are homogenous}, we assume $R = 1$ without loss of generality. For the argument of homogeneous, we refer \revision{to} \citep{drori2014performance, kim2016optimized,taylor2017smooth}. We use the set of inequalities that are handy for backtracking linesearches: 
\begin{align*}
  \cI_{\operatorname{OBL-F}_\flat} = &\biggl\{f_{k-1,0} \geq f_{k,0} + \langle g_k, x_{k-1} - x_{k} \rangle + \frac{1}{2L}\norm{g_{k-1} - g_{k}}^2\biggr\}_{k=1}^{N} 
  \\
  &\qquad   \bigcup \biggl\{f_\star \geq f_{k,0} + \langle g_k, x_\star - x_k \rangle \biggr\}_{k=0}^N.
\end{align*}
For calculating $\cR(\cA_N, \cP, \cC, \cI_{\operatorname{OBL-F}_\flat}) $ with fixed $\cA_N$, define the PEP with 
\begin{align*}
&\cR(\cA_N, \cP, \cC, \cI_{\operatorname{OBL-F}_\flat}) 
\\
&= \left(
\begin{array}{lllll}
    &\maximize \quad & f_{N,0}  -  f_\star 
    \\
    &\text{subject to} \quad &1\mkern18mu \quad \geq& \norm{x_0 - x_\star}^2 
    \\
    & &f_{k-1,0} \geq& f_{k,0} + \langle g_k, x_{k-1} - x_{k} \rangle + \frac{1}{2L}\norm{g_{k-1} - g_{k}}^2, & k \in \{1,2,\dots, N\}
    \\
    & & f_\star \mkern29mu \geq& f_{k,0} + \langle g_k, x_\star - x_k \rangle, & k \in \{0,1,\dots, N \}
    \\
    & & \revision{x_k} &\revision{\text{is following the algorithm }\cA_N.}
\end{array}\right)
\end{align*}
Using the notation of Section~\ref{ss::notation}, we reformulate \revision{the} problem of computing the risk $\cR(\cA_N, \cP, \cC, \cI_{\operatorname{OBL-F}_\flat})$ as the following SDP:
\begin{alignat*}{4}
&\maximize_{\fG, \fF_0} \quad  &&\pf_{N}^\intercal \fF_0
\\
&\mbox{subject to}\quad 
    &&1 \geq \px_0^\intercal \fG \px_0  
    \\
    & && 0 \geq (\pmb{f}_{k} - \pmb{f}_{k-1})^\intercal \fF_0 + \pmb{g}_{k}^\intercal \fG (\pmb{x}_{k-1} - \pmb{x}_{k}) + \frac{1}{2L} (\pmb{g}_{k-1} - \pmb{g}_{k})^\intercal \fG(\pmb{g}_{k-1} - \pmb{g}_{k}) , \qquad &&k \in \{1,2, \dots, N\}
    \\
    & && 0 \geq \pf_k^\intercal \fF_0 - \pg_k^\intercal \fG \px_k, \qquad &&k \in \{0,1, \dots, N\}.
    \\
    & && \revision{\fG \succeq 0, \fF_0 \geq 0}
\end{alignat*}
For above transformation, $d\geq N+2 $ is used \citep{taylor2017smooth}. \revision{The Lagrangian} of the optimization problem becomes 
\begin{align*}
    \Lambda(\fF_0&, \fG, \flam, \fbet, \tau)
    \\
    &=  - \pmb{f}_{N}^\intercal \fF_0 + \tau (\px_0^\intercal \fG \px_0 -1 )+\sum_{k=0}^{N} \beta_k \left( \pf_k^\intercal \fF_0 - \pg_k^\intercal \fG \px_k\right) 
    \\& + \sum_{k=1}^{N} \lambda_k \left( (\pmb{f}_{k} - \pmb{f}_{k-1})^\intercal \fF_0 + \pmb{g}_{k}^\intercal \fG (\pmb{x}_{k-1} - \pmb{x}_{k}) + \frac{1}{2L} (\pmb{g}_{k-1} - \pmb{g}_{k})^\intercal \fG(\pmb{g}_{k-1} - \pmb{g}_{k})\right)
\end{align*}
with dual variables $\flam = (\lambda_1, \dots, \lambda_{N}) \in \bR_{+}^{N}$, $\fbet = (\beta_0, \dots, \beta_N) \in \bR_{+}^{N+1}$, and $\tau \geq 0$. 

Then the dual formulation of PEP problem is 
\begin{equation}
\begin{alignedat}{2}
    &\maximize_{(\flam, \fbet,\tau) \geq \pmb{0}}\quad  &&{-\tau} 
    \\
     &\mbox{subject to} &&\pmb{0} = -\pf_{N} + \sum_{k=1}^{N} \lambda_k (\pmb{f}_{k} - \pf_{k-1}) +\sum_{k=0}^{N} \beta_k \pmb{f}_k
     \\
     & &&0 \preceq S(\flam, \fbet, \tau), 
\end{alignedat} \label{eqn::OBL-F-dualPEP}    
\end{equation}
where $S$ is defined as
\begin{align*}
    S(\flam, \fbet,  \tau) &=  \tau \pmb{x}_0\pmb{x}_0^\intercal + \sum_{k=0}^{N}\frac{\beta_k}{2}\left(-\pmb{g}_k \pmb{x}_k^\intercal-\pmb{x}_k \pmb{g}_k^\intercal  \right)
    \\
    &+ \sum_{k=1}^{N} \frac{\lambda_k}{2} \left(\pmb{g}_{k} (\pmb{x}_{k-1} - \pmb{x}_{k})^\intercal +   (\pmb{x}_{k-1} - \pmb{x}_{k})\pmb{g}_{k}^\intercal + \frac{1}{L} (\pmb{g}_{k-1} - \pmb{g}_{k})(\pmb{g}_{k-1} - \pmb{g}_{k})^\intercal\right).
\end{align*}
We have a strong duality argument
\begin{align*}
    \argmin_{s_{i,j}}\maximize_{\fG, \fF_0} \mkern7mu \pf_{N+1}^\intercal \fF_0  = \argmin_{h_{i,j}}\minimize_{(\flam, \fbet, \tau) \geq \pmb{0}}\mkern7mu {\tau}, 
\end{align*}
 as ORC-F's optimality proof. Remind that \eqref{eqn::OBL-F-dualPEP} finds the ``best'' proof for the algorithm. Now we investigate the optimization step for algorithm. The last part is minimizing  \eqref{eqn::OBL-F-dualPEP} with stepsize, i.e. 
\begin{alignat}{3}
    &\minimize_{h_{i,j}}&&\maximize_{(\flam, \fbet, ,\tau) \geq \pmb{0}}\quad && {\tau}  \label{eqn::OBL-F-hij-wkdual}
    \\
    & && \text{subject to} \quad &&\pmb{0} = -\pf_{N} + \sum_{k=1}^{N} \lambda_k (\pmb{f}_{k} - \pf_{k-1}) +\sum_{k=0}^{N} \beta_k \pmb{f}_k \label{eqn::OBL-F-const1}
    \\
     & && &&0 \preceq S(\flam, \fbet, \tau). \label{eqn::OBL-F-Smatrix}
\end{alignat}
We note that $\pf_i$ is \revision{a} standard unit vector mentioned in \eqref{eqn::egf} (not a variable), we can write \eqref{eqn::OBL-F-const1} as 
\begin{equation}
\begin{alignedat}{2}
 &&\left(\begin{array}{ll}
    &\beta_k =  \lambda_{k+1} - \lambda_{k}, \quad  k \in \{1, \dots, N-1\}
    \\
    &\beta_0 = \lambda_1
    \\
    &\beta_N =  1 - \lambda_{N}.
\end{array}\right) 
\end{alignedat}\label{eqn::OBL-F-constlinear}
\end{equation}
We consider \eqref{eqn::OBL-F-Smatrix} with \eqref{eqn::OBL-F-constlinear} and FSFO's $h_{i,j}$. To be specific, we substitute $\fbet$ to $\flam$ in $S(\flam, \fbet, \tau)$. To show the dependency of $S$ to $(h_{i,j})$ since $\px_k$ are represented with $(h_{i,j})$, we will explicitly write $S$ as $S(\flam, \tau;(h_{i,j}))$. Then, we get
\begin{align*}
    S(\flam,  \tau;(h_{i,j})) &=\tau \pmb{x}_0\pmb{x}_0^\intercal + \sum_{k=1}^N\frac{\lambda_k}{2L} (\pmb{g}_{k-1}-\pmb{g}_{k})(\pmb{g}_{k-1}-\pmb{g}_{k})^\intercal 
    \\&+ \sum_{k=1}^{N-1} \sum_{t=0}^{k-1} \left( \frac{\lambda_k}{2} \frac{h_{k,t}}{L} + \frac{\lambda_{k+1} - \lambda_{k}}{2}\sum_{j=t+1}^{k}\frac{h_{j,t}}{L} \right) \left(\pmb{g}_k\pmb{g}_t^\intercal + \pmb{g}_t\pmb{g}_k^\intercal\right) 
    \\&+ \sum_{t=0}^{N-1} \left( \frac{\lambda_N}{2} \frac{h_{N,t}}{L} + \frac{1 - \lambda_{N}}{2}\sum_{j=t+1}^{N}\frac{h_{j,t}}{L} \right) \left(\pmb{g}_N\pmb{g}_t^\intercal + \pmb{g}_t\pmb{g}_N^\intercal\right)
    \\& - \sum_{k=1}^{N-1} \frac{\lambda_{k+1} - \lambda_{k}}{2} \left(\pmb{x}_0\pmb{g}_k^\intercal + \pmb{g}_k\pmb{x}_0^\intercal\right) - \frac{\lambda_1}{2}\left(\pmb{x}_0\pmb{g}_0^\intercal + \pmb{g}_0\pmb{x}_0^\intercal\right) - \frac{1 - \lambda_N}{2}\left(\pmb{x}_0\pmb{g}_N^\intercal + \pmb{g}_N\pmb{x}_0^\intercal\right).
\end{align*}

Using the fact that $\px_0, \pg_i, \pf_i$ are unit vectors, we can represent $S(\flam, \tau;(h_{i,j}))$ with $\fgam(\flam) = -L\fbet = -L(\lambda_1, \lambda_2 - \lambda_1, \dots, 1-\lambda_N) = (\hat{\fgam}(\flam), \gamma_{N}(\flam))
$ and  $\tau' = 2L\tau$ as 
\begin{align*}
    S(\flam, \tau';(h_{i,j}))  &= \frac{1}{L}\left(
\begin{array}{ccc}
    \frac{1}{2}\tau' & \frac{1}{2}\hat{\fgam}(\flam)^\intercal & \frac{1}{2}\gamma_{N}(\flam) \\
    \frac{1}{2}\hat{\fgam}(\flam) & Q(\flam; (h_{i,j})) & q(\flam; (h_{i,j})) 
    \\
    \frac{1}{2}\gamma_{N}(\flam) & q(\flam; (h_{i,j}))^\intercal & \frac{\lambda_N}{2} \\
\end{array}
\right) \succeq 0.
\end{align*}
\revision{Here, $Q$ and $\fq$ are defined as 
\begin{align*}
     Q(\flam;(h_{i,j})) &=  \sum_{k=1}^{N-1}\frac{\lambda_k}{2} (\pmb{g}'_{k-1}-\pmb{g}'_{k})(\pmb{g}'_{k-1}-\pmb{g}'_{k})^\intercal + \frac{\lambda_N}{2}\pg'_{N-1}\pg_{N-1}^{'\intercal}
    \\&+ \sum_{k=1}^{N-1} \sum_{t=0}^{k-1} \left( \frac{\lambda_k}{2} h_{k,t} + \frac{\lambda_{k+1} - \lambda_{k}}{2}\sum_{j=t+1}^{k}h_{j,t} \right) \left(\pmb{g}'_k\pmb{g}_t^{'\intercal} + \pmb{g}'_t\pmb{g}_k^{'\intercal}\right)
\end{align*}
and 
\begin{align*}
    \fq(\flam;(h_{i,j})) &= -\frac{\lambda_N}{2}\pg_{N-1}' + \sum_{t=0}^{N-1} \left( \frac{\lambda_N}{2} h_{N,t} + \frac{1 - \lambda_{N}}{2}\sum_{j=t+1}^{N}h_{j,t} \right) \pmb{g}'_t
    \\
    &= \sum_{t=0}^{N-2} \left( \frac{\lambda_N}{2} h_{N,t} + \frac{1 - \lambda_{N}}{2}\sum_{j=t+1}^{N}h_{j,t} \right) \pmb{g}'_t +  \left(\frac{1}{2}{h_{N,N-1}} -\frac{\lambda_N}{2}\right)\pmb{g}'_{N-1}.
\end{align*}
where $\pg_k' = e_{k+1} \in \bR^{N+1}$.} Note that \eqref{eqn::OBL-F-hij-wkdual} is equivalent to 
\begin{alignat}{3}
    &\minimize_{h_{i,j}}&&\minimize_{(\flam, \tau') \geq \pmb{0}}\quad && {\tau'}  \label{eqn::OBL-F-hij-wkdual-2}
    \\
    & && \text{subject to} \quad && \left(
\begin{array}{ccc}
    \frac{1}{2}\tau' & \frac{1}{2}\hat{\fgam}(\flam)^\intercal & \frac{1}{2}\gamma_{N}(\flam) \\
    \frac{1}{2}\hat{\fgam}(\flam) & Q(\flam; (h_{i,j})) & q(\flam; (h_{i,j})) 
    \\
    \frac{1}{2}\gamma_{N}(\flam) & q(\flam; (h_{i,j}))^\intercal & \frac{\lambda_N}{2} \\
\end{array}
\right) \succeq 0 \nonumber
\end{alignat}
and dividing this optimized value with $2L$ gives the optimized value of \eqref{eqn::OBL-F-hij-wkdual}. Using Schur complement \citep{golub1996matrix}, \eqref{eqn::OBL-F-hij-wkdual-2} can be converted to the problem as 
\begin{alignat}{3}
    &\minimize_{h_{i,j}}&&\minimize_{(\flam,  \tau') \geq \pmb{0}}\quad && {\tau'}  \label{eqn::OBL-F-hij-final-problem}
    \\
    & && \text{subject to} \quad && \left(
\begin{array}{cc}
    Q - \frac{2\fq \fq^\intercal}{\lambda_N} & \frac{1}{2}( \hat{\fgam}(\flam) - \frac{2\fq \gam_N(\flam)}{\lambda_N}) \\
    \frac{1}{2}( \hat{\fgam}(\flam) - \frac{2\fq \gam_N(\flam)}{\lambda_N})^\intercal & \frac{1}{2} (\tau' - \frac{\gam_N(\flam)^2}{\lambda_N}) \\
\end{array}
\right) \succeq 0. \label{eqn::OBL-F-final-matrix}
\end{alignat}
So far we simplified SDP. We will have \revision{three steps}: finding variables that \revision{make} \eqref{eqn::OBL-F-final-matrix}'s left hand side zero, showing that the solution from \revision{the} first step satisfies KKT condition, and finally showing that obtained algorithm is equivalent to OBL-F$_\flat$. 
\begin{claim}
\label{claim::4}
    There is a point that makes \eqref{eqn::OBL-F-final-matrix}'s left-hand side zero. 
\end{claim}
\begin{proof}
\revision{Defining $\{r_{k,t}\}_{k=1,2, \dots, N,  t=0, \dots, k-1}$ as $$r_{k,t} = \lambda_k h_{k,t}- \frac{\gamma_k}{L}\sum_{j=t+1}^k h_{j,t}.$$
Then, if $r_{i,j}$ is determined, \citep[Theorem 5.1]{drori2014performance} indicates this uniquely determine $h_{i,j}$. We set $s_N = \frac{N(N+1)}{2}$, $T =\frac{1}{s_N +  \sqrt{s_N}} $, and set $(\lambda_k)_{k=0}^N$ and $(r_{N,k})_{k=0}^{N-1}$ as 
\begin{equation}
    \begin{aligned}
        &\lambda_{k} = \frac{k(k+1)}{2} T, \qquad k \in \{1,2, \dots, N\}
        \\
        & r_{N,k} = \frac{k+1}{\sqrt{s_N} + 1}, \qquad k \in \{0,2, \dots, N-2\}
        \\
        &r_{N, N-1} - \lambda_N = \frac{N}{\sqrt{s_N} + 1}.
    \end{aligned}\label{eqn::OBL-F-solution1}
\end{equation}
Moreover, we set \begin{equation}
\begin{aligned}
 &r_{k,t} = \frac{(k+1)(t+1)}{s_N + \sqrt{s_N}}, \qquad k \in \{1,2, \dots, N-1\}, \quad t\in \{0,1, \dots ,k-2\}
 \\
 &r_{k,k-1} = \frac{k(k+1) + s_N}{s_N + \sqrt{s_N}}, \qquad k \in \{ 1,2 \dots ,N-1\}.
\end{aligned}
 \label{eqn::OBL-F-solution2}
\end{equation}
In addition, we set $\hat{\fgam}$ as 
\begin{align*}
    \gamma_t &= \frac{\gamma_N}{\lambda_N} r_{N, t}, \qquad  t \in \{0,1, \dots, N-2\} 
    \\
    \gamma_{N-1} &= \frac{\gamma_N(r_{N, N-1} - \lambda_N)}{\lambda_N}.
\end{align*}
Lastly, we set $\tau'$ as 
\begin{align}
    \tau' =\frac{L^2}{s_N + \sqrt{s_N}},
    \label{eqn::OBL-F-solution3}
\end{align}
and $\tau = \frac{L}{2(s_N + \sqrt{s_N})} = \frac{L}{N(N+1) + \sqrt{2N(N+1)}}$. These variables make \eqref{eqn::OBL-F-final-matrix}'s left-hand side zero. }
\end{proof}

\begin{claim}
\eqref{eqn::OBL-F-solution1}, \eqref{eqn::OBL-F-solution2} and \eqref{eqn::OBL-F-solution3} are an optimal solution of \eqref{eqn::OBL-F-hij-final-problem}.
\end{claim}
\begin{proof}
Let we represent $S$ with the variable $(r_{i,j})$. We will denote this as $\fA$. To be specific, 
\begin{align*}
   \fA(\flam, \fbet, \tau',(r_{i,j})) = S(\flam, \fgam, \tau'; (h_{i,j})) = \left(
\begin{array}{ccc}
    \frac{1}{2}\tau' & -\frac{L}{2}\hat{\fbet}^\intercal & -\frac{L}{2}\bet_N \\
    -\frac{L}{2}\hat{\fbet}^\intercal & Q(\flam; (r_{i,j})) & \fq((r_{i,j})) 
    \\
    -\frac{L}{2}\bet_N & \fq((r_{i,j}))^\intercal & \frac{\lambda_N}{2} \\
\end{array}
\right)  \succeq 0 .
\end{align*}
Here, $\fbet = (\hat{\fbet}^\intercal, \beta_N)^\intercal$, 
\begin{align*}
     Q(\flam;(r_{i,j})) = \sum_{k=1}^{N-1}\frac{\lambda_k}{2} (\pmb{g}_{k-1}-\pmb{g}_{k})(\pmb{g}_{k-1}-\pmb{g}_{k})^\intercal + \frac{\lambda_N}{2}\pg_{N-1}\pg_{N-1}^\intercal+ \sum_{k=1}^{N-1} \sum_{t=0}^{k-1} \frac{r_{k,t}}{2} \left(\pmb{g}_k\pmb{g}_t^\intercal + \pmb{g}_t\pmb{g}_k^\intercal\right)
\end{align*}
and 
\begin{align*}
    \fq((r_{i,j})) &= \sum_{t=0}^{N-1} \frac{r_{N,t}}{2}\pmb{g}_t - \frac{\lambda_N}{2} \pg_{N-1}.
\end{align*}
\revision{Define a linear SDP relaxation of \eqref{eqn::OBL-F-hij-wkdual-2} as 
\begin{equation}
\begin{alignedat}{4}
    &\minimize_{r_{i,j}}&&\minimize_{(\flam, \fbet, \tau') \geq \pmb{0}}\quad && {\tau'} 
    \\
    & &&\text{subject to} \quad && \fA(\flam, \fbet, \tau', (r_{i,j})) \succeq 0.
    \\
    & && &&\fB(\flam, \fbet, \tau') = \left(\flam, \fbet, \tau'\right)\geq 0
    \\
    & && &&\fC(\flam, \fbet ) = (-\lambda_1 + \beta_0, \lambda_1 - \lambda_2 + \beta_1, \dots, \lambda_{N-1} - \lambda_N + \beta_{N-1}, -1 + \lambda_N + \beta_N) = 0.
\end{alignedat}
\label{eqn::OBL-F-relaxation}
\end{equation}
\citep[Theorem 3]{drori2014performance} indicates that
 if we prove the choice in the previous claim satisfies KKT condition of \eqref{eqn::OBL-F-relaxation}, then this is also an optimal solution for the original problem since $(r_{i,j})$ uniquely determines $h_{i,j}$}.
The Lagrangian of the minimization problem is 
\begin{align*}
    \cL(\flam, \fbet, \tau',(r_{i,j}), \fK, \fb, \fc) = \frac{1}{2}\tau' - \tr\left\{\fA(\flam, \fbet, \tau',(r_{i,j})) \fK\right\}- \fb^\intercal \fB(\flam, \fbet, \tau') - \fc^\intercal  \fC(\flam, \fbet) 
\end{align*}
and the KKT conditions of the minimization problems are 
    \begin{align*}
        &\fA(\flam, \fbet,  \tau';(r_{i,j}))\succeq 0, \fB(\flam, \fbet,  \tau') \geq 0, \fC(\flam, \fbet) =0,
        \\
        &\nabla_{(\flam, \fbet, \tau',(r_{i,j}))}\cL(\flam, \fbet,  \tau',(r_{i,j}), \fK, \fb, \fc)  = 0,
        \\
        &\fK \succeq 0, \fb \geq 0,
        \\
        &\tr\left\{\fA(\flam, \fbet, \tau',(r_{i,j})) \fK\right\}=0,  \fb^\intercal \fB(\flam, \fbet,  \tau')=0,
    \end{align*}
where $\fK$ is a symmetric matrix. Here, $\fb = (\fu, \fv, s)$. We re-index $K$'s column and row starting from -1 (so $K$'s rows and columns index are $\{-1,0,1. \dots, N\}$). Now, we will show that there exist a dual optimal solution $(\fK, \fb, \fc)$ that $(\flam, \fbet, \tau',(r_{i,j}), \fK, \fb, \fc)$ satisfies KKT condition, which proves a pair $(\flam, \fbet, \tau', (r_{i,j}))$ is an optimal solution for primal problem. The stationary condition $\nabla_{(\flam, \fbet, \tau',(r_{i,j}))}\cL(\flam, \fbet,  \tau',(r_{i,j}), \fK, \fb, \fc)  = 0$ can be rewritten as 
\begin{equation}
    \begin{aligned}
    &\frac{\partial \cL}{\partial \lambda_k} = -\frac{1}{2}\left(K_{k-1, k-1} - K_{k-1, k} - K_{k, k-1} + K_{k,k} \right) - u_k +c_{k-1} - c_k = 0, \qquad  k \in \{ 1,2, \dots, N\}
    \\
    &\frac{\partial \cL}{\partial \beta_k} = \frac{L}{2}\left(K_{-1, k}  + K_{k, -1}\right) - v_k - c_k = 0, \qquad k \in \{0,1, \dots, N\}
    \\
    &\frac{\partial \cL}{\partial \tau'} = \frac{1}{2} -\frac{1}{2}K_{-1, -1} - s =0
    \\
    &\frac{\partial \cL}{\partial r_{k,t}} = -\frac{1}{2}(K_{k,t} + K_{t,k}) = 0,\qquad k \in \{ 1,2, \dots, N\}, \quad t \in \{0,1,\dots, k-1\}.
    \end{aligned}\label{eqn::OBL-F-KKT-stationary}
\end{equation}
We already know that $\fB(\flam, \fbet, \tau') \neq 0$, we can set $\fb = 0$. Then, \eqref{eqn::OBL-F-KKT-stationary} reduces to 
    \begin{align*}
    &K_{k,t} = 0, \qquad k \in \{1,2, \dots, N, t = 0,1,\dots, k-1\}
    \\
    &-\frac{1}{2}(K_{k-1, k-1} + K_{k,k})  +c_{k-1} - c_k  = 0, \qquad k \in \{1,2, \dots, N\}
    \\
    &L K_{-1, k} - c_k = 0, \qquad k \in \{0, 1,\dots, N\}
    \\
    &K_{-1, -1}  =1.
    \end{align*}
Then, we have
\begin{equation*}
\begin{aligned}
\fK &= \left(
\begin{array}{cccccc}
    1& \frac{c_0}{L} & \frac{c_1}{L} & \dots & \frac{c_{N-1}}{L} &\frac{c_N}{L} \\
    \frac{c_0}{L} & K_{0,0} & 0 & \dots&0 & 0\\
    \vdots & \vdots & \vdots & \ddots & \vdots &\vdots  \\
    \frac{c_{N-1}}{L} & 0 & 0 & \dots & K_{N-1, N-1} & 0\\
    \frac{c_{N}}{L} & 0& 0 & \dots &0 & K_{N,N}\\
\end{array}
\right) \succeq 0
\end{aligned}
\end{equation*} and since $\tr\left\{\fA(\flam, \fbet, \tau',(r_{i,j})) \fK\right\}=0$ with  $\fA \succeq 0$, we can replace this condition by $\fA(\flam, \fbet, \tau',(r_{i,j})) \fK=0$.
Then the KKT condition for the given $(\flam, \fbet, \tau', (r_{i,j}))$ reduces to 
\begin{equation*}
    \begin{aligned}
        &\frac{1}{2}\tau' - \frac{1}{2}\fbet^\intercal\fc = 0
        \\
        &\frac{1}{2L}\tau'\fc - \frac{L}{2}\diag(K_{0,0}, \dots, K_{N-1, N-1}, K_{N,N}) \fbet = 0 
        \\
        &-\fbet \fc^\intercal + \left(\begin{array}{cc}Q & \fq\\ \fq^\intercal & \frac{1}{2}\lambda_N  \end{array}\right)\diag(K_{0,0}, \dots, K_{N-1, N-1}, K_{N,N}) = 0.
    \end{aligned}
\end{equation*}
This clearly indicates that $K_{N,N} >0$ and $c_i$ are all determined as positive, and $K \succeq 0$. \end{proof}

\begin{claim}
The obtained algorithm is OBL-F$_\flat$. 
\end{claim}

\begin{proof}
By calculating $(h_{i,j})$ of OBL-F$_\flat$, we can prove the equivalence of the obtained solution and OBL-F$_\flat$. Indeed, OBL-F$_\flat$ is obtained by using \citep{lee2021geometric}'s auxiliary sequences. Except the last step of OBL-F$_\flat$, we will show that obtained $(\hat{h}_{i,j})$ satisfies \begin{align*}
    x_0 - \sum_{i=1}^{k+1} \sum_{j=0}^{i-1} \frac{\hat{h}_{i.j}}{L} \nabla f(x_j) &=  \left(1- \frac{2}{k+3} \right)\left( x_0 - \sum_{i=1}^k \sum_{j=0}^{i-1} \frac{\hat{h}_{i,j}}{L} \nabla f(x_j) - \frac{1}{L} \nabla f(x_k)\right)
    \\&+ \frac{2}{k+3}\left(x_0 - \sum_{j=0}^k \frac{j+1}{L} \nabla f(x_j)\right),
\end{align*}
which is re-written form of OBL-F$_\flat$. Comparing $\nabla f(x_j)$'s each coefficient, we should prove 
\begin{align*}
    \sum_{i = j+1}^{k+1} \hat{h}_{i,j} &=  \left(1- \frac{2}{k+3}\right)\sum_{i = j+1}^{k} \hat{h}_{i,j} + \frac{2}{k+3} (j+1) \qquad j \in \{0,1,\dots, k-1\}
    \\
    \hat{h}_{k+1,k} &= \left(1- \frac{2}{k+3}\right) + \frac{2}{k+3} (k+1),
\end{align*}
which is exactly equal to the recursive rule of \citep[Theorem 3]{drori2014performance}. The last step of OBL-F$_\flat$ can be analyzed similarly.
\end{proof}
In sum, the algorithm's performance criterion $f(x_N) - f_\star$ is bounded as 
$$f(x_{N}) - f_\star \leq  \frac{L}{N(N+1) + \sqrt{2N(N+1)}}\norm{x_0 - x_\star}^2.$$

Overall, we showed that OBL-F$_\flat$ is the ``best'' algorithm under $\cI_{\operatorname{OBL-F}_\flat}$. 
$$\operatorname{OBL-F}_\flat = \cA^\star_N(f(x_{N})- f_\star, \norm{x_0 - x_\star} \leq R, \cI_{\operatorname{OBL-F}_\flat})$$ and
\begin{align*}
    \cR(\operatorname{OBL-F}_\flat, f(x_{N})- f_\star,& \norm{x_0 - x_\star} \leq R, \cI_{\operatorname{OBL-F}_\flat}) \nonumber
    \\
    &= \cR^*(\kA_N, f(x_{N})- f_\star, \norm{x_0 - x_\star} \leq R, \cI_{\operatorname{OBL-F}_\flat}) 
    \\& = \frac{LR^2}{N(N+1) + \sqrt{2N(N+1)}} \nonumber
\end{align*}
hold.

\subsection{Proof of $\mathcal{A}^\star$-optimality of FGM}
To obtain FGM as an $\mathcal{A}^\star$-optimal algorithm, set $f(y_{N+1}) - f_\star$ to be the performance measure and $\norm{x_0 - x_\star} \leq R$ to be the initial condition.
Since the constraints and the objective of the problem \revision{are homogenous}, we assume $R = 1$ without loss of generality. For the argument of homogeneous, we refer \revision{to} \citep{drori2014performance, kim2016optimized,taylor2017smooth}. We use the set of inequalities that are handy for randomized coordinate updates and backtracking linesearches:
\begin{align*}
  \cI_{\operatorname{FGM}} = &\biggl\{  f_{k,0} \geq f_{k+1,1} + \frac{1}{2L} \norm{g_k}^2 \biggr\}_{k=0}^N  \bigcup \biggl\{f_{k,1} \geq f_{k,0} + \langle g_k, y_{k} - x_{k} \rangle\biggr\}_{k=1}^{N} 
  \\
  &\qquad   \bigcup \biggl\{f_\star \geq f_{k,0} + \langle g_k, x_\star - x_k \rangle \biggr\}_{k=0}^N.
\end{align*}
For calculating $\cR(\cA_N, \cP, \cC, \cI_{\operatorname{FGM}}) $ with fixed $\cA_N$, define the PEP with $\cI_{\operatorname{FGM}}$ as
\begin{equation*}
\cR(\cA_N, \cP, \cC, \cI_{\operatorname{FGM}}) = \left(
\begin{array}{lllll}
    &\maximize \quad & f_{N+1,1}&  - f_\star
    \\
    &\text{subject to} \quad &1 \quad \geq& \norm{x_0 - x_\star}^2 
    \\
    & &f_{k,0} \geq& f_{k+1,1} + \frac{1}{2L} \norm{g_k}^2, \quad \qquad &k\in \{0,1, \dots, N\}
    \\
    & &f_{k,1} \geq& f_{k,0} + \langle g_k, y_{k} - x_{k} \rangle,  \qquad &k \in \{1, \dots, N\}
    \\
    & &f_\star \mkern11mu \geq& f_{k,0} + \langle g_k, x_\star - x_k \rangle , \qquad &k \in \{0,1, \dots, N\}
    \\
    & & \revision{x_k, y_k} & \revision{\text{are following the algorithm } \cA_N.}
\end{array}\right)
\end{equation*}
Using the notation of Section~\ref{ss::notation}, we reformulate \revision{the} problem of computing the risk $\cR(\cA_N, \cP, \cC, \cI_{\operatorname{FGM}})$ as the following SDP:
\begin{alignat*}{4}
&\maximize_{\fG, \fF_0, \fF_1} \quad  &&\pf_{N+1}^\intercal \fF_1
\\
&\mbox{subject to}\quad 
    &&1 \geq \px_0^\intercal \fG \px_0  
    \\
    & && 0 \geq \pf_{k+1}^\intercal\fF_1 - \pf_{k}^\intercal \fF_0 + \frac{1}{2L} \pg_k^\intercal \fG \pg_k, \qquad &&k \in \{0,1, \dots, N\}
    \\
    & && 0 \geq \pf_k^\intercal(\fF_0 - \fF_1) + \pg_k^\intercal \fG (\px_{k-1} - \px_k ) - \frac{1}{L} \pg_{k-1}^\intercal \fG \pg_k, \qquad &&k \in \{1,2, \dots, N\}
    \\
    & && 0 \geq \pf_k^\intercal \fF_0 - \pg_k^\intercal \fG \px_k, \qquad && k \in \{0,1, \dots, N\}.
    \\
    & && \revision{\fG \succeq 0, \fF_0 \geq 0, \fF_1 \geq 0.}
\end{alignat*}
For above transformation, $d\geq N+2 $ is used \citep{taylor2017smooth}. \revision{The Lagrangian} of the optimization problem becomes 
\begin{align*}
    \Lambda(\fF_0, \fF_1&, \fG, \flam, \fbet, \falp, \tau)
    \\
    &=  - \pmb{f}_{N+1}^\intercal \fF_1 + \tau (\px_0^\intercal \fG \px_0 -1 ) + \sum_{k=0}^{N} \alpha_{k}\left( \pf_{k+1}^\intercal\fF_1 - \pf_{k}^\intercal \fF_0 + \frac{1}{2L} \pg_k^\intercal \fG \pg_k
    \right)
    \\&+ \sum_{k=1}^N \lambda_k \left(\pf_k^\intercal(\fF_0 - \fF_1) + \pg_k^\intercal \fG (\px_{k-1} - \px_k ) - \frac{1}{L} \pg_{k-1}^\intercal \fG \pg_k\right) 
    \\
    & +\sum_{k=0}^{N} \beta_k \left( \pf_k^\intercal \fF_0 - \pg_k^\intercal \fG \px_k\right) 
\end{align*}
with dual variables $\flam = (\lambda_1, \dots, \lambda_N) \in \bR_{+}^{N}$, $\fbet = (\beta_0, \dots, \beta_N)\in \bR_{+}^{N+1}$, $\falp = (\alp_0, \dots, \alp_N)\in \bR_{+}^{N+1}$, and $\tau \geq 0$. 
Then the dual formulation of PEP problem is
\begin{equation}
\begin{alignedat}{2}
    &\maximize_{(\flam, \fbet, \falp, \tau) \geq \pmb{0}}\quad  &&{-\tau} 
    \\
     &\mbox{subject to} &&\pmb{0} = -\sum_{k=0}^{N} \alpha_k \pmb{f}_k + \sum_{k=1}^N \lambda_k \pmb{f}_k +\sum_{k=0}^{N} \beta_k \pmb{f}_k
     \\
     & &&\pmb{0} = -\pmb{f}_{N+1} - \sum_{k=1}^N \lambda_k \pmb{f}_k + \sum_{k=0}^{N} \alpha_k \pmb{f}_{k+1}
     \\
     & &&0 \preceq S(\flam, \fbet, \falp, \tau), 
\end{alignedat} \label{eqn::FGM-dualPEP}    
\end{equation}
where $S$ is defined as  
\begin{align*}
    S(\flam, \fbet, \falp, \tau) &=  \tau \pmb{x}_0\pmb{x}_0^\intercal+ \sum_{k=0}^{N} {\alpha_k} \left( \frac{1}{2L} \pmb{g}_{k}\pmb{g}_{k}^\intercal\right) + \sum_{k=0}^{N}\frac{\beta_k}{2}\left(-\pmb{g}_k \pmb{x}_k^\intercal-\pmb{x}_k \pmb{g}_k^\intercal  \right)
    \\
    &+ \sum_{k=1}^N \frac{\lambda_k}{2} \left(\pmb{g}_k (\pmb{x}_{k-1} - \pmb{x}_{k})^\intercal +  (\pmb{x}_{k-1} - \pmb{x}_{k})\pmb{g}_k^\intercal - \frac{1}{L} \pmb{g}_{k-1}\pmb{g}_{k}^\intercal - \frac{1}{L} \pmb{g}_{k}\pmb{g}_{k-1}^\intercal\right).
\end{align*}
We have a strong duality argument
\begin{align*}
    \argmin_{h_{i,j}}\maximize_{\fG, \fF_0, \fF_1} \mkern7mu \pf_{N+1}^\intercal \fF_1  = \argmin_{h_{i,j}}\minimize_{(\flam, \fbet, \falp,\tau) \geq \pmb{0}}\mkern7mu {\tau}, 
\end{align*}
 as ORC-F's optimality proof. Remind that \eqref{eqn::FGM-dualPEP} finds the ``best'' proof for the algorithm. Now we investigate the optimization step for algorithm. The last part is minimizing  \eqref{eqn::FGM-dualPEP} with stepsize, i.e. 
 \begin{alignat}{3}
    &\minimize_{h_{i,j}}&&\minimize_{(\flam, \fbet, \falp,\tau) \geq \pmb{0}}\quad && {\tau}  \label{eqn::FGM-hij-wkdual}
    \\
    & && \text{subject to} \quad && \pmb{0}  = -\sum_{k=0}^{N} \alpha_k \pmb{f}_k + \sum_{k=1}^N \lambda_k \pmb{f}_k +\sum_{k=0}^{N} \beta_k \pmb{f}_k \label{eqn::FGM-const1}
     \\
     & && &&\pmb{0}  = -\pmb{f}_{N+1} - \sum_{k=1}^N \lambda_k \pmb{f}_k + \sum_{k=0}^{N} \alpha_k \pmb{f}_{k+1} \label{eqn::FGM-const2}
     \\
     & && &&0 \preceq S(\flam, \fbet, \falp, \tau). \label{eqn::FGM-Smatrix}
\end{alignat}
We note that $\pf_i$ is \revision{a} standard unit vector mentioned in \eqref{eqn::egf} (not a variable), we can write \eqref{eqn::FGM-const1} and \eqref{eqn::FGM-const2} as 
\begin{equation}
\begin{alignedat}{2}
&&\left(\begin{array}{ll}
\beta_k = \alpha_k - \lambda_k  = \lambda_{k+1} - \lambda_k, \qquad & k \in \{1, \dots, N-1\}
    \\
    \beta_0 = \alpha_0 = \lambda_1
    \\
    \beta_N = \alpha_N - \lambda_N = 1 - \lambda_N. 
\end{array}\right) 
\\
 &&\left(\begin{array}{ll}
\alpha_N = 1 
    \\
    \alpha_k = \lambda_{k+1}, \quad  & k \in \{0,1, \dots, N-1\}
\end{array}\right) 
\end{alignedat}\label{eqn::FGM-constlinear}
\end{equation}
We consider \eqref{eqn::FGM-Smatrix} with \eqref{eqn::FGM-constlinear} and FSFO's $h_{i,j}$. To be specific, we substitute $\falp$ and $\fbet$ to $\flam$ in $S(\flam, \fbet, \falp, \tau)$. To show the dependency of $S$ to $(h_{i,j})$ since $\px_k$ are represented with $(h_{i,j})$, we will explicitly write $S$ as $S(\flam, \tau;(h_{i,j}))$. Then, we get
\begin{align*}
    S(\flam,  \tau;(h_{i,j})) &=\tau \pmb{x}_0\pmb{x}_0^\intercal -\sum_{k=1}^N \frac{\lambda_k}{2L}\pmb{g}_k\pmb{g}_k^\intercal  + \frac{1}{2L} \pmb{g}_{N}\pmb{g}_{N}^\intercal + \sum_{k=1}^N\frac{\lambda_k}{2L} (\pmb{g}_{k-1}-\pmb{g}_{k})(\pmb{g}_{k-1}-\pmb{g}_{k})^\intercal 
    \\&+ \sum_{k=1}^{N-1} \sum_{t=0}^{k-1} \left( \frac{\lambda_k}{2} \frac{h_{k,t}}{L} + \frac{\lambda_{k+1} - \lambda_{k}}{2}\sum_{j=t+1}^{k}\frac{h_{j,t}}{L} \right) \left(\pmb{g}_k\pmb{g}_t^\intercal + \pmb{g}_t\pmb{g}_k^\intercal\right) 
    \\&+ \sum_{t=0}^{N-1} \left( \frac{\lambda_N}{2} \frac{h_{N,t}}{L} + \frac{1 - \lambda_{N}}{2}\sum_{j=t+1}^{N}\frac{h_{j,t}}{L} \right) \left(\pmb{g}_N\pmb{g}_t^\intercal + \pmb{g}_t\pmb{g}_N^\intercal\right)
    \\&- \sum_{k=1}^{N-1} \frac{\lambda_{k+1} - \lambda_{k}}{2} \left(\pmb{x}_0\pmb{g}_k^\intercal + \pmb{g}_k\pmb{x}_0^\intercal\right) - \frac{\lambda_1}{2}\left(\pmb{x}_0\pmb{g}_0^\intercal + \pmb{g}_0\pmb{x}_0^\intercal\right) - \frac{1 - \lambda_N}{2}\left(\pmb{x}_0\pmb{g}_N^\intercal + \pmb{g}_N\pmb{x}_0^\intercal\right).
\end{align*}
Using the fact that $\px_0, \pg_i, \pf_i$ are unit vectors, we can represent $S(\flam, \tau;(h_{i,j}))$ with $\fgam(\flam) = -L\fbet = -L(\lambda_1, \lambda_2 - \lambda_1, \dots, 1-\lambda_N) = (\hat{\fgam}(\flam), \gamma_{N}(\flam)) $ and  $\tau' = 2L\tau$ as  
\begin{align*}
   S(\flam, \tau';(h_{i,j}))  &= \frac{1}{L}\left(
\begin{array}{ccc}
    \frac{1}{2}\tau' & \frac{1}{2}\hat{\fgam}(\flam)^\intercal & \frac{1}{2}\gamma_{N}(\flam) \\
    \frac{1}{2}\hat{\fgam}(\flam) & Q(\flam; (h_{i,j})) & q(\flam; (h_{i,j})) 
    \\
    \frac{1}{2}\gamma_{N}(\flam) & q(\flam; (h_{i,j}))^\intercal & \frac{1}{2} \\
\end{array}
\right) \succeq 0.
\end{align*}
\revision{Here, $Q$ and $\fq$ are defined as 
\begin{align*}
     Q(\flam;(h_{i,j})) =  &-\sum_{k=1}^{N-1} \frac{\lambda_k}{2}\pmb{g}'_k\pmb{g}_k^{'\intercal}  + \sum_{k=1}^{N-1}\frac{\lambda_k}{2} (\pmb{g}'_{k-1}-\pmb{g}'_{k})(\pmb{g}'_{k-1}-\pmb{g}'_{k})^\intercal + \frac{\lambda_N}{2}\pg'_{N-1}\pg_{N-1}^{'\intercal}
    \\&+ \sum_{k=1}^{N-1} \sum_{t=0}^{k-1} \left( \frac{\lambda_k}{2} h_{k,t} + \frac{\lambda_{k+1} - \lambda_{k}}{2}\sum_{j=t+1}^{k}h_{j,t} \right) \left(\pmb{g}'_k\pmb{g}_t^{'\intercal} + \pmb{g}'_t\pmb{g}_k^{'\intercal}\right)
\end{align*}
and 
\begin{align*}
    \fq(\flam;(h_{i,j})) &= -\frac{\lambda_N}{2}\pg'_{N-1} + \sum_{t=0}^{N-1} \left( \frac{\lambda_N}{2} h_{N,t} + \frac{1 - \lambda_{N}}{2}\sum_{j=t+1}^{N}h_{j,t} \right) \pmb{g}'_t
    \\
    &= \sum_{t=0}^{N-2} \left( \frac{\lambda_N}{2} h_{N,t} + \frac{1 - \lambda_{N}}{2}\sum_{j=t+1}^{N}h_{j,t} \right) \pmb{g}'_t +  \left(\frac{1}{2}{h_{N,N-1}} -\frac{\lambda_N}{2}\right)\pmb{g}'_{N-1}.
\end{align*}
where $\pg'_k = e_{k+1} \in \bR^{N+1}$.}
Note that \eqref{eqn::FGM-hij-wkdual} is equivalent to 
\begin{alignat}{3}
    &\minimize_{h_{i,j}}&&\minimize_{(\flam,  \tau') \geq \pmb{0}}\quad && {\tau'}  \label{eqn::FGM-hij-wkdual-2}
    \\
    & && \text{subject to} \quad && \left(
\begin{array}{ccc}
    \frac{1}{2}\tau' & \frac{1}{2}\hat{\fgam}(\flam)^\intercal & \frac{1}{2}\gamma_{N}(\flam) \\
    \frac{1}{2}\hat{\fgam}(\flam) & Q(\flam; (h_{i,j})) & q(\flam; (h_{i,j})) 
    \\
    \frac{1}{2}\gamma_{N}(\flam) & q(\flam; (h_{i,j}))^\intercal & \frac{1}{2} \\
\end{array}
\right) \succeq 0 \nonumber
\end{alignat}
and dividing this optimized value with $2L$ gives the optimized value of \eqref{eqn::FGM-hij-wkdual}. Using Schur complement \citep{golub1996matrix}, \eqref{eqn::FGM-hij-wkdual-2} can be converted to the problem as 
\begin{alignat}{3}
    &\minimize_{h_{i,j}}&&\minimize_{(\flam, \tau') \geq \pmb{0}}\quad && {\tau'}  \label{eqn::FGM-hij-final-problem}
    \\
    & && \text{subject to} \quad && 
\left(
\begin{array}{cc}
    Q - 2\fq \fq^\intercal & \frac{1}{2}( \hat{\fgam}(\flam) - 2\fq \gam_N(\flam)) \\
    \frac{1}{2}( \hat{\fgam}(\flam) - 2\fq \gam_N(\flam))^\intercal & \frac{1}{2} (\tau' - {\gam_N(\flam)^2}) \\
\end{array}
\right)
\succeq 0. \label{eqn::FGM-final-matrix}
\end{alignat}
So far we simplified SDP. We will have \revision{three steps}: finding variables that \revision{make} \eqref{eqn::FGM-final-matrix}'s left hand side zero, showing that the solution from \revision{the} first step satisfies KKT condition, and finally showing that obtained algorithm is equivalent to FGM.

\begin{claim}
\label{claim::7}
    There is a point that makes \eqref{eqn::FGM-final-matrix}'s left-hand side zero. 
\end{claim}
\begin{proof}
\revision{Defining $\{r_{k,t}\}_{k=1,2, \dots, N,  t=0, \dots, k-1}$ as $$r_{k,t} = \lambda_k h_{k,t}- \frac{\gamma_k}{L}\sum_{j=t+1}^k h_{j,t}.$$
Then, if $r_{i,j}$ is determined, \citep[Theorem 5.1]{drori2014performance} indicates this uniquely determine $h_{i,j}$. We set $(\lambda_k)_{k=0}^N$ and $(r_{N,k})_{k=0}^{N-1}$ as 
\begin{equation}
\begin{aligned}
    &\lambda_{k} = \frac{\theta_{k-1}^2}{\theta_N^2}, \qquad k \in \{ 1,2, \dots, N\}
    \\
    &r_{N,k} = \frac{\theta_k}{\theta_N}, \qquad k \in \{0,1,\dots, N-2 \}
    \\
    &r_{N, N-1} - \lambda_N = \frac{\theta_{N-1}}{\theta_N}.
\end{aligned}     \label{eqn::FGM-solution1}
\end{equation}
Moreover, we set 
\begin{equation}
\begin{aligned}
 &r_{k,t} = \frac{\theta_{k-1}\theta_{t-1}}{\theta_N^2}, \qquad k \in \{1,2, \dots, N-1\}, \quad t\in \{0,1, \dots ,k-2\}
 \\
 &r_{k,k-1} = \frac{\theta_{k-1}\theta_{k-2}}{\theta_N^2} + \frac{\theta_{k-1}^2}{\theta_N^2}, \qquad k \in \{ 1,2 \dots ,N-1\}.
\end{aligned}
 \label{eqn::FGM-solution2}
\end{equation}
In addition, we set $\hat{\fgam}$ as 
\begin{alignat*}{3}
    &\gamma_t &&= \gamma_N r_{N, t}, \qquad && t \in \{0,1, \dots, N-2 \}
    \\
    &\gamma_{N-1} &&= \gamma_N(r_{N, N-1} - \lambda_N)
    \\
    &\gamma_N = L(1-\lambda_N).
\end{alignat*}
Lastly, we set $\tau'$ as 
\begin{align}
    \tau' = \frac{L^2}{\theta_N^2},
     \label{eqn::FGM-solution3}
\end{align}
and $\tau = \frac{L}{2\theta_N^2}$. These variables make \eqref{eqn::FGM-final-matrix}'s left-hand side zero. }
\end{proof}
\begin{claim}
\eqref{eqn::FGM-solution1}, \eqref{eqn::FGM-solution2} and \eqref{eqn::FGM-solution3} are an optimal solution of \eqref{eqn::FGM-hij-final-problem}.
\end{claim}
\begin{proof}
Let we represent $S$ with the variable $(r_{i,j})$. We will denote this as $\fA$. To be specific, 
\begin{align*}
   \fA(\flam, \fbet, \falp, \tau',(r_{i,j})) = S(\flam, \fgam, \tau'; (h_{i,j})) = \left(
\begin{array}{ccc}
    \frac{1}{2}\tau' & -\frac{L}{2}\hat{\fbet}^\intercal & -\frac{L}{2}\bet_N \\
    -\frac{L}{2}\hat{\fbet}^\intercal & Q(\flam; (r_{i,j})) & \fq((r_{i,j})) 
    \\
    -\frac{L}{2}\bet_N & \fq((r_{i,j}))^\intercal & \frac{1}{2} \\
\end{array}
\right)  \succeq 0 .
\end{align*}
Here, $\fbet = (\hat{\fbet}^\intercal, \beta_N)^\intercal$, \begin{align*}
     Q(\flam;\revision{(r_{i,j})}) &-\sum_{k=1}^{N-1} \frac{\lambda_k}{2}\pmb{g}_k\pmb{g}_k^\intercal  + \sum_{k=1}^{N-1}\frac{\lambda_k}{2} (\pmb{g}_{k-1}-\pmb{g}_{k})(\pmb{g}_{k-1}-\pmb{g}_{k})^\intercal + \frac{\lambda_N}{2}\pg_{N-1}\pg_{N-1}^\intercal+ \sum_{k=1}^{N-1} \sum_{t=0}^{k-1} \frac{r_{k,t}}{2} \left(\pmb{g}_k\pmb{g}_t^\intercal + \pmb{g}_t\pmb{g}_k^\intercal\right)
\end{align*}
and 
\begin{align*}
    \fq((r_{i,j})) &= \sum_{t=0}^{N-1} \frac{r_{N,t}}{2}\pmb{g}_t - \frac{\lambda_N}{2} \pg_{N-1}.
\end{align*}

\revision{Define a linear SDP relaxation of \eqref{eqn::FGM-hij-wkdual-2} as 
\begin{equation}
\begin{alignedat}{4}
    &\minimize_{r_{i,j}}&&\minimize_{(\flam, \fbet, \falp, \tau') \geq \pmb{0}}\quad && {\tau'} 
    \\
    & &&\text{subject to} \quad && \fA(\flam, \fbet, \falp, \tau', (r_{i,j})) \succeq 0.
    \\
    & && &&\fB(\flam, \fbet, \falp, \tau') = \left(\flam, \fbet, \falp, \tau'\right) \geq 0
    \\
    & && &&\fC(\flam, \fbet, \falp) = (-\alpha_0 + \beta_0, -\alpha_1 + \lambda_1 + \beta_1, \dots, -\alpha_N + \lambda_N + \beta_N) = 0
    \\
    & && &&\fD(\flam, \fbet, \falp) = (-\lambda_1 + \alpha_0, -\lambda_2 + \alpha_1, \dots, -\lambda_N + \alpha_{N-1}, \alpha_N -1) = 0.
\end{alignedat}
\label{eqn::FGM-relaxation}
\end{equation}
\citep[Theorem 3]{drori2014performance} indicates that
 if we prove the choice in the previous claim satisfies KKT condition of \eqref{eqn::FGM-relaxation}, then this is also an optimal solution for the original problem.} \revision{The Lagrangian} of the minimization problem is 
\begin{align*}
    \cL(\flam, &\fbet, \falp, \tau',(r_{i,j}), \fK, \fb, \fc, \fd) 
    \\
    &= \frac{1}{2}\tau' - \tr\left\{\fA(\flam, \fbet, \falp, \tau',(r_{i,j})) \fK\right\}- \fb^\intercal \fB(\flam, \fbet, \falp, \tau') - \fc^\intercal  \fC(\flam, \fbet, \falp) - \fd^\intercal\fD(\flam, \fbet, \falp)
\end{align*}
and the KKT conditions of the minimization problems are 
    \begin{align*}
        &\fA(\flam, \fbet, \falp, \tau';(r_{i,j}))\succeq 0, \fB(\flam, \fbet, \falp, \tau') \geq 0, \fC(\flam, \fbet, \falp) =0, \fD(\flam, \fbet, \falp) = 0,
        \\
        &\nabla_{(\flam, \fbet, \falp, \tau',(r_{i,j}))}\cL(\flam, \fbet, \falp, \tau',(r_{i,j}), \fK, \fb, \fc, \fd)  = 0,
        \\
        &\fK \succeq 0, \fb \geq 0,
        \\
        &\tr\left\{\fA(\flam, \fbet, \falp, \tau',(r_{i,j})) \fK\right\}=0,  \fb^\intercal \fB(\flam, \fbet, \falp, \tau')=0,
    \end{align*}
where $\fK$ is a symmetric matrix. Here, $\fb = (\fu, \fv, \fw, s)$. We re-index $K$'s column and row starting from -1 (so $K$'s rows and columns index are $\{-1,0,1. \dots, N\}$). Now, we will show that there exist a dual optimal solution $(\fK, \fb, \fc, \fd)$ that $(\flam, \fbet, \falp, \tau',(r_{i,j}), \fK, \fb, \fc, \fd)$ satisfies KKT condition, which proves a pair $(\flam, \fbet, \falp, \tau', (r_{i,j}))$ is an optimal solution for primal problem. The stationary condition $\nabla_{(\flam, \fbet, \falp, \tau',(r_{i,j}))}\cL(\flam, \fbet, \falp, \tau',(r_{i,j}), \fK, \fb, \fc, \fd)  = 0$ can be rewritten as
\begin{equation}
    \begin{alignedat}{4}
    &\frac{\partial \cL}{\partial \lambda_k} &&= -\frac{1}{2}\left(K_{k-1, k-1} - K_{k-1, k} - K_{k, k-1}\right) - u_k - c_k +d_{k-1} = 0, \qquad k \in \{1,2, \dots, N\}
    \\
    &\frac{\partial \cL}{\partial \beta_k} &&= \frac{L}{2}\left(K_{-1, k}  + K_{k, -1}\right) - v_k - c_k = 0, \qquad k \in \{0,1, \dots, N\}
    \\
    &\frac{\partial \cL}{\partial \alpha_k} &&= -w_k + c_k - d_k = 0, \qquad k \in \{0, 1, \dots, N\}
    \\
    &\frac{\partial \cL}{\partial \tau'} &&= \frac{1}{2} -\frac{1}{2}K_{-1, -1} - s =0
    \\
    &\frac{\partial \cL}{\partial r_{k,t}} &&= -\frac{1}{2}(K_{k,t} + K_{t,k}) = 0, \qquad k \in \{ 1,2, \dots, N\}, \quad t \in \{0,1,\dots, k-1\}.
    \end{alignedat}\label{eqn::FGM-KKT-stationary}
\end{equation}
We already know that $\fB(\flam, \fbet, \falp, \tau') \neq 0$, we can set $\fb = 0$. Then, \eqref{eqn::FGM-KKT-stationary} reduces to 
    \begin{alignat*}{4}
    &K_{k,t} = 0, \qquad k \in \{ 1,2, \dots, N\}, \quad t \in \{0,1,\dots, k-1\}
    \\
    &-\frac{1}{2}{K_{k-1, k-1}}  - c_k + d_{k-1}= 0, \qquad k \in \{ 1,2, \dots, N\}
    \\
    &L K_{-1, k} - c_k = 0, \qquad k \in \{ 0, 1,\dots, N\}
    \\
    & c_k - d_k= 0,\qquad  k \in \{ 0,1, \dots, N    \}
     \\
    &K_{-1, -1} =1.
    \end{alignat*}
Then, we have
\begin{equation*}
\begin{aligned}
\fK &= \left(
\begin{array}{cccccc}
    1& \frac{c_0}{L} & \frac{c_1}{L} & \dots & \frac{c_{N-1}}{L} &\frac{c_N}{L} \\
    \frac{c_0}{L} & 2c_0 - 2c_1 & 0 & \dots&0 & 0\\
    \vdots & \vdots & \vdots & \ddots & \vdots &\vdots  \\
    \frac{c_{N-1}}{L} & 0 & 0 & \dots & 2c_{N-1} - 2c_N & 0\\
    \frac{c_{N}}{L} & 0& 0 & \dots &0 & K_{N,N}\\
\end{array}
\right) \succeq 0
\end{aligned}
\end{equation*} and since $\tr\left\{\fA(\flam, \fbet, \falp, \tau',(r_{i,j})) \fK\right\}=0$ with  $\fA \succeq 0$, we can replace this condition by $\fA(\flam, \fbet, \falp, \tau',(r_{i,j})) \fK=0$.
Then the KKT condition for the given $(\flam, \fbet, \falp, \tau', (r_{i,j}))$ reduces to 
\begin{equation*}
    \begin{aligned}
        &\frac{1}{2}\tau' - \frac{1}{2}\fbet^\intercal\fc = 0
        \\
        &\frac{1}{2L}\tau'\fc - \frac{L}{2}\diag(2c_0 - 2c_1, \dots, 2c_{N-1} - 2c_N, K_{N,N}) \fbet = 0 
        \\
        &-\fbet \fc^\intercal + \left(\begin{array}{cc}Q & \fq\\ \fq^\intercal & \frac{1}{2}  \end{array}\right)\diag(2c_0 - 2c_1, \dots, 2c_{N-1} - 2c_N, K_{N,N}) = 0
    \end{aligned}
\end{equation*}
and this is equivalent to 
    \begin{align*}
        &\sum_{i=0}^N \theta_i c_i = L^2
        \\
        &c_i = (2c_i - 2c_{i+1}) \theta_i \qquad \text{for } i = 0, 1, \dots, N-1 
        \\
        &c_N = K_{N,N} \theta_N.
    \end{align*}
This clearly indicates that $K_{N,N} >0$ and $c_i$ are all determined as positive, and 
\begin{equation*}
\begin{aligned}
\fK &= \left(
\begin{array}{cccccc}
    1& \frac{c_0}{L} & \frac{c_1}{L} & \dots & \frac{c_{N-1}}{L} &\frac{c_N}{L} \\
    \frac{c_0}{L} & \frac{c_0}{\theta_0} & 0 & \dots&0 & 0\\
    \vdots & \vdots & \vdots & \ddots & \vdots &\vdots  \\
    \frac{c_{N-1}}{L} & 0 & 0 & \dots & \frac{c_{N-1}}{\theta_{N-1}} & 0\\
    \frac{c_{N}}{L} & 0& 0 & \dots &0 & \frac{c_N}{\theta_N}\\
\end{array}
\right) \succeq 0
\end{aligned}
\end{equation*}
since $\sum_{i=0}^N \theta_i c_i = L^2$. 
\end{proof}
\begin{claim}
The obtained algorithm is FGM. 
\end{claim}

\begin{proof}
By calculating $(h_{i,j})$ of FGM, we can prove the equivalence of the obtained solution and FGM. Indeed, FGM is obtained by using \citep{lee2021geometric}'s auxiliary sequences.
We will show that obtained $(\hat{h}_{i,j})$ satisfies \begin{align*}
    x_0 - \sum_{i=1}^{k+1} \sum_{j=0}^{i-1} \frac{\hat{h}_{i.j}}{L} \nabla f(x_j) &= \left(1 - \frac{1}{\theta_{k+1}}\right)  \left( x_0 - \sum_{i=1}^k \sum_{j=0}^{i-1} \frac{\hat{h}_{i,j}}{L} \nabla f(x_j) - \frac{1}{L} \nabla f(x_k)\right)
    \\&+ \frac{1}{\theta_{k+1}}\left(x_0 - \sum_{j=0}^k \frac{\theta_j}{L} \nabla f(x_j)\right),
\end{align*}
which is re-written form of FGM. Comparing $\nabla f(x_j)$'s each coefficient, we should prove 
\begin{align*}
    \sum_{i = j+1}^{k+1} \hat{h}_{i,j} &= \left(1 - \frac{1}{\theta_{k+1}}\right) \sum_{i = j+1}^{k} \hat{h}_{i,j} + \frac{1}{\theta_{k+1}}\theta_j \qquad j \in \{0,1,\dots, k-1\}
    \\
    \hat{h}_{k+1,k} &= \left(1 - \frac{1}{\theta_{k+1}}\right) + \frac{1}{\theta_{k+1}} \theta_k,
\end{align*}
which is exactly equal to the recursive rule of \citep[Theorem 3]{drori2014performance}. 
\end{proof}

In sum, the algorithm's performance criterion $f(y_{N+1}) - f_\star$ is bounded as 
$$f(y_{N+1}) - f_\star \leq \frac{L}{2\theta_N^2}\norm{x_0 - x_\star}^2.$$

Overall, we showed that FGM is the ``best'' algorithm under $\cI_{\operatorname{FGM}}$. 
$$\operatorname{FGM} = \cA^\star_N(f(y_{N+1})- f_\star, \norm{x_0 - x_\star} \leq R, \cI_{\operatorname{FGM}})$$ and
\begin{align*}
    \cR(\operatorname{FGM}, f(y_{N+1})- f_\star,& \norm{x_0 - x_\star} \leq R, \cI_{\operatorname{FGM}}) \nonumber
    \\
    &= \cR^*(f(y_{N+1})- f_\star, \norm{x_0 - x_\star} \leq R, \cI_{\operatorname{FGM}}) 
    \\& = \frac{LR^2}{2\theta_N^2}\nonumber
\end{align*}
hold. 
\subsection{Conjecture of $\mathcal{A}^\star$-optimality of OBL-G$_\flat$}
We give a conjecture for $\cA^\star$-optimality of OBL-G$_\flat$. Set $\norm{\nabla f(x_N)}^2$ to be the performance measure and $f(x_0) - f_\star \leq \frac{1}{2}LR^2 $ to be the initial condition.  We use the set of inequalities that are handy for backtracking linesearches. Since the constraints and the objective of the problem \revision{are homogenous}, we assume $R = 1$ without loss of generality. For the argument of homogeneous, we refer \revision{to}\citep{drori2014performance, kim2016optimized,taylor2017smooth}. We use the set of inequalities that are handy for backtracking linesearches: 
\begin{align*}
  \cI_{\operatorname{OBL-G}_\flat} = &\biggl\{f_{k-1,0} \geq f_{k,0} + \langle g_k, x_{k-1} - x_{k} \rangle + \frac{1}{2L}\norm{g_{k-1,0} - g_k}^2\biggr\}_{k=1}^{N} 
  \\
  &\qquad   \bigcup \biggl\{f_{N,0} \geq f_{k,0} + \langle g_k, x_N - x_k \rangle \biggr\}_{k=0}^{N-1} \bigcup \biggl\{ f_{N,0} \geq f_\star + \frac{1}{2}\norm{g_{N, 0}}^2\biggr\}.
\end{align*}
For calculating $\cR(\cA_N, \cP, \cC, \cI_{\operatorname{OBL-G}_\flat}) $ with fixed $\cA_N$, define the PEP with $\cI_{\operatorname{OBL-G}_\flat}$ as
\begin{align*}
&\cR(\cA_N, \cP, \cC, \cI_{\operatorname{OBL-G}_\flat}) 
\\ &= \left(
\begin{array}{lllll}
    &\maximize \quad & \norm{g_N}^2 & 
    \\
    &\text{subject to} \quad &0 \quad  &\geq  f_{0,0}  - \frac{1}{2}L 
    \\
    & &f_{k-1,0} &\geq  f_{k,0} + \langle g_k, x_{k-1} - x_{k} \rangle + \frac{1}{2L}\norm{g_{k-1} - g_{k}}^2,  k\in \{1, 2, \dots, N\}
    \\
    & &f_{N,0} &\geq f_{k,0} + \langle g_k, x_N - x_k \rangle ,  k \in \{0,1, \dots, N-1\}
    \\
    & &f_{N, 0}& \geq  f_\star + \frac{1}{2L}\norm{g_N}^2
    \\
    & & \revision{x_k} & \revision{\text{is following the algorithm } \cA_N.}
\end{array}\right)
\end{align*}
Using the notation of Section~\ref{ss::notation}, we reformulate \revision{the} problem of computing the risk $\cR(\cA_N, \cP, \cC, \cI_{\operatorname{OBL-G}_\flat})$ as the following SDP:
\begin{alignat*}{4}
&\maximize_{\fG, \fF_0} \quad  &&\pg_N^\intercal \fG \pg_N
\\
&\mbox{subject to}\quad 
    &&0 \geq \pf_0^\intercal F_0 - \frac{1}{2}L 
    \\
    & && 0 \geq (\pmb{f}_{k} - \pmb{f}_{k-1})^\intercal \fF_0 + \pmb{g}_{k}^\intercal \fG (\pmb{x}_{k-1} - \pmb{x}_{k}) + \frac{1}{2L} (\pmb{g}_{k-1} - \pmb{g}_{k})^\intercal \fG(\pmb{g}_{k-1} - \pmb{g}_{k}), \qquad &&k \in \{1,2, \dots, N\}
    \\
    & &&0 \geq (\pf_k^\intercal - \pf_N)^\intercal \fF_0  + \pg_k^\intercal\fG\px_N -\pg_k^\intercal \fG \px_k, \qquad &&k \in \{0, 1, \dots, N-1\}
    \\
    & && 0 \geq -\pf_N^\intercal \fF_0 + \frac{1}{2L}\pg_N^\intercal \fG \pg_N
    \\
    & && \revision{\fG \succeq 0, \fF_0 \geq 0}.
\end{alignat*}
For above transformation, $d\geq N+2 $ is used \citep{taylor2017smooth}. \revision{The Lagrangian} of the optimization problem becomes 
\begin{align*}
    \Lambda(\fF_0&, \fG, \flam, \fbet, \tau, c)
    \\
    &=  - \pg_N^\intercal \fG \pg_N + \tau\left(\pmb{f}_0^\intercal \fF_0 - \frac{1}{2}L\right)+c\left( -\pmb{f}_N^\intercal \fF_0 + \frac{1}{2L} \pmb{g}_N^\intercal \fG \pmb{g}_N\right) 
    \\& + \sum_{k=1}^{N} \lambda_k \left( (\pmb{f}_{k} - \pmb{f}_{k-1})^\intercal \fF_0 + \pmb{g}_{k}^\intercal \fG (\pmb{x}_{k-1} - \pmb{x}_{k}) + \frac{1}{2L} (\pmb{g}_{k-1} - \pmb{g}_{k})^\intercal \fG(\pmb{g}_{k-1} - \pmb{g}_{k})\right)
    \\
    &  +\sum_{k=0}^{N-1} \beta_k \left((\pmb{f}_k- \pmb{f}_N)^\intercal \fF_0 + \pmb{g}_k^\intercal \fG \pmb{x}_N - \pmb{g}_k^\intercal \fG \pmb{x}_k\right).
\end{align*}
with dual variables $\flam = (\lambda_1, \dots, \lambda_{N}) \in \bR_{+}^{N}$, $\fbet = (\beta_0, \dots, \beta_{N-1}) \in \bR_{+}^{N}$, and $\tau, c \geq 0$.

Then the dual formulation of PEP problem is 
\begin{equation}
\begin{alignedat}{2}
    &\maximize_{(\flam, \fbet, \tau,c) \geq \pmb{0}}\quad  &&{-\frac{1}{2}L\tau} 
    \\
     &\mbox{subject to} &&\pmb{0} = \tau\pf_{0}  - c\pf_N + \sum_{k=1}^{N} \lambda_k (\pmb{f}_{k} - \pf_{k-1}) +\sum_{k=0}^{N-1} \beta_k (\pmb{f}_k- \pf_N)
     \\
     & &&0 \preceq S(\flam, \fbet, c),
\end{alignedat} \label{eqn::OBL-G-dualPEP}    
\end{equation}
where $S$ is defined as 
\begin{align*}
    S(\flam, \fbet,  c) &=  - \pg_N \pg_N^\intercal +  \frac{c}{2L} \pg_N \pg_N^\intercal  + \sum_{k=1}^{N} \frac{\lambda_k}{2} \left(\pmb{g}_{k} (\pmb{x}_{k-1} - \pmb{x}_{k})^\intercal +   (\pmb{x}_{k-1} - \pmb{x}_{k})\pmb{g}_{k}^\intercal + \frac{1}{L} (\pmb{g}_{k-1} - \pmb{g}_{k})(\pmb{g}_{k-1} - \pmb{g}_{k})^\intercal\right)
    \\
    & + \sum_{k=0}^{N}\frac{\beta_k}{2}\left(\pg_k \px_N^\intercal + \px_N \pg_k^\intercal -\pmb{g}_k \pmb{x}_k^\intercal-\pmb{x}_k \pmb{g}_k^\intercal  \right).
\end{align*}

We have a strong duality argument 
\begin{align*}
    \argmin_{h_{i,j}}\maximize_{\fG, \fF_0} \mkern7mu \pg_N^\intercal \fG \pg_N  = \argmin_{h_{i,j}}\minimize_{(\flam, \fbet, \tau, c) \geq \pmb{0}}\mkern7mu {\frac{1}{2}L\tau},
\end{align*}
as ORC-F's optimality proof. Remind that \eqref{eqn::OBL-G-dualPEP} finds the ``best'' proof for the algorithm. Now we investigate the optimization step for algorithm. The last part is minimizing  \eqref{eqn::OBL-G-dualPEP} with stepsize, i.e. 
\begin{alignat}{3}
    &\minimize_{h_{i,j}}&&\minimize_{(\flam, \fbet, \tau, c) \geq \pmb{0}}\quad && {\tau}  \label{eqn::OBL-G-hij-wkdual}
    \\
     & &&\mbox{subject to} &&\pmb{0} = \tau\pf_{0}  - c\pf_N + \sum_{k=1}^{N} \lambda_k (\pmb{f}_{k} - \pf_{k-1}) +\sum_{k=0}^{N-1} \beta_k (\pmb{f}_k- \pf_N)  \label{eqn::OBL-G-const1}
     \\
     & && &&0 \preceq S(\flam, \fbet,  c). 
     \label{eqn::OBL-G-Smatrix}
\end{alignat}
We note that $\pf_i$ is \revision{a} standard unit vector mentioned in \eqref{eqn::egf} (not a variable), we can write \eqref{eqn::OBL-G-const1} as 
\begin{equation}
\begin{alignedat}{2}
&&\left(\begin{array}{ll}
&\tau - \lambda_1 + \beta_0 = 0    \\
    &\lambda_k - \lambda_{k+1} + \beta_k = 0, \qquad k \in \{1, \dots, N-1\}
\\
    &- c + \lambda_N - \sum_{k=0}^{N-1} \beta_k = 0 . 
\end{array}\right) 
\end{alignedat}\label{eqn::OBL-G-constlinear}
\end{equation}
We consider \eqref{eqn::OBL-G-Smatrix} with \eqref{eqn::OBL-G-constlinear} and FSFO's $h_{i,j}$. To be specific, we represent the dependency of $S$ to $(h_{i,j})$ since $\px_k$ are represented with $(h_{i,j})$. Then, we get
\begin{equation}
\begin{aligned}
    S(\flam, \fbet, c;(h_{i,j})) 
    &= - \pmb{g}_N \pmb{g}_N^\intercal + \frac{c}{2L}\pmb{g}_N \pmb{g}_N^\intercal
    \\
    &+ \sum_{k=1}^{N} \frac{\lambda_k}{2}\left(  \pmb{g}_{k}\left(\sum_{t=0}^{k-1} \frac{h_{k, t}}{L} \pmb{g}_t\right)^\intercal+ \left(\sum_{t=0}^{k-1} \frac{h_{k, t}}{L} \pmb{g}_t\right)\pmb{g}_{k}^\intercal + \frac{1}{L}  (\pmb{g}_{k-1} - \pmb{g}_{k})(\pmb{g}_{k-1} - \pmb{g}_{k})^\intercal
    \right)
    \\
    & + \sum_{k=0}^{N-1}\frac{ \beta_k}{2} \left( -\left(\sum_{j=k+1}^N \sum_{t=0}^{j-1} \frac{h_{j,t}}{L} \pmb{g}_t\right)\pmb{g}_{k}^\intercal -\pmb{g}_{k}\left(\sum_{j=k+1}^N \sum_{t=0}^{j-1} \frac{h_{j,t}}{L} \pmb{g}_t\right)^\intercal \right).
\end{aligned} \label{eqn::OBL-G-final-matrix}
\end{equation}

\begin{claim}
\label{claim::10}
    There is a point that makes \eqref{eqn::OBL-G-final-matrix}'s right-hand side as zero. 
\end{claim}
\begin{proof}
By calculating $2LS(\flam, \fbet,  c;(h_{i,j}))$'s $\pg_i \pg_j^\intercal$ coefficients, we have
\begin{equation}
\begin{aligned}
    &\lambda_1 - 2\beta_0 \sum_{l=1}^N h_{l,0} \qquad 
    \\
    &\lambda_{i} + \lambda_{i+1} - 2\beta_i \sum_{l=i+1}^N h_{l,i}, \qquad i \in \{1, 2, \dots, N-1\}
    \\
    &\lambda_{N} + c - 2L
    \\
    &\lambda_i (h_{i,i-1} -1 ) - \beta_i \sum_{l=i+1}^N h_{l,i-1}- \beta_{i-1} \sum_{l=i+1}^N h_{l,i}, \qquad i \in \{1,2, \dots, N-1\}
    \\
    &\lambda_{N}(h_{N,N-1} -1 )
    \\
    &\lambda_i h_{i,j} - \beta_i  \sum_{l=i+1}^N h_{l,j} - \beta_j  \sum_{l=i+1}^N h_{l,i}, i \in \{2, 3, \dots, N-1\}, \qquad j\in \{0,1, \dots, i-2\}
    \\
    &\lambda_N h_{N,j}, \qquad  j \in \{0,1,\dots, N-2\}.
\end{aligned}\label{eqn::OBL-G-Sij}
\end{equation}
For finding a solution of $S = 0$, we will set $\lambda_i$ as  $$\lambda_{N-k+1} = \frac{1}{k(k+1)} \lambda.$$ \eqref{eqn::OBL-G-constlinear} can be written as 
\begin{align*}
    &\tau - \frac{1}{N(N+1)}\lambda + \beta_0 = 0
    \\
    &\frac{1}{(N-k+1)(N-k+2)}\lambda - \frac{1}{(N-k)(N-k+1)}\lambda + \beta_k = 0, \qquad k \in \{1, \dots, N-1\}, 
    \\
    &- c + \frac{1}{2}\lambda - \sum_{k=0}^{N-1} \beta_k = 0 . 
\end{align*}
Therefore, we have
\begin{equation}
    \begin{aligned}
        &\beta_k = \frac{2}{(N-k)(N-k+1)(N-k+2)}\lambda, \qquad k  \in \{1,2, \dots, N-1\}
        \\
        &\beta_0 = \frac{1}{N(N+1)}\lambda - \tau = \hat{\beta_0} \lambda
        \\
        &c = \tau.
    \end{aligned}\label{eqn::OBL-G-constlinear3}
\end{equation}
\newpage 

Then, with \eqref{eqn::OBL-G-Sij} and \eqref{eqn::OBL-G-constlinear3}, we get
\begingroup
\allowdisplaybreaks
\begin{equation}
\begin{aligned}
    &0 =  \frac{1}{N(N+1)}\lambda - 2\hat{\beta_0}\lambda \sum_{l=1}^N h_{l,0} 
    \\
    &0 = \frac{2}{(N-i)(N-i+2)}\lambda  -  \frac{4}{(N-i)(N-i+1)(N-i+2)}\lambda \sum_{l=i+1}^N h_{l,i}, \qquad i \in \{1, 2, \dots, N-1\}
    \\
    &0 = \frac{1}{2}\lambda + \tau - 2L 
    \\
    &0 =  \frac{1}{(N-i+1)(N-i+2)}\lambda (h_{i,i-1} -1 ) - \frac{2}{(N-i)(N-i+1)(N-i+2)}\lambda  \sum_{l=i+1}^N h_{l,i-1}
    \\ &\qquad \qquad \qquad - \frac{2}{(N-i+1)(N-i+2)(N-i+3)}\lambda  \sum_{l=i+1}^N h_{l,i}, \qquad i \in \{2, \dots, N-1\}
    \\
    &0 =  \frac{1}{N(N+1)}\lambda (h_{1, 0} -1 ) - \frac{2}{(N-1)N(N+1)}\lambda  \sum_{l=2}^N h_{l,0} - \hat{\beta_0}\lambda \sum_{l=2}^N h_{l,1}
    \\
    &0 = \frac{1}{2}\lambda(h_{N,N-1} -1 ) \qquad 
    \\
    &0 = \frac{1}{(N-i+1)(N-i+2)}\lambda h_{i,j} - \frac{2}{(N-i)(N-i+1)(N-i+2)} \lambda\sum_{l=i+1}^N h_{l,j} 
    \\ &\qquad \qquad \qquad - \frac{2}{(N-j)(N-j+1)(N-j+2)}\lambda \sum_{l=i+1}^N h_{l,i}, \qquad  i \in \{2, 3, \dots, N-1\}, \qquad j\in \{1, \dots, i-2\}
    \\
    &0 = \frac{1}{(N-i+1)(N-i+2)}\lambda h_{i,0} - \frac{2}{(N-i)(N-i+1)(N-i+2)} \lambda\sum_{l=i+1}^N h_{l,0} 
    \\ &\qquad \qquad \qquad - \hat{\beta_0}\lambda\sum_{l=i+1}^N h_{l,i}, \qquad  i \in \{2, 3, \dots, N-1\}
    \\
    &0 = \frac{1}{2}\lambda h_{N,j}.
\end{aligned}\label{eqn::OBL-G-Sij'}
\end{equation}
\endgroup
Last equation indicates $h_{N,j} = 0$ for all $j = 0,1, \dots, N-2$, and by fifth equation, $h_{N, N-1} = 1$ also holds. By first and second equation, we have 
\begin{equation}
\begin{aligned}
 \sum_{l=1}^N h_{l,0} &= \frac{1}{2\hat{\beta_0}N(N+1)}
 \\
 \sum_{l=i+1}^N h_{l,i} &= \frac{N-i+1}{2},\qquad i \in \{1, 2, \dots, N-1\}.
\end{aligned}\label{eqn::OBL-G-sumhles}
\end{equation}
Since, \eqref{eqn::OBL-G-Sij'}'s forth equation is equivalent to 
\begin{equation}
    \begin{aligned}
        &0 =  (h_{i,i-1} -1 ) - \frac{2}{N-i}\sum_{l=i+1}^N h_{l,i-1}- \frac{2}{N-i+3} \sum_{l=i+1}^N h_{l,i}, \qquad i \in \{2, \dots, N-1 \}   
        \\
         &0 =  \frac{1}{N(N+1)}\lambda (h_{1, 0} -1 ) - \frac{2}{(N-1)N(N+1)}\lambda  \sum_{l=2}^N h_{l,0} - \hat{\beta_0}\lambda \sum_{l=2}^N h_{l,1}.
    \end{aligned}\label{eqn::OBL-G-Sij-4}    
\end{equation}
Combining with \eqref{eqn::OBL-G-sumhles} and \label{eqn::OBL-G-Sij-4}, we have 
\begin{align*}
    0 =  (h_{i,i-1} -1 ) - \frac{2}{N-i}\left(\frac{N-i+2}{2} - h_{i,i-1}\right)- \frac{2}{N-i+3}\frac{N-i+1}{2}
\end{align*}
for $i = 2, 3, \dots, N-1$ and 
\begin{align*}
    0 =  (h_{1,0} -1 ) - \frac{2}{N-1} \left(\frac{1}{2\hat{\beta_0}N(N+1)}
    - h_{1,0}\right)- \hat{\beta_0}N(N+1)  \frac{N}{2}.
\end{align*}
Therefore, we have 
\begin{equation}
\begin{aligned}
    &h_{1,0} = \frac{N-1}{N+1} + \frac{1}{N(N+1)^2 \hat{\beta_0}} + \frac{N^2 (N-1)}{2} \hat{\beta_0} \\
    &h_{i, i-1} = \frac{3(N-i+1)}{N-i+3}.
\end{aligned}\label{eqn::OBL-G-hiim1}
\end{equation}
With \eqref{eqn::OBL-G-hiim1}, \eqref{eqn::OBL-G-Sij'}'s seventh equation is equivalent to 
\begin{align}
  &0 = \frac{1}{(N-i+1)(N-i+2)} h_{i,j} - \frac{2}{(N-i)(N-i+1)(N-i+2)} \sum_{l=i+1}^N h_{l,j}  - \frac{N-i+1}{(N-j)(N-j+1)(N-j+2)} \label{eqn::OBL-G-six1}
\end{align} 
for $i = 2,3, \dots ,N-1$ and $j = 1,\dots, i-2$. For $i = 1, 2, \dots N-2$, and $j = 1, \dots, i-1$, we have
\begin{align}
  &0 = \frac{1}{(N-i)(N-i+1)} h_{i+1,j} - \frac{2}{(N-i-1)(N-i)(N-i+1)} \sum_{l=i+2}^N h_{l,j}  - \frac{N-i}{(N-j)(N-j+1)(N-j+2)}.\label{eqn::OBL-G-six1'}
\end{align}
With $(N-i+2)$\eqref{eqn::OBL-G-six1} $-$ $(N-i-1)$\eqref{eqn::OBL-G-six1'}, we have 
\begin{align}
    0 = \frac{1}{N-i+1} h_{i,j} - \frac{1}{N-i} h_{i+1,j} - \frac{4(N-i) + 2}{(N-j)(N-j+1)(N-j+2)}
    \label{eqn::OBL-G-hijdiff}
\end{align}
for $i = 2, 3, \dots, N-2$ and $j = 1, \dots, i-2$, and putting $i = N-1$ in \eqref{eqn::OBL-G-six1},      
\begin{align}
  h_{N-1,j} = \frac{12}{(N-j)(N-j+1)(N-j+2)}
  \label{eqn::OBL-G-hNm1}
\end{align} 
for $j = 1, \dots ,N-2$. \eqref{eqn::OBL-G-Sij'}'s eighth equation is equivalent to 
\begin{align}
  &\hat{\beta_0} \frac{N-i+1}{2} = \frac{1}{(N-i+1)(N-i+2)} h_{i,0} - \frac{2}{(N-i)(N-i+1)(N-i+2)} \sum_{l=i+1}^N h_{l,0} \label{eqn::OBL-G-seven1}
\end{align} 
for $i = 2,3, \dots, N-1$. For $i = 1, 2, \dots, N-2$, we have 
\begin{align}
  &\hat{\beta_0}\frac{N-i}{2} = \frac{1}{(N-i)(N-i+1)} h_{i+1,0} - \frac{2}{(N-i-1)(N-i)(N-i+1)} \sum_{l=i+2}^N h_{l,0} \label{eqn::OBL-G-seven1'}
\end{align} 
With $(N-i+2)$\eqref{eqn::OBL-G-seven1} $-$ $(N-i-1)$\eqref{eqn::OBL-G-seven1'}, we have 
\begin{align*}
    \hat{\beta_0}(2(N-i) + 1) = \frac{1}{N-i+1} h_{i,0} - \frac{1}{N-i} h_{i+1,0}
\end{align*}
for $i = 2, 3, \dots, N-2$, which indicates 
\begin{align}
    h_{i,0} = \left(\frac{1}{N-1} h_{2,0} - \hat{\beta_0} (i-2)(2N - i)\right) (N-i + 1)
    \label{eqn::OBL-G-hi0}
\end{align}
for $i = 3, \dots, N-1$. Putting $i = 2$ in \eqref{eqn::OBL-G-seven1} and using \eqref{eqn::OBL-G-hi0}, we have 
\begin{align*}
    h_{2,0} &= \frac{2}{N-2} \sum_{l=3}^N h_{l,0} + \hat{\beta_0} \frac{N(N-1)^2}{2}
    \\
    &= \frac{(N-3)N}{(N-2)(N-1)} h_{2,0} - \hat{\beta_0} \frac{(N-3)N(N+1)}{2} + \hat{\beta_0} \frac{N(N-1)^2}{2}
\end{align*}
which indicates $h_{2,0} = N(N-1)(N-2)\hat{\beta_0}$. With \eqref{eqn::OBL-G-hi0}, 
\begin{align*}
    h_{i,0} = \left(N(N-2) - (i-2)(2N - i)\right) (N-i + 1)\hat{\beta_0} 
\end{align*}
for $i = 2, \dots, N-1$. Moreover, by \eqref{eqn::OBL-G-sumhles}, we have
\begin{align*}
    \sum_{i=1}^N h_{i,0} &= h_{1,0} + \sum_{i=2}^N h_{i,0}
    \\
    &= \frac{N-1}{N+1} + \frac{1}{N(N+1)^2 \hat{\beta_0}} + \frac{N^2 (N-1)}{2} \hat{\beta_0} + \frac{(N-2)(N-1)N(N+1)}{4}\hat{\beta_0}
    \\
    &=\frac{1}{2\hat{\beta_0}N(N+1)} 
\end{align*} 
which indicates $\hat{\beta_0} = \frac{2\left(\sqrt{\frac{N(N+1)}{2}} -1\right)}{(N-1)N(N+1)(N+2)}$, and all $h_{i,0}$ is determined. Using \eqref{eqn::OBL-G-hijdiff} and \eqref{eqn::OBL-G-hNm1}, we can also derive $$h_{i,j} = \frac{2(N-i)(N-i+1)(N-i+2)}{(N-j)(N-j+1)(N-j+2)}$$
for $i = 2, \dots,N$ and $j = 1,\dots, i-2$. 
In sum, 
\begin{equation*}
\begin{aligned}
H &= \left(
\begin{array}{cccccc}
    \frac{N + \sqrt{2N(N+1)}}{N+2}&  &  & &  & \\
    (N-2)(N-1)N \hat{\beta_0} & \frac{3(N-1)}{N+1} & & &\\
    (N-3)(N-2)(N-1) \hat{\beta_0} & \frac{2(N-3)(N-2)(N-1)}{(N-1)(N)(N+1)} & \frac{3(N-2)}{N} & & \\    
    \vdots & \vdots & \vdots & \ddots & &  \\
    1*2*3 \hat{\beta_0} & \frac{2*1*2*3}{(N-1)(N)(N+1)} & \frac{2*1*2*3}{(N-2)(N-1)(N)} & \dots & \frac{3*2}{4} & \\
    0 & 0& 0 & \dots &0 & \frac{3 *1}{3}\\
\end{array}
\right)
\end{aligned}
\end{equation*}
satisfies $S = 0$ with $\tau, \hat{\beta_0}, \lambda$ obtained above. 
\end{proof}
Since there cannot find the relationship similar to \citep[Theorem 3]{drori2014performance}, we cannot find the optimality of OBL-G$_\flat$. However, \eqref{eqn::OBL-G-hij-wkdual} is bi-convex over $(h_{i,j})$ and $(\flam, \fbet, \tau, c)$, so far each given $N$ we numerically solved and suspect that obtained solution in the previous claim is the solution for \eqref{eqn::OBL-G-hij-wkdual}. 

\begin{claim}
The obtained algorithm has the same $(h_{i,j})$ with OBL-G$_\flat$. 
\end{claim}
\begin{proof}
By calculating $(h_{i,j})$ of OBL-G$_\flat$, we can prove the equivalence of the obtained solution and OBL-G$_\flat$. Indeed, OBL-G$_\flat$ is obtained by using \citep{lee2021geometric}'s auxiliary sequences.
We will show that obtained $(\hat{h}_{i,j})$ satisfies \begin{align*}
    x_0 - \sum_{i=1}^{k+1} \sum_{j=0}^{i-1} \frac{\hat{h}_{i.j}}{L} \nabla f(x_j) &= \frac{N-k-2}{N-k+2} \left( x_0 - \sum_{i=1}^k \sum_{j=0}^{i-1} \frac{\hat{h}_{i,j}}{L} \nabla f(x_j) - \frac{1}{L} \nabla f(x_k)\right)
    \\&+ \frac{4}{N-k+2}\left(x_0 - \sum_{j=0}^k \frac{N-k+1}{2L} \nabla f(x_j)\right),
\end{align*}
which is re-written form of OBL-G$_\flat$. Comparing $\nabla f(x_j)$'s each coefficient, we should prove 
\begin{align*}
    \sum_{i = j+1}^{k+1} \hat{h}_{i,j} &= \frac{N-k-2}{N-k+2}\sum_{i = j+1}^{k} \hat{h}_{i,j} +  \frac{4}{N-k+2} \frac{N-j+1}{2} \qquad j \in \{0,1,\dots, k-1\}
    \\
    \hat{h}_{k+1,k} &= \frac{N-k-2}{N-k+2} +  \frac{4}{N-k+2}\frac{N-k+1}{2},
\end{align*}
which can be easily checked with the matrix $H$. 
\end{proof}

\end{document}